\documentclass[twocolumn]{svjour3}

\usepackage[latin1]{inputenc}
\usepackage{times}

\usepackage{psfrag, amsmath, amssymb, amscd, amsfonts, latexsym, epsf, graphicx}
\usepackage{natbib}
\usepackage{hyperref}

\usepackage[labelformat=simple]{subcaption}

\def\bbbr{{\mathbb R}} 
\def\bbbz{{\mathbb Z}}

\def\sampl{\scriptsize\mbox{sampl}}
\def\normsampl{\scriptsize\mbox{normsampl}}
\def\intdisc{\scriptsize\mbox{int}}
\def\disc{\scriptsize\mbox{disc}}

\def\norm{\scriptsize\mbox{norm}}

\def\scaleestrel{\scriptsize\mbox{scaleest,rel}}
\def\scaleref{\scriptsize\mbox{ref}}
\def\blob{\scriptsize\mbox{blob}}
\def\edge{\scriptsize\mbox{edge}}
\def\ridge{\scriptsize\mbox{ridge}}

\def\hybrnormsampl{\scriptsize\mbox{hybr-sampl}}
\def\hybrint{\scriptsize\mbox{hybr-int}}

\journalname{arXiv preprint}

\begin{document}

\title{\bf Approximation properties relative to continuous scale space for hybrid discretisations of Gaussian derivative operators%
\thanks{The support from the Swedish Research Council 
              (contract 2022-02969) is gratefully acknowledged. }}

\titlerunning{Approximation properties for hybrid discretisations of Gaussian derivative operators}

\author{Tony Lindeberg}

\institute{Tony Lindeberg \at
              Computational Brain Science Lab,
              Division of Computational Science and Technology,
              KTH Royal Institute of Technology,
              SE-100 44 Stockholm, Sweden.
              \email{tony@kth.se},
              ORCID: 0000-0002-9081-2170.}

\date{Received: date / Accepted: date}

\maketitle

\begin{abstract}
This paper presents an analysis of properties of two hybrid discretisation methods
for Gaussian derivatives, based on convolutions with either the normalised sampled Gaussian kernel or the integrated Gaussian kernel followed by central differences.
The motivation for studying these discretisation methods is that in
situations when multiple spatial derivatives of different orders are
needed at the same scale level, they can be computed significantly
more efficiently, compared to more direct derivative approximations
based on explicit convolutions with either sampled Gaussian derivative
kernels or integrated Gaussian derivative kernels.
  
%While these computational benefits do also hold for the genuinely discrete approach for computing discrete analogues of Gaussian derivatives, based on convolution with the discrete analogue of the Gaussian kernel followed by central differences, the underlying mathematical primitives for the discrete analogue of the Gaussian kernel, in terms of modified Bessel functions of integer order, may not be available in certain frameworks for image processing, such as when performing deep learning based on scale-parameterised filters in terms of Gaussian derivatives, with learning of the scale levels. The hybrid discretisations studied in this paper do, from this perspective, offer a computationally more efficient way of implementing deep networks based on Gaussian derivatives for such use cases.
  
We characterise the properties of these hybrid discretisation methods in terms of quantitative performance measures, concerning the amount of spatial smoothing that they imply, as well as the relative consistency of scale estimates obtained from scale-invariant feature detectors with automatic scale selection, with an emphasis on the behaviour for very small values of the scale parameter, which may differ significantly from corresponding results obtained from the fully continuous scale-space theory, as well as between different types of discretisation methods.

The presented results are intended as a guide, when designing as well as interpreting the experimental results of scale-space algorithms that operate at very fine scale levels.

\end{abstract}

\begin{figure*}[hbtp]

  \begin{center}
    \begin{tabular}{ccc}
      {\em Continuous Gauss\/} & {\em  Hybrid normalized sampled Gauss\/}
      & {\em Hybrid integrated Gauss\/} \\
      \includegraphics[width=0.22\textwidth]{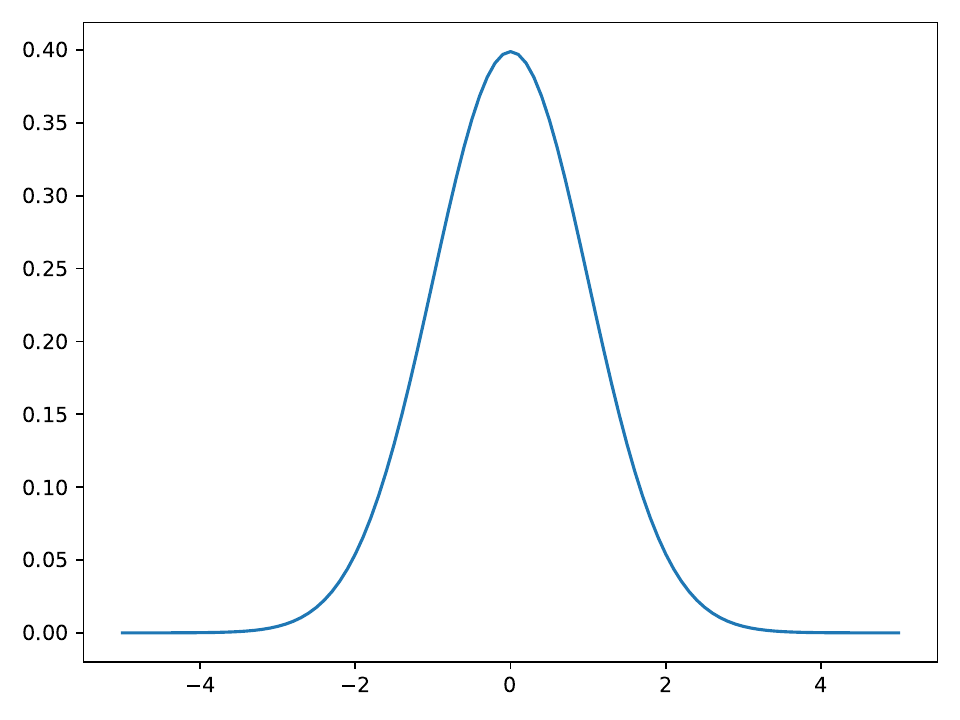} 
      & \includegraphics[width=0.22\textwidth]{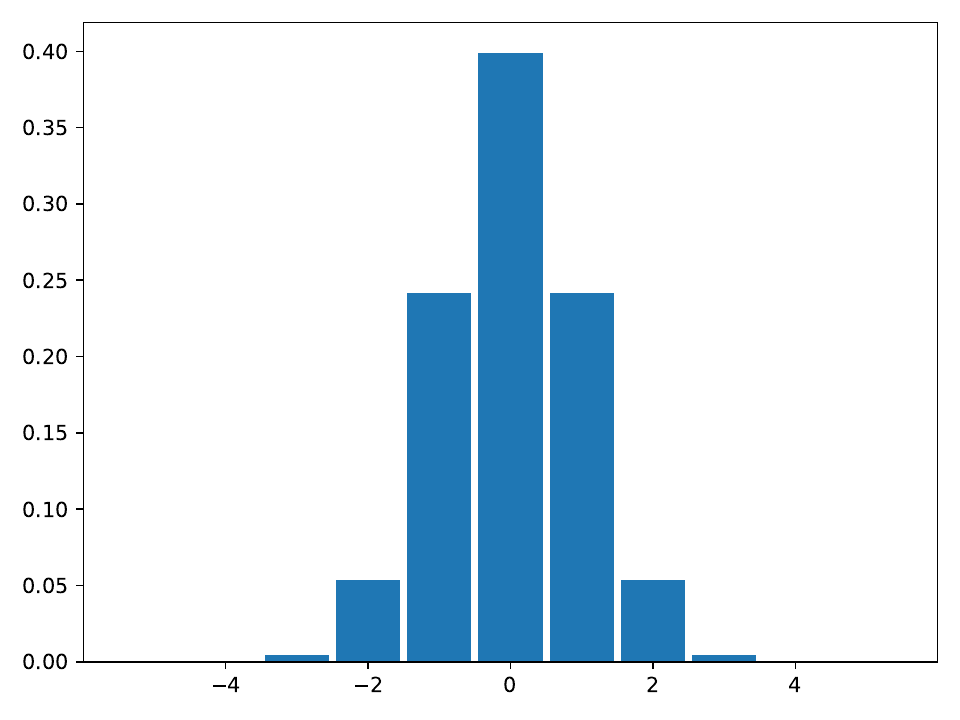}
      & \includegraphics[width=0.22\textwidth]{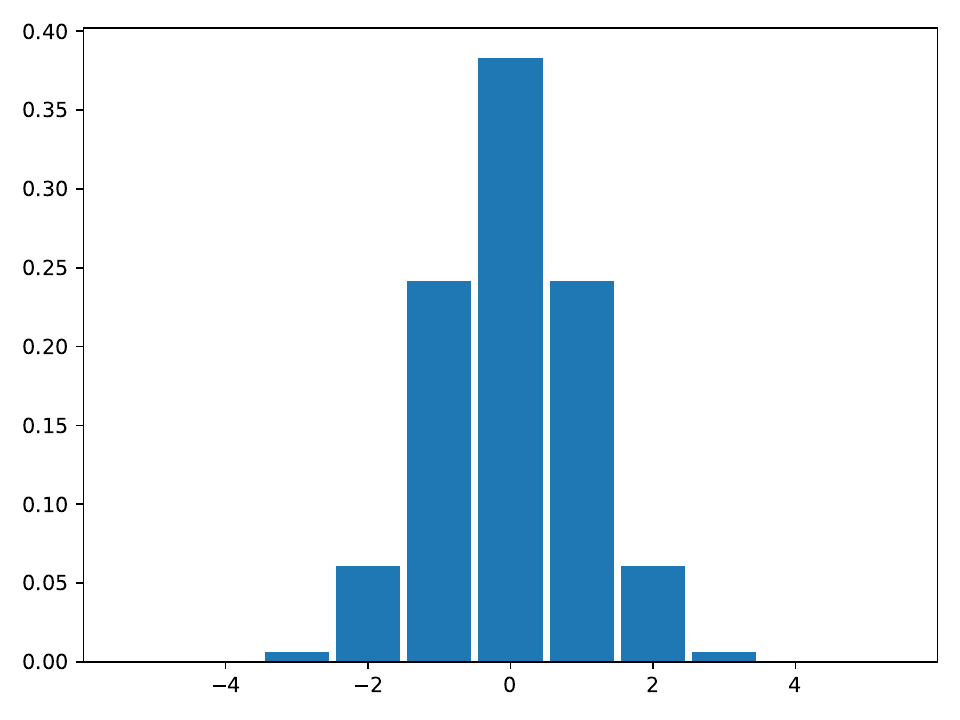}
\\
      \includegraphics[width=0.22\textwidth]{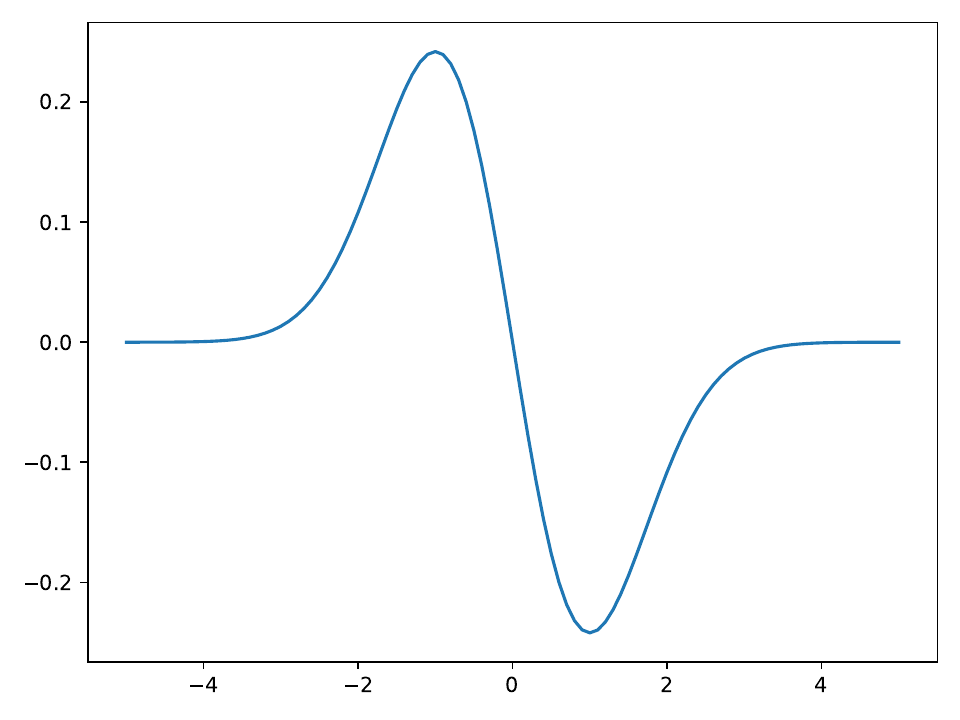}
      & \includegraphics[width=0.22\textwidth]{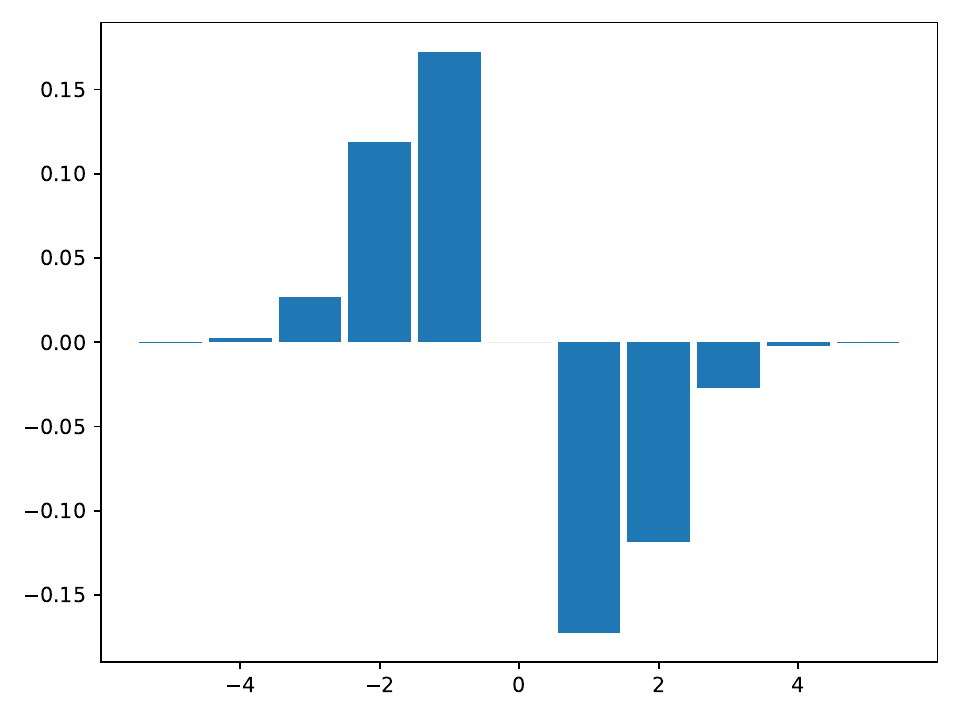}
      & \includegraphics[width=0.22\textwidth]{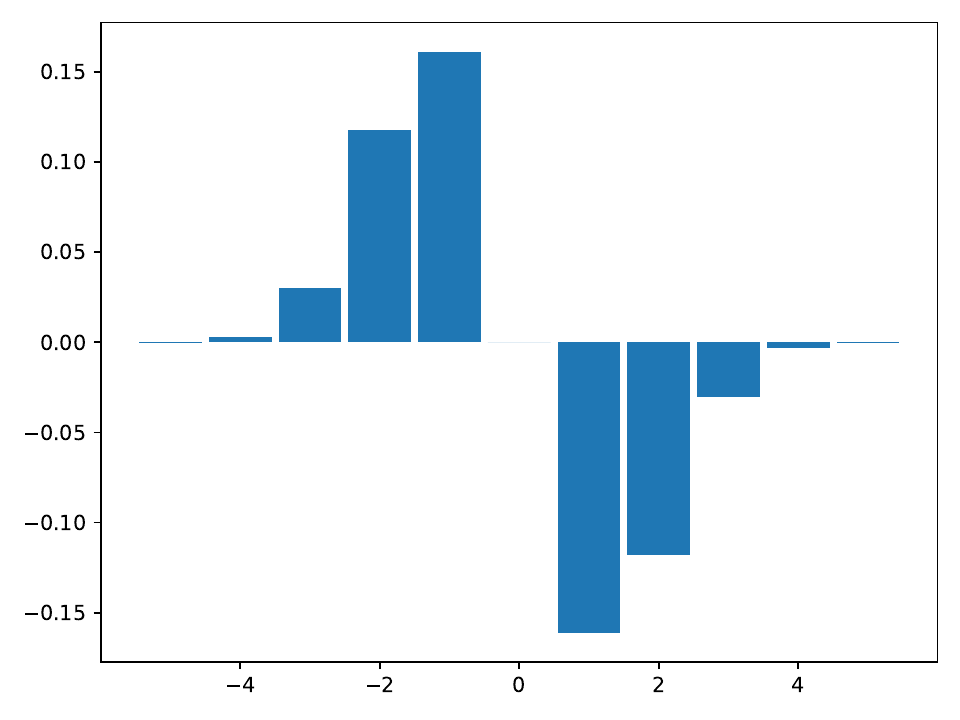}
\\
      \includegraphics[width=0.22\textwidth]{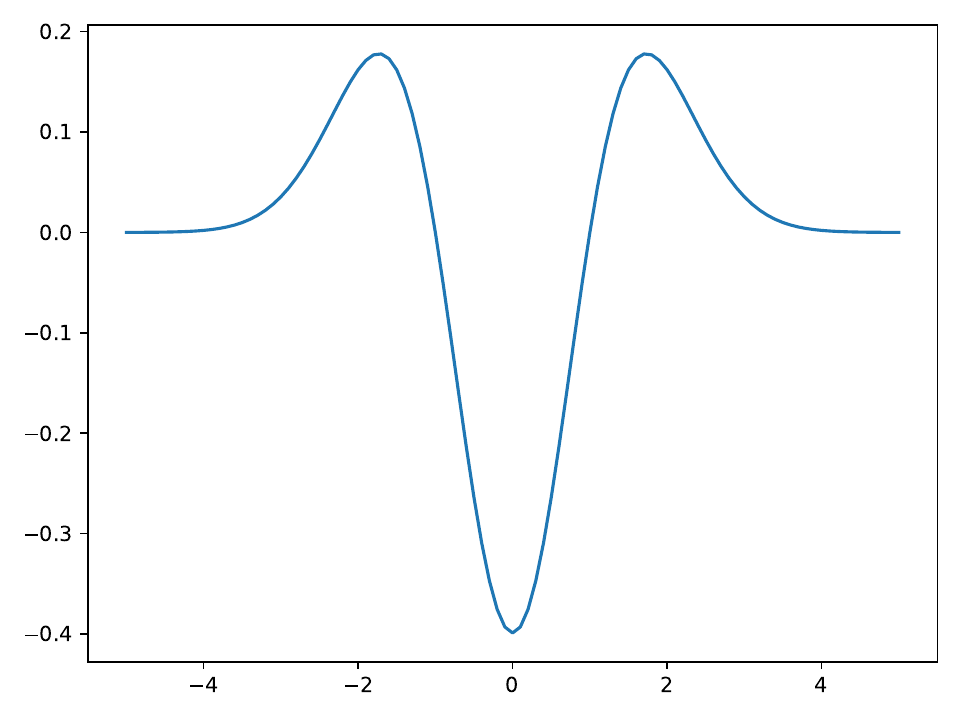}
      & \includegraphics[width=0.22\textwidth]{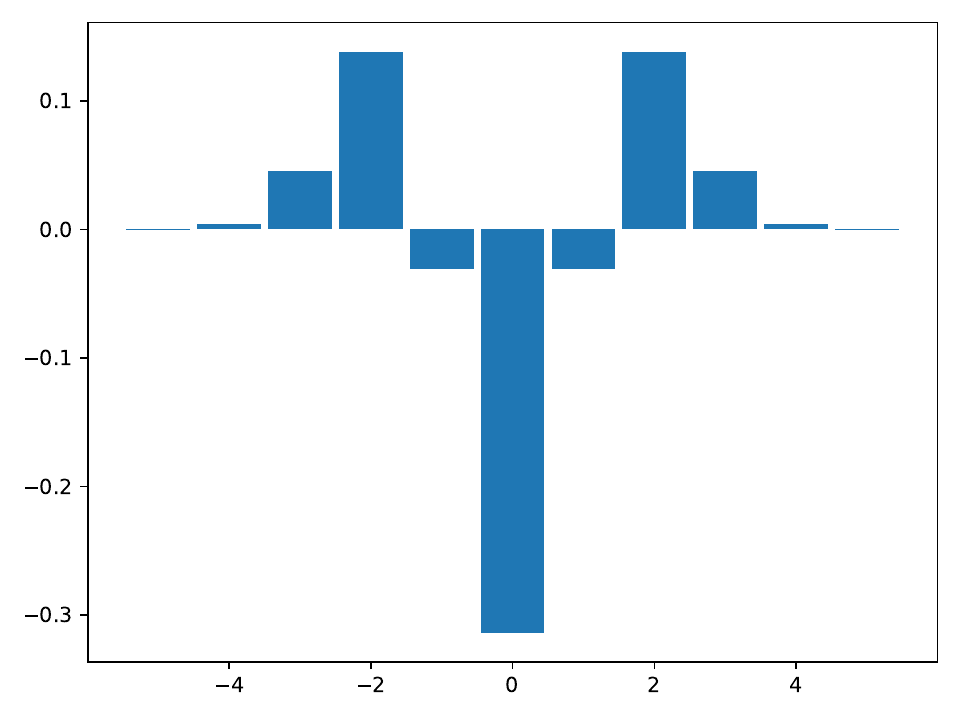}
      & \includegraphics[width=0.22\textwidth]{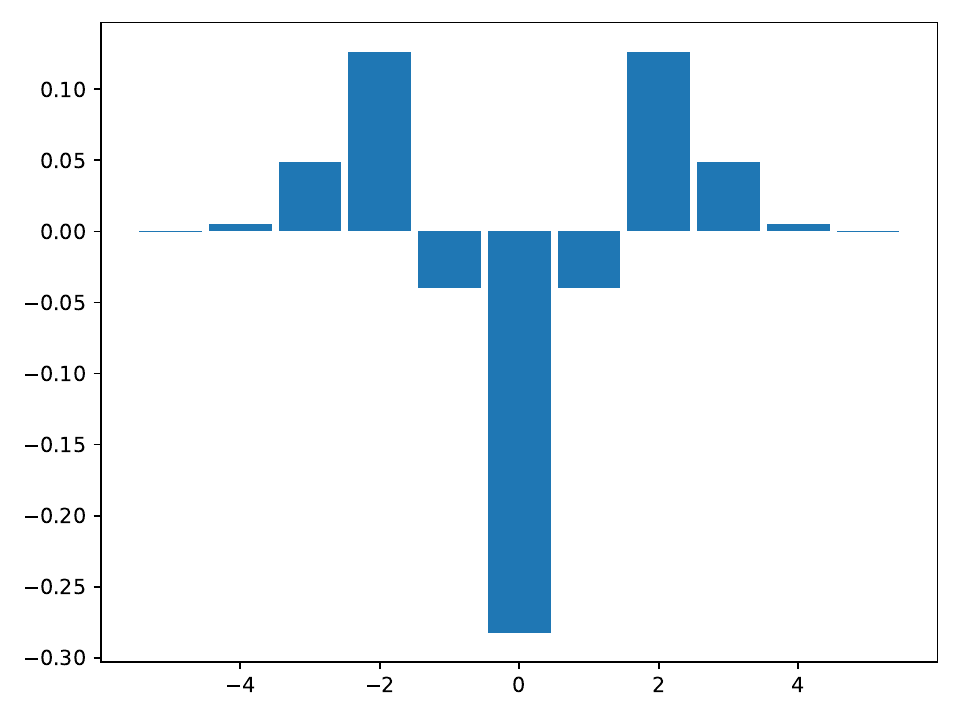}
\\
      \includegraphics[width=0.22\textwidth]{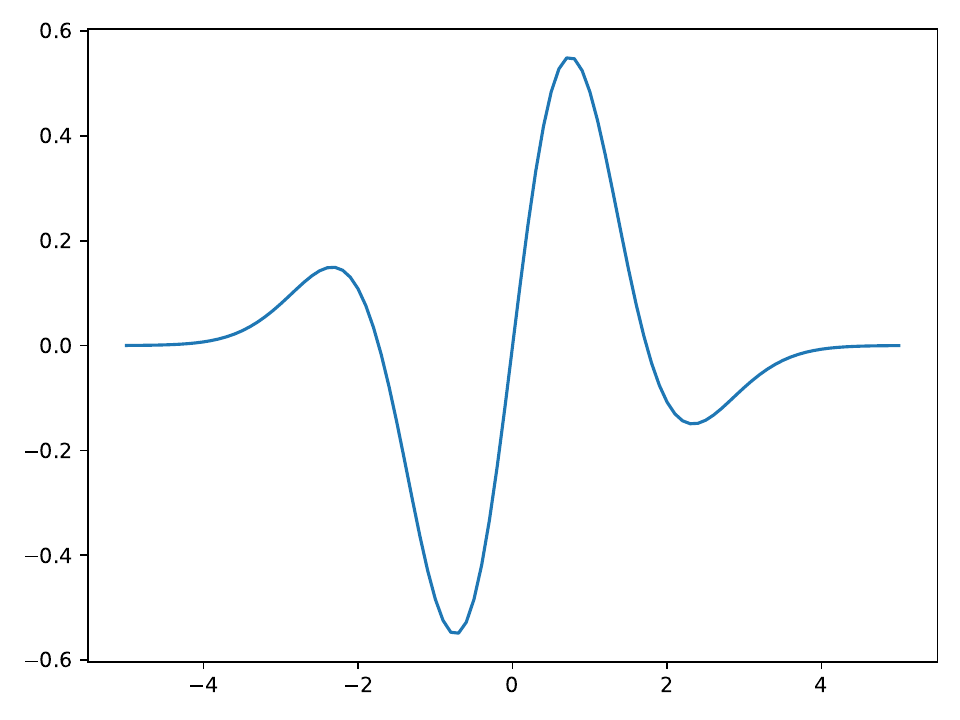}
      & \includegraphics[width=0.22\textwidth]{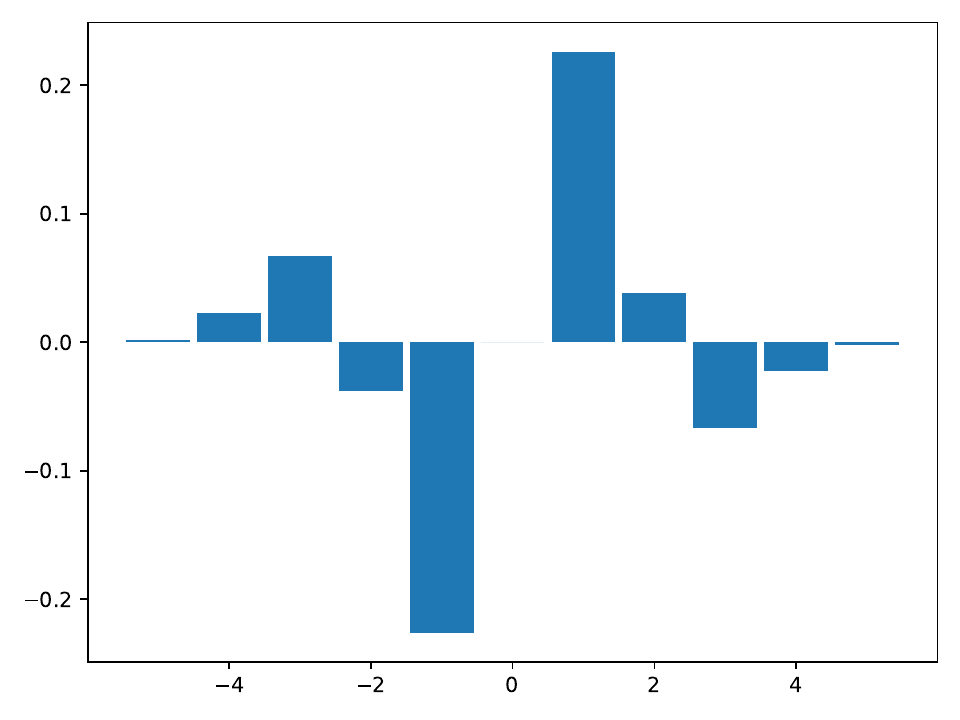}
      & \includegraphics[width=0.22\textwidth]{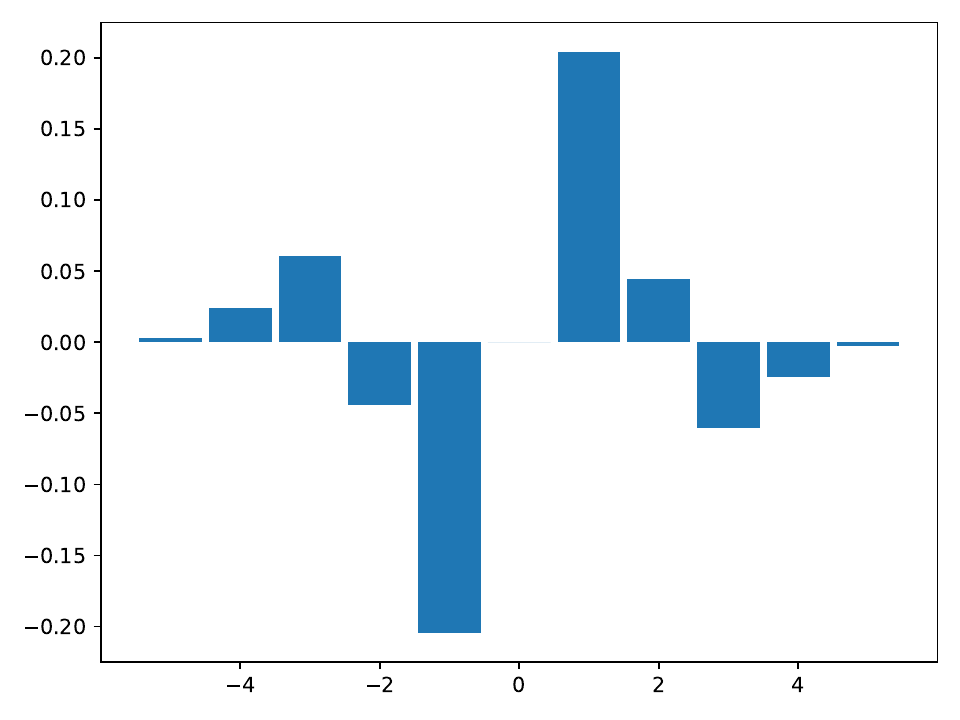}
\\
      \includegraphics[width=0.22\textwidth]{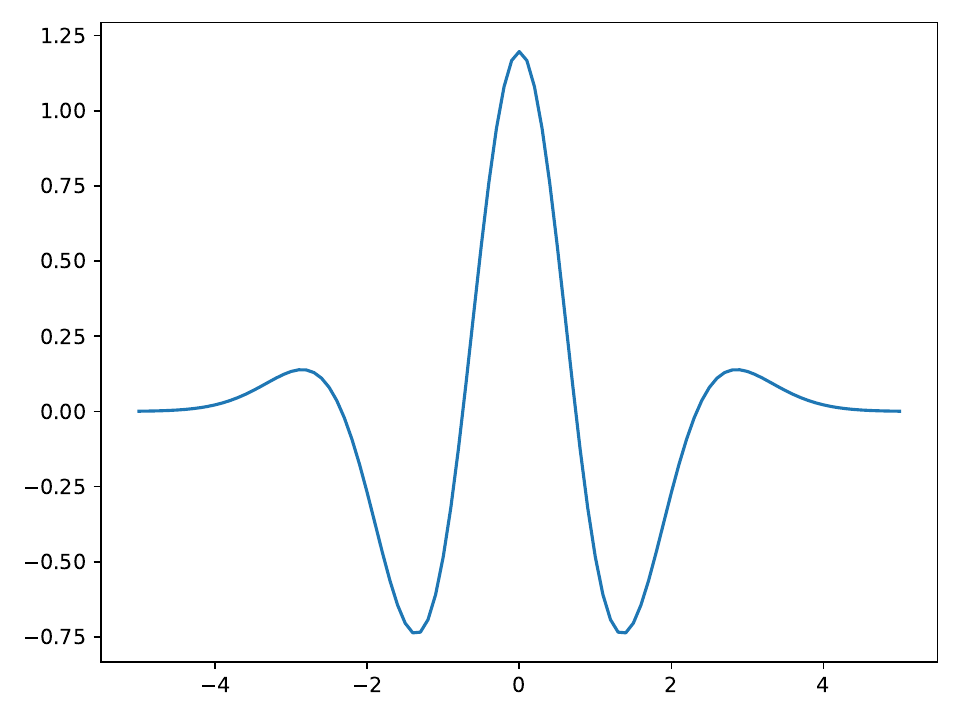}
      & \includegraphics[width=0.22\textwidth]{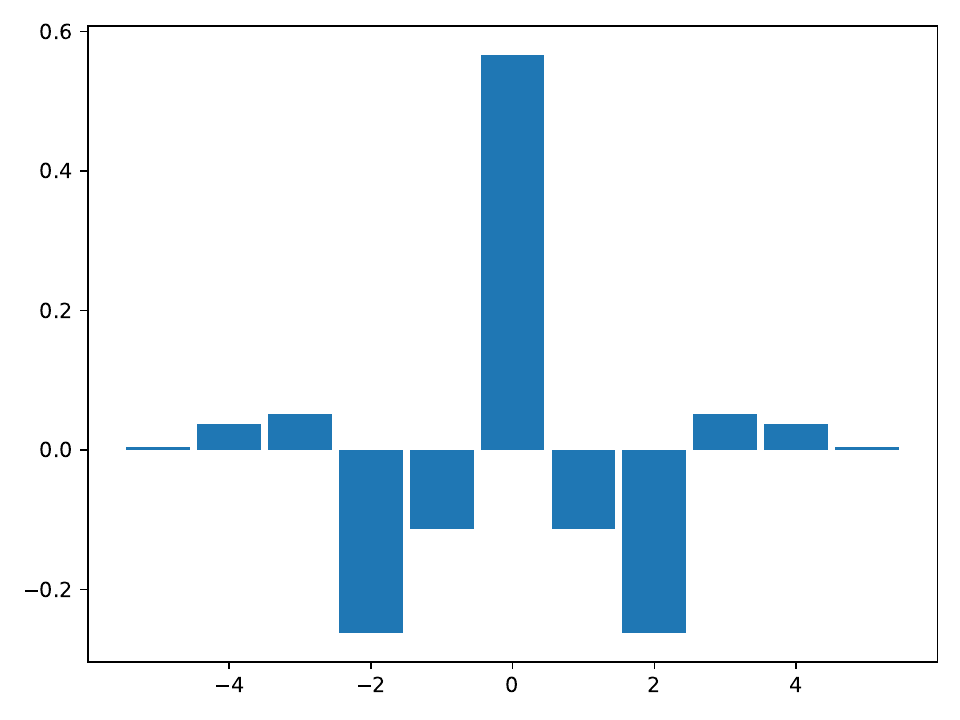}
      & \includegraphics[width=0.22\textwidth]{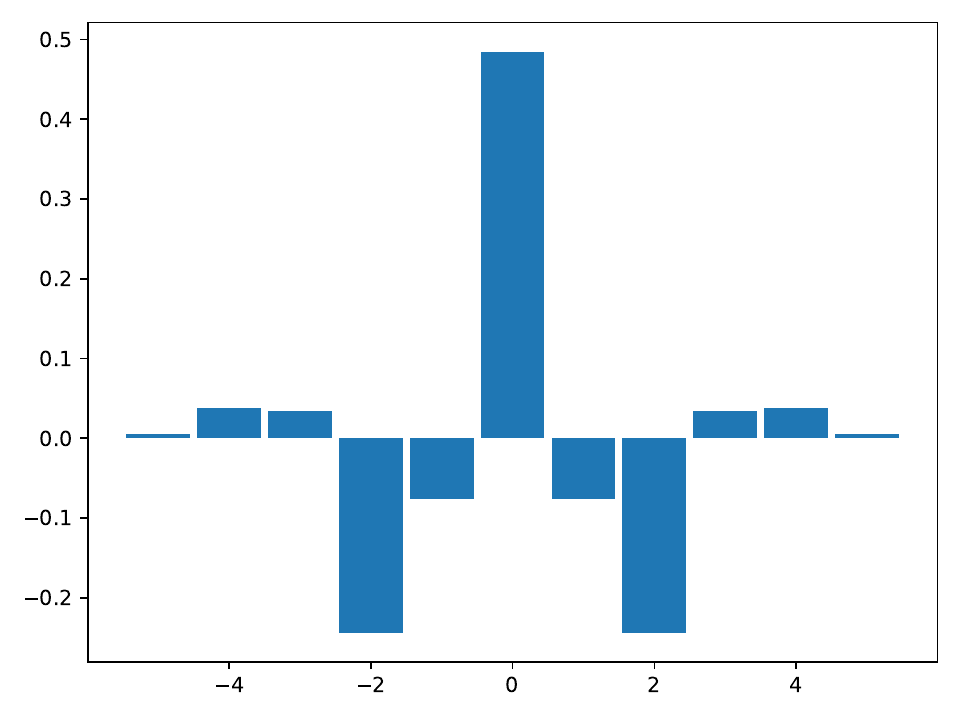}
 \\
    \end{tabular}
  \end{center}

  \caption{Graphs of the main types of Gaussian smoothing kernels as well as of the
    equivalent convolution kernels for the 
    hybrid discretisations of Gaussian derivative operators
    considered specially in this paper, here at the
    scale $\sigma = 1$, with the raw
    smoothing kernels in the top row and the order of spatial differentiation
    increasing downwards up to order 4:
    {\bf (left column)} continuous Gaussian kernel and continuous Gaussian derivatives,
    {\bf (middle column)} normalised sampled Gaussian kernel and central differences
    applied to the normalised sampled Gaussian kernel,
    {\bf (right column)} integrated Gaussian kernel and central differences
    applied to the  integrated Gaussian kernel.
    Note that the scaling of the vertical axis may vary between the
    different subfigures.
    ({\bf Horizontal axes:} the 1-D spatial coordinate $x \in [-5, 5]$.)
    (Graphs of the regular sampled Gaussian derivative kernels, the regular integrated Gaussian derivative kernels and the discrete analogues of Gaussian derivatives up to order 4 are shown in Figure~1 in (Lindeberg \citeyear{Lin24-JMIV}).)}
  \label{fig-kernel-graphs}
\end{figure*}

\section{Introduction}
\label{sec-intro}

When analysing image data by automated methods,
a fundamental constraints originates from
the fact that natural images may contain different types of structures
at different spatial scales. For this reason, the notion of
% multi-scale representations such as
scale-space representation
(Iijima \citeyear{Iij62};  Witkin \citeyear{Wit83}; Koenderink \citeyear{Koe84};
Koenderink and van Doorn \citeyear{KoeDoo87-BC}, \citeyear{KoeDoo92-PAMI};
 Lindeberg \citeyear{Lin93-Dis}, \citeyear{Lin94-SI}, \citeyear{Lin10-JMIV};
 Florack \citeyear{Flo97-book};
 Weickert {\em et al.\/}\ \citeyear{WeiIshImi99-JMIV};
 ter Haar Romeny \citeyear{Haa04-book})
%and
%pyramid representation (Burt and Adelson \citeyear{BA83-COM},
%Crowley and Stern \citeyear{Cro84-dolp},
%Simoncelli {\em et al.\/} \citeyear{SimFreAdeHee92-IT},  
%Simoncelli and Freeman \citeyear{SimFre-ICIP95},
%Lindeberg and Bretzner \citeyear{LinBre03-ScSp},
%Crowley and Riff \citeyear{CroRif03-ScSp},
%Lowe \citeyear{Low04-IJCV})
has been developed to process the image
data at multiple scales, in such a way that different types of image
features can be obtained depending on the spatial extent of the
image operators.
Specifically, according to both theoretical and empirical findings
in the area of scale-space theory,
% In classical computer vision,
Gaussian derivative responses, or 
approximations thereof, can be used as a powerful basis for %specifically
expressing a rich variety of feature detectors, in terms of provably scale-covariant or
scale-invariant image operations, that can in an automated manner
handle variabilities in scale, caused by varying the distance between
the observed objects and the camera
(Lindeberg \citeyear{Lin97-IJCV}, \citeyear{Lin97-IJCV},
\citeyear{Lin13-ImPhys}, \citeyear{Lin12-JMIV}, \citeyear{Lin15-JMIV},
\citeyear{Lin21-EncCompVis};
Bretzner and Lindeberg \citeyear{BL97-CVIU};
Chomat {\em et al.\/} \citeyear{ChoVerHalCro00-ECCV};
Lowe \citeyear{Low04-IJCV};
Bay {\em et al.\/} \citeyear{BayEssTuyGoo08-CVIU}).
More recently, Gaussian derivative operators have also been used as
as a basis, to formulate
parameterised mathematical primitives, to be used as the layers in deep networks
(Jacobsen {\em et al.\/} \citeyear{JacGemLouSme16-CVPR};
Worrall and Welling \citeyear{WorWel19-NeuroIPS};
Lindeberg \citeyear{Lin20-JMIV}, \citeyear{Lin22-JMIV};
Pintea {\em et al.\/} \citeyear{PinTomGoeLooGem21-IP};
Sangalli {\em et al.\/} \citeyear{SanBluVelAng22-BMVC};
Penaud-Polge {\em et al.\/} \citeyear{PenVelAng22-ICIP};
Yang {\em et al.\/} \citeyear{YanDasMah23-arXiv};
Gavilima-Pilataxi and Ibarra-Fiallo \citeyear{GavIva23-ICPRS}).

When to implement the underlying Gaussian derivative operators in scale-space
theory in practice, special attention does, however, need to be taken
concerning the fact that most of the scale-space formulations are
based on continuous signals or images (see Appendix~\ref{app-theor-scsp} for a conceptual background),
while real-world signals and images are discrete
(Lindeberg \citeyear{Lin90-PAMI}, \citeyear{Lin93-Dis}, \citeyear{Lin24-JMIV}).
%except for the genuine discrete scale-space theories in (REFS).
Thus, when discretising the composed effect of the underlying
Gaussian smoothing operation and the following derivative computations,
it is essential to ensure that the
desirable properties of the theoretically well-founded scale-space representations
are to a sufficiently good degree of approximation transferred to the discrete implementation. Simultaneously, the amount of necessary computations needed
for the implementation may often constitute a limiting factor, when to choose
an appropriate discretisation method for expressing the actual algorithms,
that are to operate on the discrete data to be analysed.

While one may argue that at sufficiently coarse scale levels, it ought to be the case that the choice of discretisation method should not significantly affect the quality of the output of a scale-space algorithm, at very fine scale levels, on the other hand,
the properties of a discretised implementation of notions from scale-space theory may depend strongly on the actual choice of a discretisation method.

The subject of this article, is to perform a more detailed analysis of a class of hybrid discretisation methods,
based on convolution with either the normalised sampled Gaussian kernel or the integrated Gaussian kernel, followed by
computations of discrete derivative approximations by central difference operators,
and specifically characterise the degree of approximation of continuous expressions
in scale-space theory, that these discretisation give rise to,
see Figure~\ref{fig-kernel-graphs} for examples of graphs of equivalent convolution
kernels corresponding to these discretisation methods.
This class of discretisation methods was outlined among extensions to future work
in Section~6.1 in
(Lindeberg \citeyear{Lin24-JMIV}), and was also complemented
with a description about their theoretical 
properties in Footnote~13 in
(Lindeberg \citeyear{Lin24-JMIV}).
There were, however, no further in-depth characterisations of the approximation properties of these discretisations, with regard to what results they lead to
in relation to corresponding results from the continuous scale-space theory.

The main goal of this paper is to address this topic in terms of a set of quantitative performance measures, intended to be of general applicability for different types of visual tasks. Specifically, we will perform comparisons to the
other main types of discretisation methods considered in
(Lindeberg \citeyear{Lin24-JMIV}),
based on either (i)~explicit convolutions with sampled Gaussian derivative kernels,
(ii)~explicit convolutions with integrated Gaussian derivative kernels, or
(iii)~convolution with the discrete analogue of the Gaussian kernel, followed by computations of discrete derivative approximations by central difference operators.

A main rationale for studying this class of hybrid discretisations
is that in situations when multiple Gaussian derivative responses of different orders are needed at the same scale level,
these hybrid discretisations imply substantially lower amounts of computations,
compared to explicit convolutions with either sampled Gaussian derivative kernels or integrated Gaussian derivative kernels for each order of differentiation.
The reason for this better computational efficiency, which also holds for
the discretisation approach based on convolution with the discrete analogue
of the Gaussian kernel followed by central differences, is that the spatial smoothing part of the operation, which is performed over a substantially larger number of input data than the small-support central difference operators, can be shared between
the different orders of differentiation.

A further rationale for studying these hybrid discretisations is that in certain applications, such as the use of Gaussian derivative operators in deep learning architectures (Jacobsen {\em et al.\/} \citeyear{JacGemLouSme16-CVPR},
Lindeberg \citeyear{Lin21-SSVM,Lin22-JMIV},
Pintea {\em et al.\/} \citeyear{PinTomGoeLooGem21-IP},
Sangalli {\em et al.\/} \citeyear{SanBluVelAng22-BMVC},
Penaud-Polge {\em et al.\/} \citeyear{PenVelAng22-ICIP},
Gavilima-Pilataxi and Ibarra-Fiallo \citeyear{GavIva23-ICPRS}),
the modified Bessel functions of integer order, as used as the
underlying mathematical primitives in the discrete analogue of the
Gaussian kernel, may, however, not be fully available in the framework
used for implementing the image processing operations. For this
reason, the hybrid discretisations may, for efficiency reasons,
constitute an interesting alternative to using discretisations in
terms of either sampled Gaussian derivative kernels or integrated
Gaussian derivative kernels, when to implement certain tasks, such as
learning of the scale levels by backpropagation, which usually require
full availability of the underlying mathematical primitives in the scale-parameter\-ised filter family with regard to the deep learning framework, to be able to propagate the gradients between the layers in the deep learning architecture.

Deliberately, the scope of this paper is therefore to complement 
the in-depth treatment of discretisations of the Gaussian smoothing
operation and the Gaussian derivative operators, and as a specific
complement to the outline of the hybrid discretisations in the future
work in (Lindeberg \citeyear{Lin24-JMIV}).

We will therefore not consider other theoretically well-founded discretisations
of scale-space operations (Wang \citeyear{Wan00-TIP},
Lim and Stiehl \citeyear{LimSti03-ScSp},
Tschirsich and Kuijper \citeyear{TscKui15-JMIV},
Slav{\'\i}k and Stehl{\'\i}k \citeyear{SlaSte15-JMathAnalAppl},
Rey-Otero and Delbracio \citeyear{OteDel16-IPOL}).
Nor will we consider alternative approaches in terms
of pyramid representations (Burt and Adelson \citeyear{BA83-COM},
Crowley and Stern \citeyear{Cro84-dolp},
Simoncelli {\em et al.\/} \citeyear{SimFreAdeHee92-IT},  
Simoncelli and Freeman \citeyear{SimFre-ICIP95},
Lindeberg and Bretzner \citeyear{LinBre03-ScSp},
Crowley and Riff \citeyear{CroRif03-ScSp},
Lowe \citeyear{Low04-IJCV}),
Fourier-based implementations,
splines (Unser {\em et al.\/} \citeyear{UnsAldEde91-PAMI,UnsAldEde93-SP},
Wang and Lee \citeyear{WanLee98-PAMI},
Bouma {\em et al.\/} \citeyear{BouVilBesRomGer07-ScSp},
Bekkers \citeyear{Bek20-ICLR},
Zheng {\em et al.\/} \citeyear{ZheGonYouTao22-IJCV}),
recursive filters (Deriche \citeyear{Der92-ICIP},
Young and van~Vliet \citeyear{YouVlie95-SP},
van Vliet {\em et al.\/} \citeyear{VliYouVer98-PR},
Geusebroek {\em et al.\/} \citeyear{GeuSmeWei03-TIP},
Farneb{\"a}ck and Westin \citeyear{FarWes06-JMIV},
Charalampidis \citeyear{Cha16-SP}),
or specific wavelet theory
(Mallat \citeyear{Mal89-PAMI,Mal89-ASSP,Mal99-book,Mal16-RoySoc},
Daubechies \citeyear{Dau92-book},
Meyer \citeyear{Mey92-book},
Teolis \citeyear{Teo95-book},
Debnath and Shah \citeyear{DebSha02-book}).
%although the notion of Gaussian
%derivatives can also be considered as a family of wavelets.

Instead, we will focus on a selection of five specific methods, for implementing
Gaussian derivative operations in terms of purely discrete convolution operations, and
then with the emphasis on the behaviour for very small values of the
scale parameter. This focus on very fine scale levels is particularly
motivated from requirements regarding deep learning, where deep
learning architectures often prefer to base their decisions on very
fine-scale information, and specifically below the rule of thumb in
classical computer vision, of not going trying to go to scale levels below a standard
deviation $\sigma = \sqrt{s}$ of the Gaussian derivative kernels
below, say $1/\sqrt{2}$ or 1, in units of the grid spacing.

\section{Methods}

The notion of scale-space representation is very general, and applies
to wide classes of signals. % over different dimensions.
Based on the separability property when computing Gaussian derivative responses
in arbitrary dimensions (see Appendices~\ref{app-sep-gauss-smooth-gauss-ders}--\ref{app-def-1D-model-probl}),
%for a theoretical background),
% regarding this separability property),
%as well as the different types of spatial
%discretisations considered here, based on the in-depth treatment of
%different discretisations for the Gaussian smoothing and the Gaussian
%derivative operators in Lindeberg (\citeyear{Lin24-JMIV})),
we will in this section focus on discretisations of 1-D Gaussian derivative kernels.
By separable extension over multiple dimensions,
% the methodology described below
this methodology can then be applied to signals, images and
video over arbitrary numbers of dimensions.

\subsection{Discretisation methods for 1-D Gaussian derivative operators}
\label{sec-discr-methods}

Given the definition of a scale-space representation of a one-dimensional
continuous signal
(Iijima \citeyear{Iij62};
 Witkin \citeyear{Wit83};
 Koenderink \citeyear{Koe84};
 Koenderink and van Doorn \citeyear{KoeDoo87-BC,KoeDoo92-PAMI};
 Lindeberg \citeyear{Lin93-Dis,Lin94-SI,Lin10-JMIV};
 Florack \citeyear{Flo97-book};
 Sporring {\em et al.\/}\ \citeyear{SpoNieFloJoh96-SCSPTH};
 Weickert {\em et al.\/}\ \citeyear{WeiIshImi99-JMIV};
 ter Haar Romeny \citeyear{Haa04-book}),
 the 1-D Gaussian kernel is defined according to
 (for $x \in \bbbr$):
\begin{equation}
  g(x;\; s) = \frac{1}{\sqrt{2 \pi s}} \, e^{-x^2/2s},
\end{equation}
where the parameter $s \in \bbbr_+$ is referred to as the scale parameter,
and any 1-D Gaussian derivative kernel for spatial differentiation
order $\alpha$ is defined according to
\begin{equation}
  g_{x^{\alpha}}(x;\; s) = \partial_{x^{\alpha}} g(x;\; s),
\end{equation}
with the associated computation of Gaussian derivative responses
from any 1-D input signal $f(x)$, in turn, defined according to
\begin{equation}
  \label{eq-1D-cont-gauss-der-conv}
   L_{x^{\alpha}}(x;\; s) = \int_{u \in \bbbr} g_{x^{\alpha}}(u;\, s) \, f(x - u) \, du.
 \end{equation}
Let us first consider the following ways of approximating the Gaussian
convolution operation for discrete data, based on convolutions with
either (for $n \in \bbbz)$:
\begin{itemize}
\item
  the sampled Gaussian kernel defined according to
  (see also Appendix~\ref{app-sampl-gauss-kern})
  \begin{equation}
  \label{eq-sampl-gauss}
  T_{\sampl}(n;\; s) = g(n;\; s),
\end{equation}
\item
  the normalised sampled Gaussian kernel defined according to
  (see also Appendix~\ref{app-norm-sampl-gauss-kern})
  \begin{equation}
  \label{eq-def-norm-sampl-gauss}
  T_{\normsampl}(n;\; s) = \frac{g(n;\; s)}{\sum_{m \in \bbbz} g(m;\; s)},
\end{equation}
\item
  the integrated Gaussian kernel defined according to
  (Lindeberg \citeyear{Lin93-Dis} Equation~(3.89))
  (see also Appendix~\ref{app-int-gauss-kern})
\begin{equation}
  \label{eq-def-int-gauss-kern}
   T_{\intdisc}(n;\; s) = \int_{x = n - 1/2}^{n + 1/2} g(x;\; s) \, dx,
 \end{equation}
\item
  or the discrete analogue of the Gaussian kernel defined according to
  (Lindeberg \citeyear{Lin90-PAMI} Equation~(19))
  (see also Appendix~\ref{app-disc-gauss})
   \begin{equation}
  \label{eq-disc-gauss}
   T_{\disc}(n;\; s) = e^{-s} I_n(s),
 \end{equation}
where $I_n(s)$ denotes the modified Bessel functions of integer order
 (Abramowitz and Stegun \citeyear{AS64}).
\end{itemize}
Then, we consider the following previously studied methods for discretising the computation of Gaussian derivative operators, in terms of either:
\begin{itemize}
\item
  convolutions with sampled Gaussian derivative kernels according to
   (see also Appendix~\ref{app-sampl-gaussder-kern})
  \begin{equation}
  \label{eq-sampl-gauss-der}
  T_{\sampl,x^{\alpha}}(n;\, s) = g_{x^{\alpha}}(n;\, s),
\end{equation}
\item
  convolutions with integrated Gaussian derivative kernels according to
  (Lindeberg \citeyear{Lin24-JMIV} Equation~(54))
   (see also Appendix~\ref{app-int-gaussder-kern})
  \begin{equation}
  \label{eq-def-int-gauss-der}
  T_{\intdisc,x^{\alpha}}(n;\, s) = \int_{x = n - 1/2}^{n + 1/2} g_{x^{\alpha}}(x;\; s) \, dx,
\end{equation}
\item
  the genuinely discrete approach corresponding to convolution with the discrete analogue of the Gaussian kernel $T_{\disc}(n;\; s)$ according to
  (\ref{eq-disc-gauss}) followed by central difference operators $\delta_{x^{\alpha}}$, thus corresponding to the equivalent discrete approximation kernel
  (Lindeberg \citeyear{Lin93-JMIV} Equation~(58))
   (see also Appendix~\ref{app-disc-gauss-ders})
\begin{equation}
  \label{eq-disc-der-gauss}
  T_{\disc,x^{\alpha}}(n;\; s) = (\delta_{x^{\alpha}} T_{\disc})(n;\; s).
\end{equation}
\end{itemize}
Here, the central difference operators are for orders 1 and 2 defined according to
\begin{align}
  \begin{split}
    \delta_x & = (-\tfrac{1}{2}, 0, +\tfrac{1}{2}),
   \end{split}\\
 \begin{split}
    \delta_{xx} & = (+1, -2, +1),
   \end{split}
\end{align}
and for higher values of $\alpha$ according to:
\begin{equation}
  \label{eq-def-cent-diff-op-arb-order}
  \delta_{x^{\alpha}}
  = \left\{
         \begin{array}{ll}
           \delta_x (\delta_{xx})^i & \mbox{if $\alpha = 1 + 2 i$,} \\
            (\delta_{xx})^i & \mbox{if $\alpha = 2 i$,}
          \end{array}
        \right.
\end{equation}
for integer $i$, where the special cases $\alpha = 3$ and $\alpha = 4$
then correspond to the difference operators
\begin{align}
  \begin{split}
    \delta_{xxx} & = (-\tfrac{1}{2}, +1, 0, -1, +\tfrac{1}{2}),
   \end{split}\\
 \begin{split}
    \delta_{xxxx} & = (+1, -4, +6, -4, +1).
   \end{split}
\end{align}
In addition to the above, already studied discretisation methods in
(Lindeberg \citeyear{Lin24-JMIV}),
we will here specifically consider the properties of the following
hybrid methods, in terms of either:
\begin{itemize}
\item
  the hybrid approach corresponding to convolution with the normalised
  sampled Gaussian kernel $T_{\normsampl}(n;\; s)$ according to
  (\ref{eq-def-norm-sampl-gauss}) followed by central difference
  operators $\delta_{x^{\alpha}}$ according to  (\ref{eq-def-cent-diff-op-arb-order}),
  thus corresponding to the equivalent discrete approximation kernel
  (Lindeberg \citeyear{Lin24-JMIV} Equation~(116))
  (see also Appendix~\ref{app-hybr-disc-normsamplgaussdiff})
\begin{equation}
      \label{eq-hybr-normsampl-disc-der}
      T_{\hybrnormsampl,x^{\alpha}}(n;\; s)
       = (\delta_{x^{\alpha}} T_{\normsampl})(n;\; s),
\end{equation}
\item
 the hybrid approach corresponding to convolution with the integrated
 Gaussian kernel $T_{\intdisc}(n;\; s)$ according to
 (\ref{eq-def-int-gauss-kern}) followed by central difference
 operators $\delta_{x^{\alpha}}$ according to (\ref{eq-def-cent-diff-op-arb-order}),
 thus corresponding to the equivalent discrete approximation kernel
 (Lindeberg \citeyear{Lin24-JMIV} Equation~(117))
 (see also Appendix~\ref{app-hybr-disc-intgaussdiff})
 \begin{equation}
      \label{eq-hybr-int-disc-der}
      T_{\hybrint,x^{\alpha}}(n;\; s)
       = (\delta_{x^{\alpha}} T_{\intdisc})(n;\; s).
\end{equation}        
\end{itemize}
A motivation for introducing these hybrid discretisation methods
(\ref{eq-hybr-normsampl-disc-der}) and (\ref{eq-hybr-int-disc-der}),
based on convolutions with the normalised sampled Gaussian kernel
(\ref{eq-def-norm-sampl-gauss}) or the integrated Gaussian kernel
(\ref{eq-def-int-gauss-kern}) followed by central difference operators
of the form (\ref{eq-def-cent-diff-op-arb-order}), is that these
discretisation methods have substantially better computational
efficiency, compared to explicit convolutions with either the
sampled Gaussian derivative kernels
(\ref{eq-sampl-gauss-der}) or the integrated Gaussian derivative
kernels (\ref{eq-def-int-gauss-der}), in situations when spatial
derivatives of multiple orders $\alpha$ are needed at the same scale level.

The reason for this is that the same spatial smoothing stage can then
be shared between the computations of discrete derivative
approximations for the different orders of spatial differentiation,
thus implying that these hybrid methods will be as computationally
efficient as the genuinely discrete approach, based on convolution
with the discrete analogue of the Gaussian kernel
(\ref{eq-disc-gauss}) followed by central differences of the form
(\ref{eq-def-cent-diff-op-arb-order}), and corresponding to equivalent
convolution kernels of the form (\ref{eq-disc-der-gauss}).

\subsubsection{Quantitative measures of approximation properties
  relative to continuous scale space}

To measure how well the above discretisation methods for the Gaussian
derivative operators 
reflect properties of the underlying continuous Gaussian derivatives,
we will consider quantifications in terms of the following
spatial spread measure, defined from spatial variance $V$ of the absolute value of each
equivalent discrete derivative approximation kernel
(Lindeberg \citeyear{Lin24-JMIV} Equation~(80)):
  \begin{equation}
    \label{eq-def-rel-scale-err-gauss-ders}
     \sqrt{V(|T_{x^{\alpha}}(\cdot;\; s)|)},
\end{equation}
where the variance $V(h(\cdot))$ of a non-negative continuous 1-D function $h(x)$ is defined as
\begin{equation}
  V(h(\cdot)) %= \\
  = \frac{\int_{x \in \bbbr} x^2 \, h(x) \, dx}{\int_{x \in \bbbr} h(x) \, dx}
  - \left(
        \frac{\int_{x \in \bbbr} x \, h(x) \, dx}{\int_{x \in \bbbr} h(x) \, dx}
     \right)^2.
\end{equation}
To furthermore more explicitly quantify the deviation from the corresponding fully continuous spatial spread measures $\sqrt{V(|g_{x^{\alpha}}(\cdot;\; s)|)}$, we also
consider the following measures of the offsets of the spatial spread measures
(Lindeberg \citeyear{Lin24-JMIV} Equation~(81)):
\begin{equation}
    \label{eq-spat-spread-meas-offset}
    O_{\alpha}(s)
    = \sqrt{V(|T_{x^{\alpha}}(\cdot;\; s)|)} - \sqrt{V(|g_{x^{\alpha}}(\cdot;\; s)|)},
  \end{equation}
where the variance $V(h(\cdot))$ of a non-negative discrete function $h(n)$ is defined as
\begin{equation}
  V(h(\cdot)) %= \\
  = \frac{\sum_{n \in \bbbz} n^2 \, h(n)}{\sum_{n \in \bbbz} h(n)}
  - \left(
    \frac{\sum_{n \in \bbbz} n \, h(n)}{\sum_{n \in \bbbz} h(n)}
  \right)^2.
\end{equation}

\subsection{Methodology for characterising the resulting consistency properties over scale in terms of the accuracy of the scale estimates obtained from integrations with scale selection algorithms}
\label{sec-meth-consist-scales}

To perform a further evaluation of the hybrid discretisation method
to consistently process input data over multiple scales,
we will characterise the abilities of these methods in a context
of feature detection with automatic scale selection
(Lindeberg \citeyear{Lin97-IJCV,Lin97-IJCV,Lin21-EncCompVis}),
where hypotheses about local
appropriate scale levels are determined from local extrema
over scale of scale-normalised derivative responses.

For this purpose, we follow a similar methodology as used in 
(Lindeberg \citeyear{Lin24-JMIV} Section~4).
Thus, with the theory in Section~\ref{sec-discr-methods} %in this paper
now applied to
2-D image data, by separable application of the 1-D theory
along each image dimension, we consider scale-normalised derivative operators
defined according to (Lindeberg \citeyear{Lin97-IJCV,Lin98-IJCV})
(for $(x, y) \in \bbbr^2$ and $s \in \bbbr_+$):
\begin{equation}
  \label{eq-def-sc-norm-ders}
  \partial_{\xi} = s^{\gamma/2} \, \partial_x, \quad\quad
  \partial_{\eta} = s^{\gamma/2} \, \partial_y, 
\end{equation}
with $\gamma > 0$ being a scale normalisation power,
that is chosen specially adapted for each feature detection task.

\subsubsection{Scale-invariant feature detectors with automatic scale selection}
\label{sec-sc-inv-feat-det}

Specifically, we will evaluate the performance of
the following types of scale-invariant feature detectors,
defined from the spatial derivatives $L_x(x, y;\; s)$,
$L_y(x, y;\; s)$, $L_{xx}(x, y;\; s)$, $L_{xy}(x, y;\; s)$ and
$L_{yy}(x, y;\; s)$ up to order 2 of the 2-D scale-space
representation $L(x, y;\; s)$ of any 2-D image $f(x, y)$
at scale $s$, obtained by convolution with 2-D Gaussian kernels
$g_{2D}(x, y;\; s)$ for different values of $s$:
\begin{itemize}
\item
  {\em interest point detection} from scale-space extrema
  (extrema over both space $(x, y)$ and scale $s$) of
  the scale-normalised Laplacian operator
  (Lindeberg \citeyear{Lin97-IJCV} Equation~(30))
\begin{equation}
  \label{eq-sc-norm-lapl}
  \nabla_{\norm}^2 L = s \, (L_{xx} + L_{yy}),
\end{equation}
or the scale-normalised determinant of the Hessian operator
(Lindeberg \citeyear{Lin97-IJCV} Equation~(31))
\begin{equation}
  \label{eq-sc-norm-det-hess}  
  \det {\mathcal H}_{\norm} L = s^2 \, (L_{xx} \, L_{yy} - L_{xy}^2),
\end{equation}
where we here choose the scale normalisation parameter $\gamma = 1$, such that the
selected scale level for a Gaussian blob of size $s_0 \in \bbbr_+$
\begin{equation}
  \label{eq-def-2D-gauss-blob-s0}
  f_{\blob,s_0}(x, y)
  = g_{2D}(x, y;\; s_0)
  = \frac{1}{2 \pi s_0} \, e^{-\frac{x^2 + y^2}{2s_0}}
\end{equation}
will for both Laplacian and determinant of the Hessian interest point detection be equal to the size of the blob
(Lindeberg \citeyear{Lin97-IJCV}, Equations~(36) and~(37)):
\begin{align}
  \begin{split}
    \label{eq-lapl-sc-sel}
    (\hat{x}, \hat{y};\; \hat{s}) 
    & = \operatorname{argmin}_{(x, y;\; s)}(\nabla_{\norm}^2 L)(x, y;\; s)
  \end{split}\nonumber\\
  \begin{split}
    &
    = (0, 0;\; s_0),
  \end{split}\\
  \begin{split}
    \label{eq-det-hess-sc-sel}
    (\hat{x}, \hat{y}; \hat{s}) 
    & = \operatorname{argmax}_{(x, y;\; s)}(\det {\mathcal H}_{\norm} L)(x, y;\; s)
  \end{split}\nonumber\\
  \begin{split}
    &
    = (0, 0;\; s_0),
  \end{split}
\end{align}
% with the scale-normalized Laplacian response $(\nabla_{\norm} ^2 L)(x, y;\; s)$
%at any image point $(x, y)$ and scale $s$ according to (\ref{eq-sc-norm-lapl})
%and the scale-normalized determinant of the Hessian response
%$(\det {\mathcal H}_{\norm} L)(x, y;\; s)$ at any image point $(x, y)$
% and scale $s$ according to (\ref{eq-sc-norm-det-hess}),
\item
  {\em edge detection} from combined
  \begin{itemize}
  \item
  maxima of the gradient magnitude in the spatial gradient direction $e_v$,
  reformulated such that the second-order directional derivative in
  the gradient direction $L_{vv}$ is zero and the third-order directional derivative in the gradient direction $L_{vvv}$ is negative (Lindeberg \citeyear{Lin98-IJCV} Equation~(8)), and
\item
  maxima over scale of the
  scale-normalised gradient magnitude $L_{v,\norm}$ according to
(Lindeberg \citeyear{Lin98-IJCV} Equation~(15))
\begin{equation}
  \label{eq-sc-norm-grad-magn}
  L_{v,\norm} = s^{\gamma/2} \sqrt{L_x^2 + L_y^2},
\end{equation}
\end{itemize}
where we here set the scale normalisation parameter $\gamma$ to
(Lindeberg \citeyear{Lin98-IJCV}, Equation~(23))
\begin{equation}
   \label{eq-gamma-edge-det}
  \gamma_{\edge} = \frac{1}{2},
\end{equation}
such that the selected scale level $\hat{s}$ for
an idealised model of a diffuse edge
(Lindeberg \citeyear{Lin98-IJCV} Equation~(18))
\begin{multline}
  \label{eq-def-ideal-blur-edge-s0}  
  f_{\edge,s_0}(x, y)
  =  \int_{u = - \infty}^x g_{1D}(u;\; s_0) \, du = \\
  = \int_{u = - \infty}^x
        \frac{1}{\sqrt{2 \pi s_0}} \, e^{-\frac{u^2}{2 s_0}} \, du
\end{multline}
will be equal to the amount of diffuseness $s_0 \in \bbbr_+$ of that diffuse edge
\begin{equation}
  \label{eq-grad-magn-sc-sel}
  \hat{s} = \operatorname{argmax}_s L_{v,\norm}(0, 0;\; s)  = s_0,
\end{equation}
%with the scale-normalized gradient magnitude response
%$L_{v,\norm}(0, 0;\; s)$ at the image center  $(0, 0)$ according
%to (\ref{eq-sc-norm-grad-magn}),
\item
  {\em ridge detection} from combined
  \begin{itemize}
  \item
    zero-crossings of the first-order directional derivative
    $L_{p}$ in the first principal curvature direction $e_p$ of the Hessian matrix, such that $L_{p} = 0$
    (Lindeberg \citeyear{Lin98-IJCV} Equations~(42)), and
  \item
  minima over scale of the scale-normalised ridge strength in terms the
scale-normalised second-order derivative $L_{pp,\norm}$ in the direction $e_p$ according to (Lindeberg \citeyear{Lin98-IJCV} Equation~(47)):
\begin{multline}
  \label{eq-sc-norm-ridge-meas}
  L_{pp,\norm}
  = s^{\gamma} L_{pp} = \\
  = s^{\gamma}
      \left(
        L_{xx} + L_{yy} - \sqrt{(L_{xx} - L_{yy})^2 + 4 L_{xy}^2}
      \right),
    \end{multline}
  \end{itemize}
  where we here choose the scale normalisation parameter $\gamma$ as
(Lindeberg \citeyear{Lin98-IJCV}, Equation~(56))
\begin{equation}
  \label{eq-gamma-ridge-det}
  \gamma_{\ridge} = \frac{3}{4},
\end{equation}
such that the selected scale level $\hat{s}$ for a Gaussian ridge model of the form
\begin{equation}
  \label{eq-def-ideal-ridge-s0}  
  f_{\ridge,s_0}(x, y)
  = g_{1D}(x;\; s_0)
  = \frac{1}{\sqrt{2 \pi s_0}} \, e^{-\frac{x^2}{2 s_0}} 
\end{equation}
will be equal to the width $s_0 \in \bbbr_+$ of that idealised ridge model
\begin{equation}
   \label{eq-princ-curv-sc-sel}
  \hat{s} = \operatorname{argmax}_s L_{pp,\norm}(0, 0;\; s) = s_0.
\end{equation}
%with the scale-normalized ridge response $L_{pp,\norm}(0, 0;\; s)$
%at the image center $(0, 0)$ according
%to (\ref{eq-sc-norm-ridge-meas}).
\end{itemize}
A common property of all these scale-invariant feature detectors is, thus, that
the selected scale levels $\hat{s}$ obtained from local extrema over scale
will reflect characteristic%
\footnote{The notion of ``characteristic scale'' refers to a scale,
  that reflects a characteristic length in the image data, in a
  similar way as the notion of characteristic length is used in the
  areas of physics. See Appendix~\ref{app-char-scale}
  for further details.}
scales $s_0$ in the input data
(Lindeberg \citeyear{Lin21-EncCompVis}).
By evaluating discretisation methods of Gaussian derivatives with respect
to such scale selection properties, we therefore have a way of formulating
a well-defined proxy task, for evaluating how well the different types
of discretisation methods lead to appropriate consistency properties over
scales for the numerical implementations of Gaussian derivative operators.

Figures~\ref{fig-lapl-intpt-det}--\ref{fig-princcurv-ridge-det} in Appendix~\ref{sec-visualisation} provide visualisations of
the conceptual steps involved when defining these scale estimates $\hat{s}$.

\subsubsection{Quantitative measures for characterising the accuracy of
  the scale estimates obtained from the scale selection methodology}
\label{sec-quant-meas-sc-sel-perf}

To quantify the performance of the different discretisation methods
with regard to the above scale selection tasks,
we will
\begin{itemize}
\item
  compute the selected scale levels $\hat{\sigma} = \sqrt{\hat{s}}$
  for different values of
the characteristic scale $s_0$ in the image data, measured in dimension length
$\sigma_0 = \sqrt{s_0}$, and
\item
  quantify the deviations from the reference in terms of the relative error measure
  (Lindeberg \citeyear{Lin24-JMIV} Equation~(107))
    \begin{equation}
    \label{eq-def-sc-sel-rel-sc-err}
    E_{\scaleestrel}(\sigma) = \frac{\hat{\sigma}}{\hat{\sigma}_{\scaleref}} - 1,
    \end{equation}
\end{itemize}
under variations of the characteristic scale $s_0$ in the
input image, where the deviations between the selected scale levels
$\hat{\sigma} = \sqrt{\hat{s}}$ and the reference value $\sigma_0 = \sqrt{s_0}$ are to be interpreted
as the results of discretisation errors.

When generating discrete model signals for the different discrete
approximation methods, we use as conceptually close discretisation
methods for the input model signals (Gaussian blobs for interest point detection,
diffuse edges for edge detection, or Gaussian ridges for ridge
detection) as for the discrete approximations of
Gaussian derivatives, according to Appendix~\ref{sec-def-model-signals}.

% \footnote{Provided that we would not want to expand the experiments,
%   to all possible combinations of discretisation methods regarding both
%   the input data and image operations, which would then also make the
%   analysis and the interpretation steps more complex, this should,
%   however, constitute an appropriate matching regarding the
%   definition of the input data and the selection of discrete approximation method.}

\begin{figure*}[hbpt]

  \begin{center}

     \begin{subfigure}[t]{0.45\textwidth}
        \begin{tabular}{c}
            {\em Spatial spread measures for 1st-order derivative kernels\/} \\
            \includegraphics[width=0.97\textwidth]{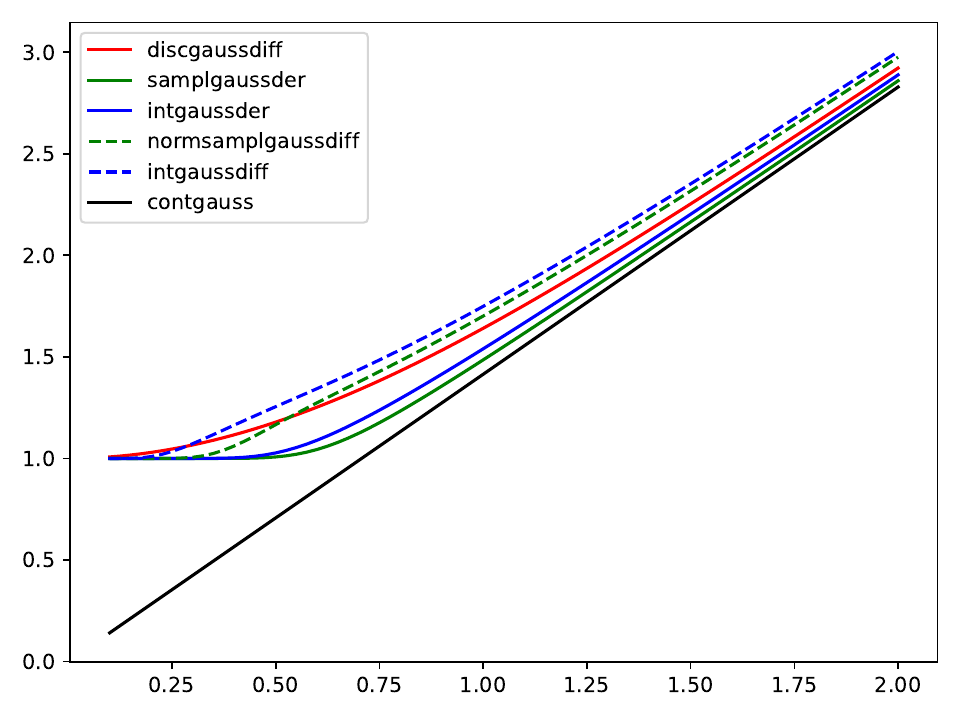}
         \end{tabular}
         \caption{Case: $\alpha = 1$.
    Note that the spatial spread measures for the spatial smoothing kernels combined with central differences are delimited from below by the
    spatial variance of the absolute value of the first-order
    central difference operator $|\delta_x|$, which is $\sqrt{V(|\delta_x|)} = 1$.}
      \label{fig-gauss-der1-relscaleerr}
     \end{subfigure}%
     ~\quad\quad~
    \begin{subfigure}[t]{0.45\textwidth}
       \begin{tabular}{c}
          {\em Spatial spread measures for 2nd-order derivative kernels\/} \\
          \includegraphics[width=0.97\textwidth]{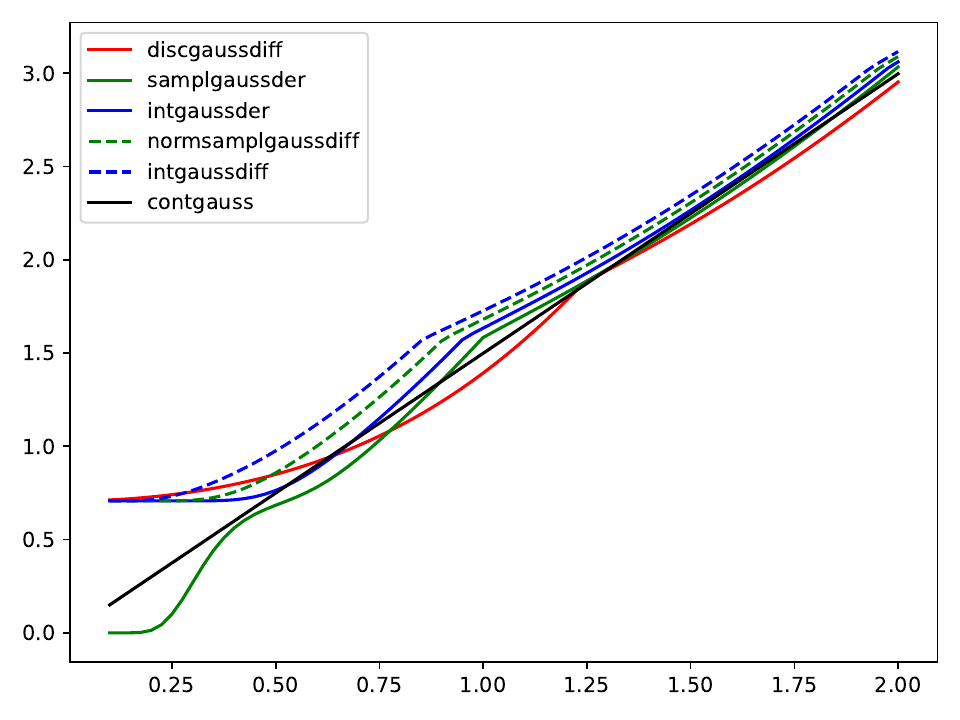}
        \end{tabular}
        \caption{Case: $\alpha = 2$. 
    Note that the spatial spread measures for the spatial smoothing kernels combined with central differences are delimited from below by the
    spatial variance of the absolute value of the second-order
    central difference operator $|\delta_{xx}|$,
    which is $\sqrt{V(|\delta_{xx}|)} = 1/\sqrt{2}$.}
      \label{fig-gauss-der2-relscaleerr}
    \end{subfigure}

    \bigskip
    \bigskip    

    \begin{subfigure}[t]{0.45\textwidth}
       \begin{tabular}{c}
          {\em Spatial spread measures for 3rd-order derivative kernels\/} \\
          \includegraphics[width=0.97\textwidth]{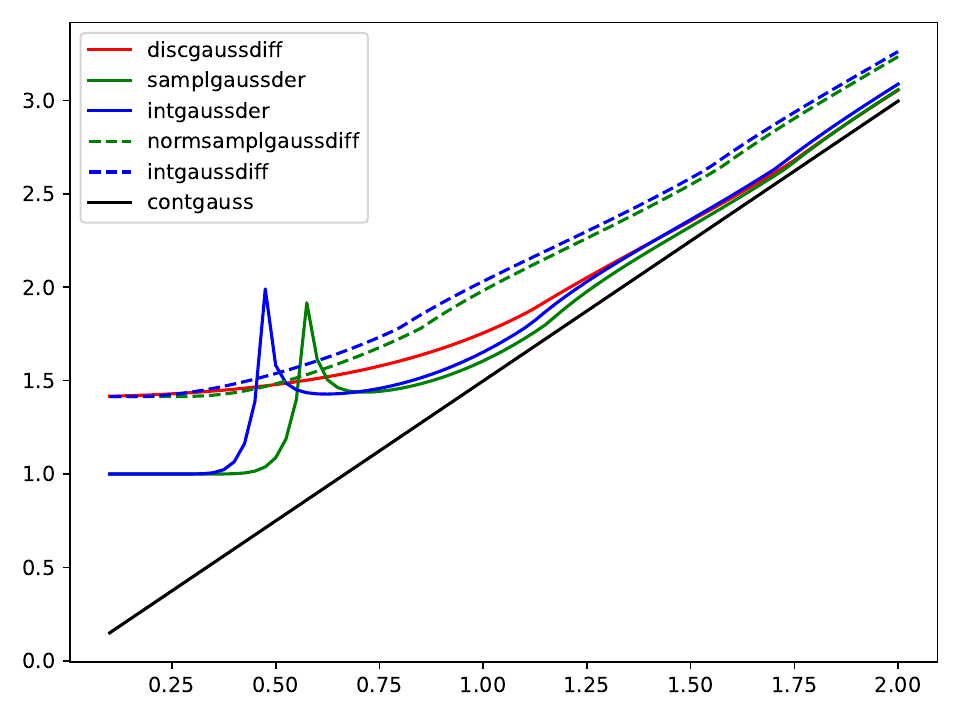}
       \end{tabular}
       \caption{Case: $\alpha = 3$.  
    Note that the spatial spread measures for the spatial smoothing kernels combined with central differences are delimited from below by the
    spatial variance of the absolute value of the third-order
    central difference operator $|\delta_{xxx}|$,
    which is $\sqrt{V(|\delta_{xxx}|)} = \sqrt{2}$.}
     \label{fig-gauss-der3-relscaleerr}
    \end{subfigure}%
     ~\quad\quad~
    \begin{subfigure}[t]{0.45\textwidth}
       \begin{tabular}{c}
           {\em Spatial spread measures for 4th-order derivative kernels\/} \\
           \includegraphics[width=0.97\textwidth]{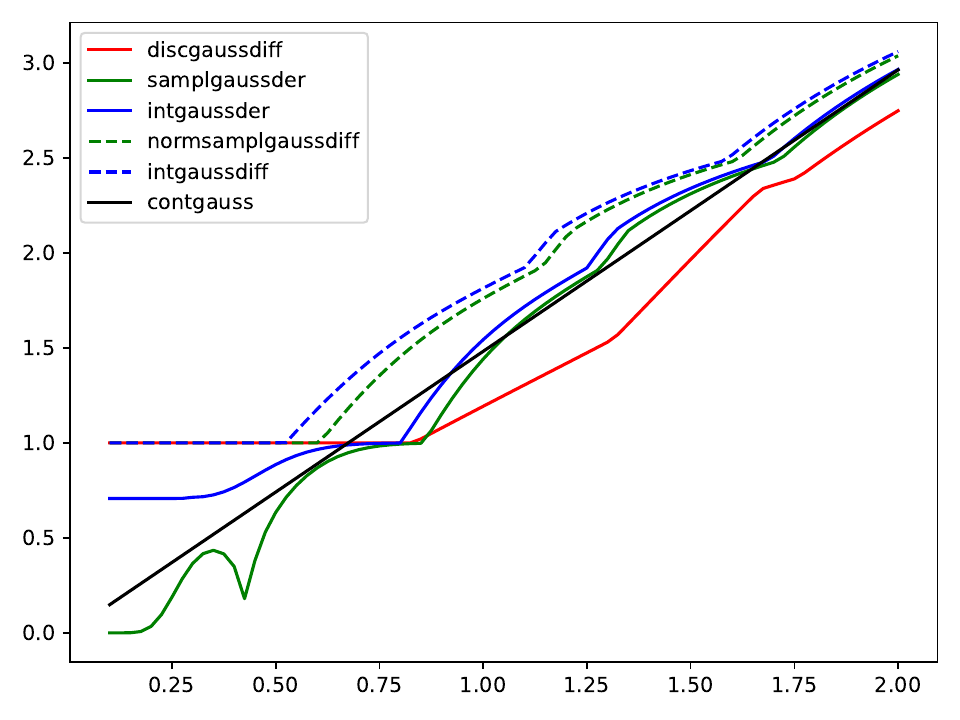}
       \end{tabular}
       \caption{Case: $\alpha = 4$. 
    Note that the spatial spread measures for the spatial smoothing kernels combined with central differences are delimited from below by the
    spatial variance of the absolute value of the fourth-order
    central difference operator $|\delta_{xxxx}|$,
    which is $\sqrt{V(|\delta_{xxxx}|)} = 1$.}
     \label{fig-gauss-der4-relscaleerr}
     \end{subfigure}
   \end{center}
   
   \caption{Graphs of the {\em spatial spread measure\/}
    $\sqrt{V(|T_{x^{\alpha}}(\cdot;\; s)|)}$ according to
    (\ref{eq-def-rel-scale-err-gauss-ders}) for
    different discrete approximations of Gaussian
    derivative kernels of order $\alpha$: {\bf (i)}~for either discrete analogues of Gaussian
    derivative kernels $T_{\disc,x^{\alpha}}(n;\; s)$
    according to (\ref{eq-disc-der-gauss}), corresponding to convolutions with the discrete analogue of the Gaussian kernel $T_{\disc}(n;\; s)$ according to
    (\ref{eq-disc-gauss}) followed by central differences
    according to (\ref{eq-def-cent-diff-op-arb-order}), {\bf (ii)}~sampled Gaussian
    derivative kernels $T_{\sampl,x^{\alpha}}(n;\, s)$ according to
    (\ref{eq-sampl-gauss-der}), {\bf (iii)}~integrated Gaussian derivative
    kernels $T_{\intdisc,x^{\alpha}}(n;\, s)$ according to
    (\ref{eq-def-int-gauss-der}), {\bf (iv)}~the hybrid discretisation kernel
    $T_{\hybrnormsampl,x^{\alpha}}(n;\; s)$ according to (\ref{eq-hybr-normsampl-disc-der}), corresponding to convolution with the normalised sampled Gaussian kernel $T_{\normsampl}(n;\; s)$ according to (\ref{eq-def-norm-sampl-gauss}) 
    followed by central differences according to (\ref{eq-def-cent-diff-op-arb-order}), and
    {\bf (v)}~the hybrid discretisation kernel $T_{\hybrint,x^{\alpha}}(n;\; s)$
    according to (\ref{eq-hybr-int-disc-der}), corresponding to convolution with the integrated Gaussian kernel $T_{\intdisc}(n;\; s)$ according to
    (\ref{eq-def-int-gauss-kern}) followed by central differences according to (\ref{eq-def-cent-diff-op-arb-order}).
    ({\bf Horizontal axes:} Scale parameter in units of
    $\sigma = \sqrt{s} \in [0.1, 2]$.)}
  \label{fig-spat-spread-meas-with-hybr-discr}
\end{figure*}

\begin{figure*}[hbpt]

  \begin{center}

     \begin{subfigure}[t]{0.45\textwidth}
        \begin{tabular}{c}
            {\em Spread measure offsets for 1st-order derivative kernels\/} \\
            \includegraphics[width=0.97\textwidth]{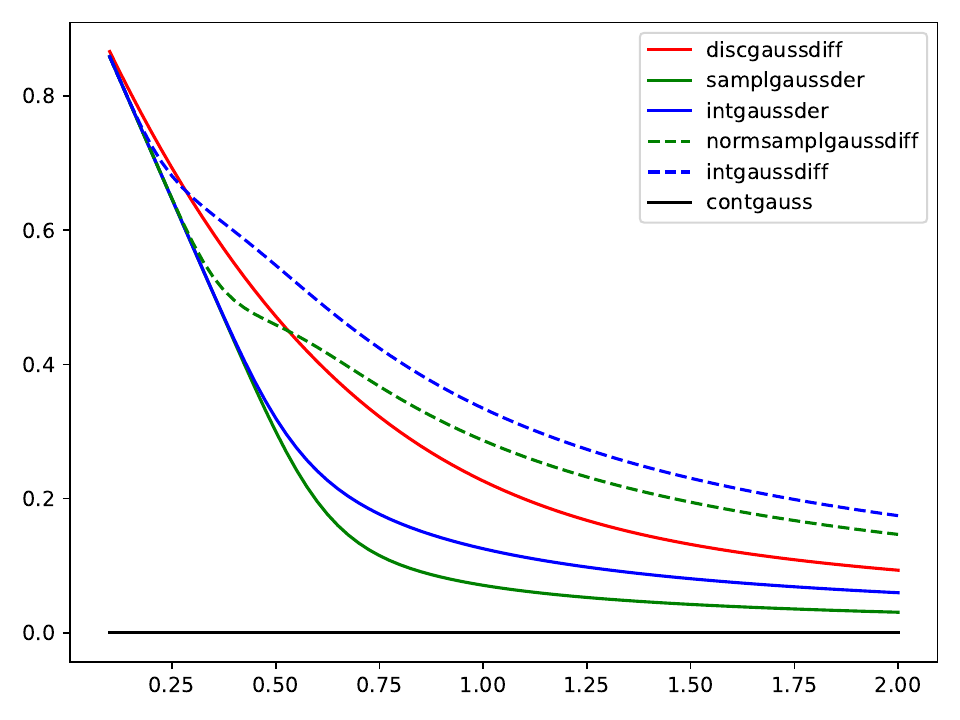}
         \end{tabular}
         \caption{Case: $\alpha = 1$.}
      \label{fig-gauss-der1-relscaleerr-diff}
     \end{subfigure}%
     ~\quad\quad~
    \begin{subfigure}[t]{0.45\textwidth}
       \begin{tabular}{c}
          {\em Spread measure offsets for 2nd-order derivative kernels\/} \\
          \includegraphics[width=0.97\textwidth]{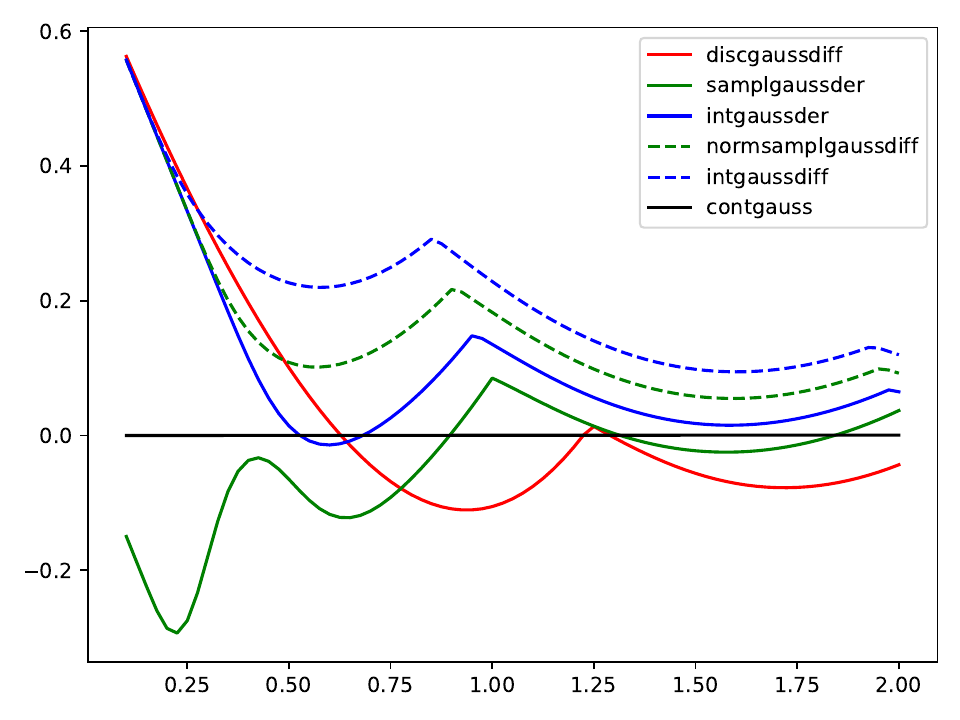}
        \end{tabular}
        \caption{Case: $\alpha = 2$.}
      \label{fig-gauss-der2-relscaleerr-diff}
    \end{subfigure}

    \bigskip
    \bigskip    

    \begin{subfigure}[t]{0.45\textwidth}
       \begin{tabular}{c}
          {\em Spread measure offsets for 3rd-order derivative kernels\/} \\
          \includegraphics[width=0.97\textwidth]{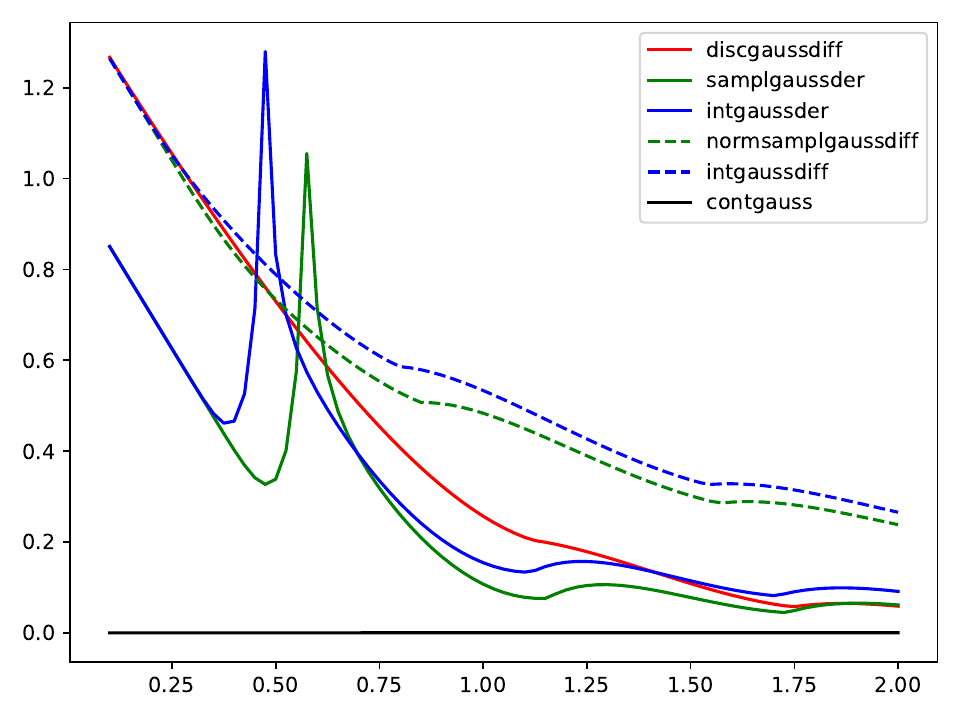}
       \end{tabular}
       \caption{Case: $\alpha = 3$.}
     \label{fig-gauss-der3-relscaleerr-diff}
    \end{subfigure}%
     ~\quad\quad~
    \begin{subfigure}[t]{0.45\textwidth}
       \begin{tabular}{c}
           {\em Spread measure offsets for 4th-order derivative kernels\/} \\
           \includegraphics[width=0.97\textwidth]{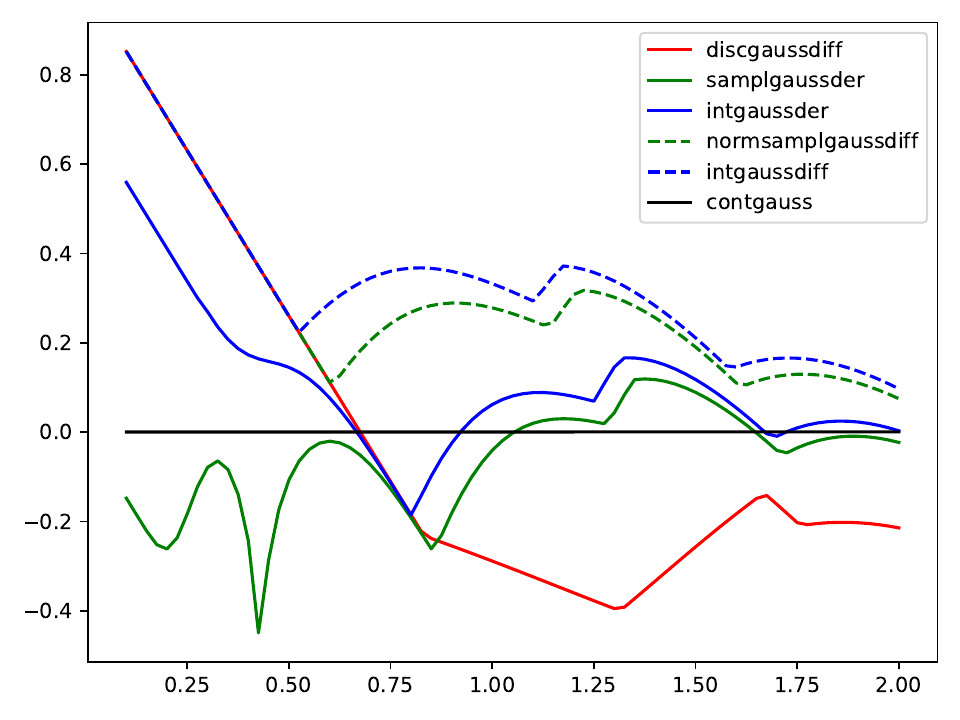}
       \end{tabular}
       \caption{Case: $\alpha = 4$.}
     \label{fig-gauss-der4-relscaleerr-diff}
     \end{subfigure}
   \end{center}
    
    \caption{Graphs of the {\em spatial spread measure offset\/}
     $O_{\alpha}(s)$, relative to the spatial spread of a continuous
     Gaussian kernel, according to (\ref{eq-spat-spread-meas-offset}), for
    different discrete approximations of Gaussian
    derivative kernels of order $\alpha$:
{\bf (i)}~for either discrete analogues of Gaussian
    derivative kernels $T_{\disc,x^{\alpha}}(n;\; s)$
    according to (\ref{eq-disc-der-gauss}), corresponding to convolutions with the discrete analogue of the Gaussian kernel $T_{\disc}(n;\; s)$ according to
    (\ref{eq-disc-gauss}) followed by central differences according to
    (\ref{eq-def-cent-diff-op-arb-order}), {\bf (ii)}~sampled Gaussian
    derivative kernels $T_{\sampl,x^{\alpha}}(n;\, s)$ according to
    (\ref{eq-sampl-gauss-der}), {\bf (iii)}~integrated Gaussian derivative
    kernels $T_{\intdisc,x^{\alpha}}(n;\, s)$ according to
    (\ref{eq-def-int-gauss-der}), {\bf (iv)}~the hybrid discretisation kernel
    $T_{\hybrnormsampl,x^{\alpha}}(n;\; s)$ according to (\ref{eq-hybr-normsampl-disc-der}), corresponding to convolution with the normalised sampled Gaussian kernel $T_{\normsampl}(n;\; s)$ according to (\ref{eq-def-norm-sampl-gauss}) 
    followed by central differences according to
    (\ref{eq-def-cent-diff-op-arb-order}), and
    {\bf (v)}~the hybrid discretisation kernel $T_{\hybrint,x^{\alpha}}(n;\; s)$
    according to (\ref{eq-hybr-int-disc-der}), corresponding to convolution with the integrated Gaussian kernel $T_{\intdisc}(n;\; s)$ according to
    (\ref{eq-def-int-gauss-kern}) followed by central differences according to (\ref{eq-def-cent-diff-op-arb-order}).
    ({\bf Horizontal axes:} Scale parameter in units of
    $\sigma = \sqrt{s} \in [0.1, 2]$.)}
    \label{fig-spat-spread-offset-meas-with-hybr-discr}
\end{figure*}

\section{Results}
\label{sec-results}

\subsection{Characterisation of the effective amount of spatial smoothing in discrete approximations of Gaussian derivatives in terms of spatial spread measures}

Figures~\ref{fig-spat-spread-meas-with-hybr-discr}--\ref{fig-spat-spread-offset-meas-with-hybr-discr}
show the graphs of computing the spatial spread measure
$\sqrt{V(|T_{x^{\alpha}}(\cdot;\; s)|)} $ according to
(\ref{eq-def-rel-scale-err-gauss-ders}) as well as the offset measure
$O_{\alpha}(s)$ according to (\ref{eq-spat-spread-meas-offset}) over
an interval of finer scale values $\sigma = \sqrt{s} \in [0.1, 2]$, for
each one of the different discretisation methods described in
Section~\ref{sec-discr-methods}.

As can be seen from
these graphs:
\begin{itemize}
\item
  The agreement with the underlying fully continuous spread measures for the continuous Gaussian derivative kernels is substantially better for the genuinely sampled or integrated Gaussian derivative kernels than for the hybrid discretisations based on combining either the normalised sampled Gaussian kernel or the integrated Gaussian kernel with central difference operators.

%  The hybrid discretisation kernels, corresponding to convolutions with either the normalised sampled Gaussian kernel or the integrated Gaussian kernel followed by central differences, generally have substantially larger offsets from the underlying continuous theory than the more direct discretisations in terms of either sampled Gaussian derivative kernels or integrated Gaussian kernels.

  In situations when multiple derivatives of different orders $\alpha$ are to be computed at the same scale levels, the hybrid discretisation methods are, however, as previously mentioned, computationally much more efficient, implying that the introduction of the hybrid discretisations implies a trade-off between the accuracy in terms of the overall amount of spatial smoothing of the equivalent discrete filters and the computational efficiency of the implementation.
\item
  The agreement with the underlying fully continuous spread measures
  for the continuous Gaussian derivative kernels is substantially
  better for the genuinely discrete analogue of Gaussian derivative
  operators, obtained by first convolving the input data with the
  discrete analogue of the Gaussian kernel and then applying central
  difference operators to the spatially smoothing input data, compared
  to using any of the hybrid discretisations. % corresponding to first smoothing the input data with either the normalised sampled Gaussian kernel or the integrated Gaussian kernel, and then applying central difference operators to the spatially smoothed input data.

  If, for efficiency reasons, a discretisation method is to be chosen, based on combining a first stage of spatial smoothing with a following application of central difference operators, the approach based on using spatial smoothing with the discrete analogue of the Gaussian kernel, in most of the cases, leads to better agreement with the underlying continuous theory, compared to using either the normalised sampled Gaussian kernel or the integrated Gaussian kernel in the first stage of spatial smoothing.

  As previously stated, the hybrid discretisation methods may,
  however, anyway be warranted in situations where the underlying
  modified Bessel functions $I_n(s)$ are not fully available in the computational environment, where the discrete filtering operations are to be implemented, such as when performing learning of the scale levels in deep networks based on Gaussian derivative operators.
\end{itemize}
A further general implication of these results is that, depending on what discretisation method is chosen for discretising the computation of Gaussian derivative responses at fine scales, different values of the spatial scale parameter $s$ will be needed, to obtain a comparable amount of spatial smoothing of the input data for the different
discretisation methods.

\subsection{Characterisation of the approximation properties relative
  to continuous scale space in terms of the scale levels selected by
  scale selection algorithms}
\label{sec-results-scsel}

We will next characterise the approximation properties of the scale estimates
for the different benchmark tasks outlined in
Section~\ref{sec-sc-inv-feat-det} with regard to the quantitative
measures defined in Section~\ref{sec-quant-meas-sc-sel-perf}.

For generating the input data, we used 50 logarithmically scale values
$\sigma_0 \in [1/3, 3]$. When performing the scale selection step,
we searched over a range of 80 logarithmically scale levels
$\sigma_0 \in [0.1, 5]$, and accumulated the variability
over scale at the image center for each differential
feature detector, and additionally performed parabolic
interpolation over the logarithmic scale values, to localise the
scale estimates to higher accuracy.
This very dense sampling of the scale levels is far beyond what is usually
needed in actual image processing or computer vision algorithms, but was chosen here
in order to essentially eliminate the effects of discrete sampling issues in the
scale direction.

\begin{figure*}[hbpt]

  \begin{center}
     \begin{subfigure}[t]{0.45\textwidth}
        \begin{tabular}{c}
           {\em Selected scales for Laplacian scale selection \/} \\
           \includegraphics[width=0.97\textwidth]{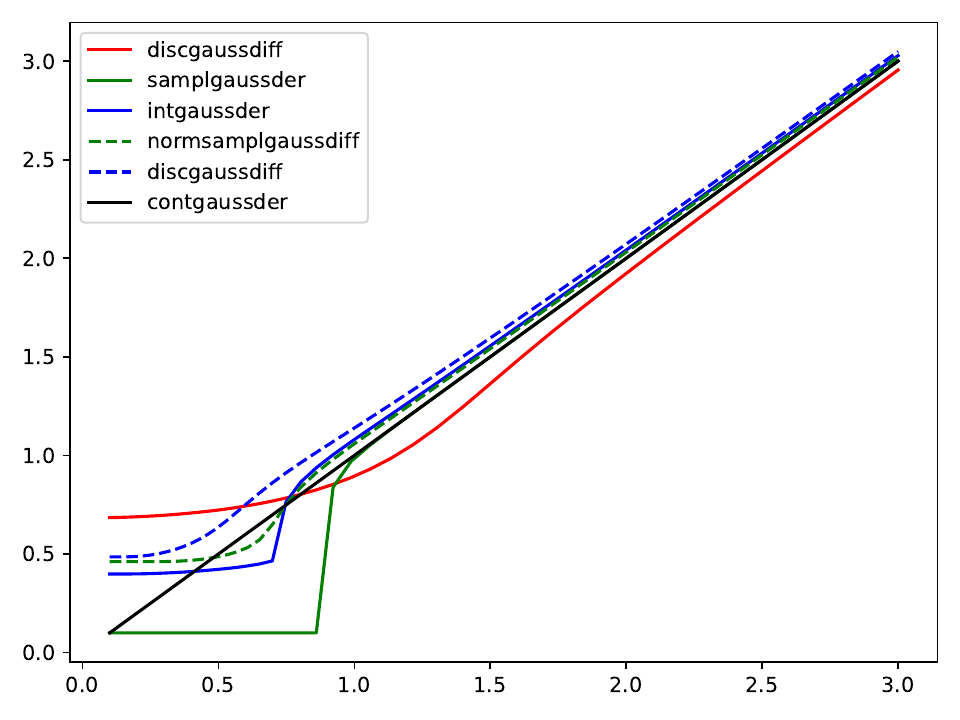}
        \end{tabular}
        \caption{Selected scales \, $\hat{\sigma}$.}
        %\caption{Selected scales $\hat{\sigma} = \sqrt{\hat{s}}$ from
        % $\argmin \nabla_{\norm}^2 L$.}
         \label{fig-sel-sc-Lapl-scsel-gauss-blob}
     \end{subfigure}%
     ~\quad\quad~
    \begin{subfigure}[t]{0.45\textwidth}
       \begin{tabular}{c}
            {\em Selected scales for gradient magnitude scale selection \/} \\
            \includegraphics[width=0.97\textwidth]{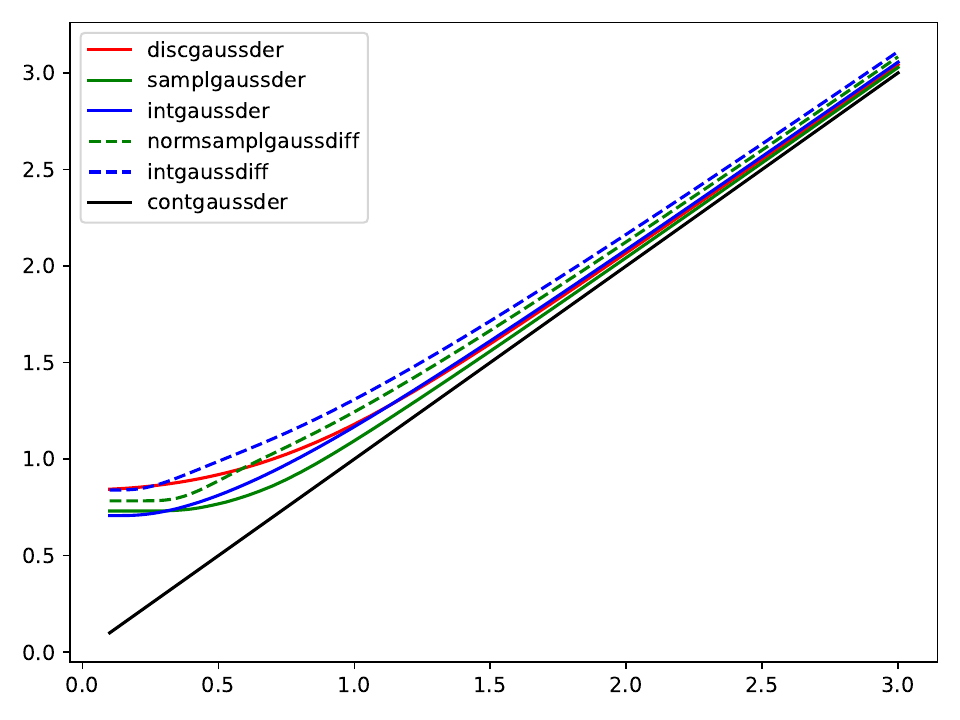}
        \end{tabular}
  \caption{Selected scales \, $\hat{\sigma}$.}
  %\caption{Selected scales $\hat{\sigma} = \sqrt{\hat{s}}$ from
  %$\argmax L_{v,\norm}$.}
    \end{subfigure}

    \bigskip
    \bigskip    

    \begin{subfigure}[t]{0.45\textwidth}
      \begin{tabular}{c}
          {\em Relative scale error for Laplacian scale selection\/} \\
          \includegraphics[width=0.97\textwidth]{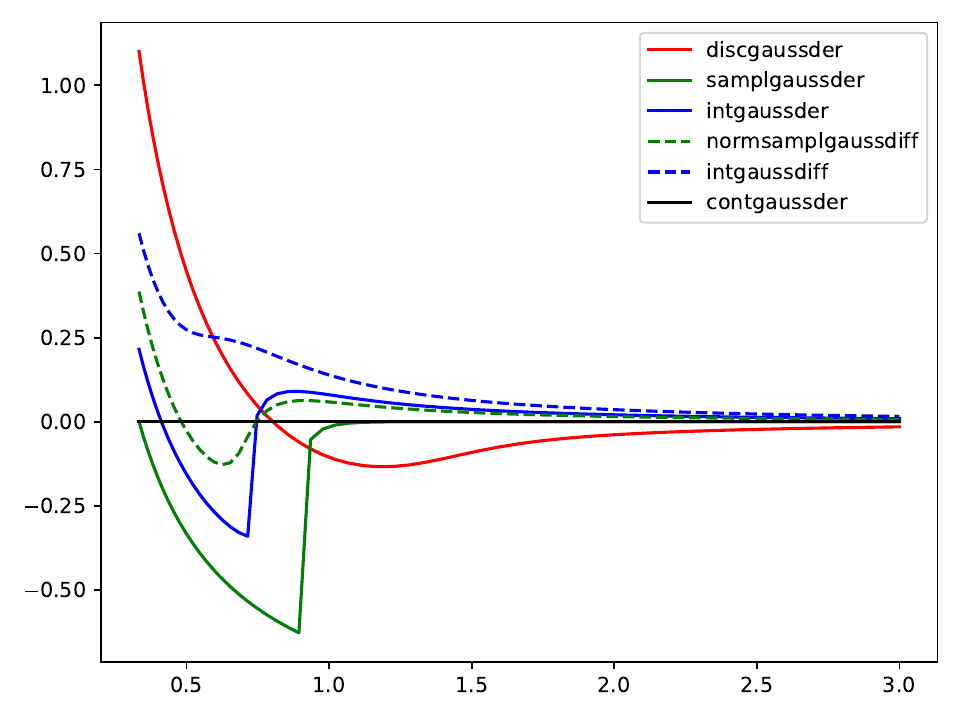}
      \end{tabular}
      \caption{Relative scale estimation error \, $E_{\scaleestrel}(\sigma)$.}
     %\caption{Relative error
    %$E_{\scaleestrel}(\sigma)$ from $\hat{\sigma} = \argmin \nabla_{\norm}^2 L$.}
        \label{fig-rel-sc-err-Lapl-scsel-gauss-blob}
    \end{subfigure}%
    ~\quad\quad~%
    \begin{subfigure}[t]{0.45\textwidth}
       \begin{tabular}{c}
          {\em Relative scale error for gradient magnitude scale selection\/} \\
           \includegraphics[width=0.97\textwidth]{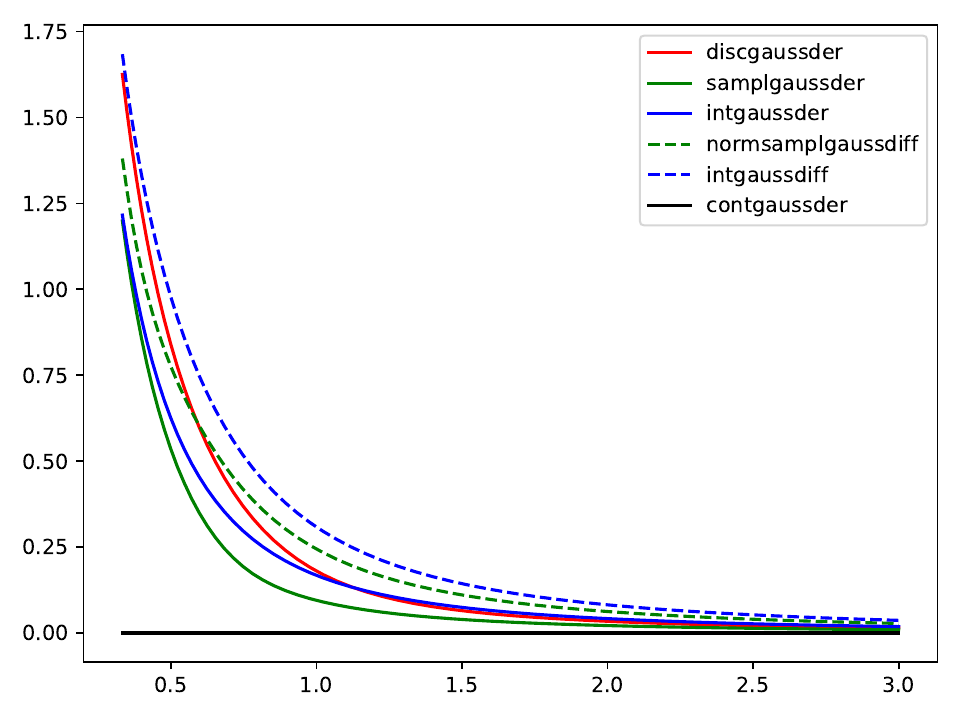}
    \end{tabular}
    \caption{Relative scale estimation error \, $E_{\scaleestrel}(\sigma)$.}
     %\caption{Relative error
     % $E_{\scaleestrel}(\sigma)$ from $\hat{\sigma} = \argmax L_{v,\norm}$.}
    \label{fig-rel-sc-err-gradmagn-scsel-diff-edge}
     \end{subfigure}
  \end{center}
  
  \caption{Graphs of the {\em selected scales\/} $\hat{\sigma} =
    \sqrt{\hat{s}}$ as well as the {\em relative scale estimation error\/}
    $E_{\scaleestrel}(\sigma)$, according to (\ref{eq-def-sc-sel-rel-sc-err}),
    when {\bf (left column)} applying scale selection from local extrema over scale of the
    {\em scale-normalised Laplacian\/} response according to (\ref{eq-sc-norm-lapl})
    to a set of Gaussian blobs of different size $\sigma_{\scaleref} =  \sigma_0$,
    for  different discrete approximations of the Gaussian
    derivative kernels or {\bf (right column)}  when applying scale selection from local extrema over scale of the
    {\em scale-normalised gradient magnitude\/} response according to (\ref{eq-sc-norm-grad-magn})
    to a set of diffuse step edges of different width $\sigma_{\scaleref} =  \sigma_0$,
    for either {\bf (i)}~discrete analogues of Gaussian
    derivative kernels $T_{\disc,x^{\alpha}}(n;\; s)$
    according to (\ref{eq-disc-der-gauss}), {\bf (ii)}~sampled Gaussian
    derivative kernels $T_{\sampl,x^{\alpha}}(n;\, s)$ according to
    (\ref{eq-sampl-gauss-der}), {\bf (iii)}~integrated Gaussian derivative
    kernels $T_{\intdisc,x^{\alpha}}(n;\, s)$ according to
    (\ref{eq-def-int-gauss-der}), {\bf (iv)}~the hybrid discretisation method
    corresponding the equivalent convolution kernels
    $T_{\hybrnormsampl,x^{\alpha}}(n;\; s)$ according to
    (\ref{eq-hybr-normsampl-disc-der}) or
    {\bf (v)}~the hybrid discretisation method
    corresponding the equivalent convolution kernels
    $T_{\hybrint,x^{\alpha}}(n;\; s)$ according to
    (\ref{eq-hybr-int-disc-der}).
    ({\bf Horizontal axes:} Reference scale
    $\sigma_{\scaleref} = \sigma_0 \in [1/3, 3]$.)}
  \label{fig-rel-sc-err-Lapl-edge-scsel-gauss-blob}
\end{figure*}

\begin{figure*}[hbtp]

  \begin{center}
     \begin{subfigure}[t]{0.45\textwidth}
        \begin{tabular}{c}
           {\em Selected scales for detHessian scale selection \/} \\
           \includegraphics[width=0.97\textwidth]{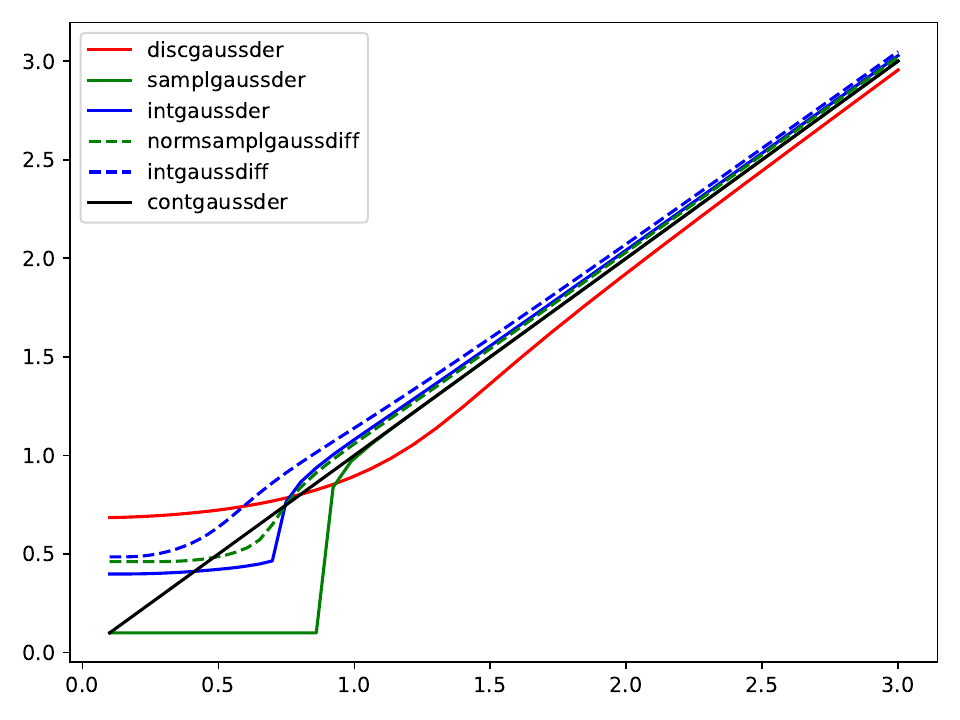}
        \end{tabular}
        \caption{Selected scales \, $\hat{\sigma}$.}
        %\caption{Selected scales $\hat{\sigma} = \sqrt{\hat{s}}$ from
         % $\argmax \det {\mathcal H} L_{\norm}$.}
         \label{fig-sel-sc-detHess-scsel-gauss-blob}
     \end{subfigure}%
     ~\quad\quad~
    \begin{subfigure}[t]{0.45\textwidth}
       \begin{tabular}{c}
            {\em Selected scales for principal curvature scale selection \/} \\
            \includegraphics[width=0.97\textwidth]{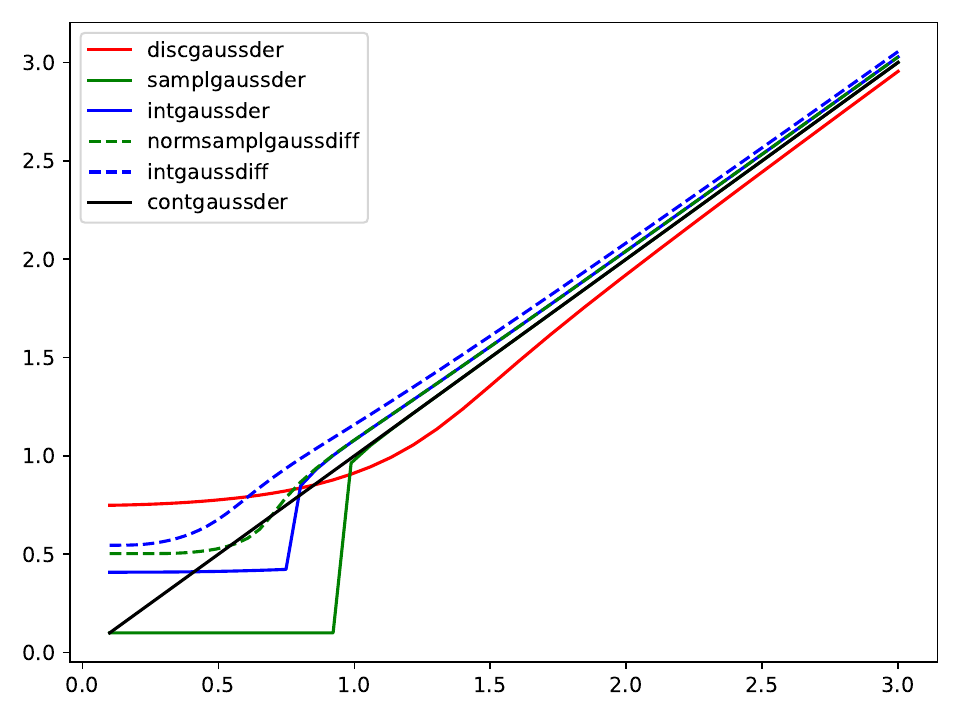}
         \end{tabular}
  \caption{Selected scales \, $\hat{\sigma}$.}
  %\caption{Selected scales $\hat{\sigma} = \sqrt{\hat{s}}$ from
  %$\argmin L_{pp,\norm}$.}
    \end{subfigure}

    \bigskip
    \bigskip    

    \begin{subfigure}[t]{0.45\textwidth}
      \begin{tabular}{c}
          {\em Relative scale error for detHessian selection\/} \\
          \includegraphics[width=0.97\textwidth]{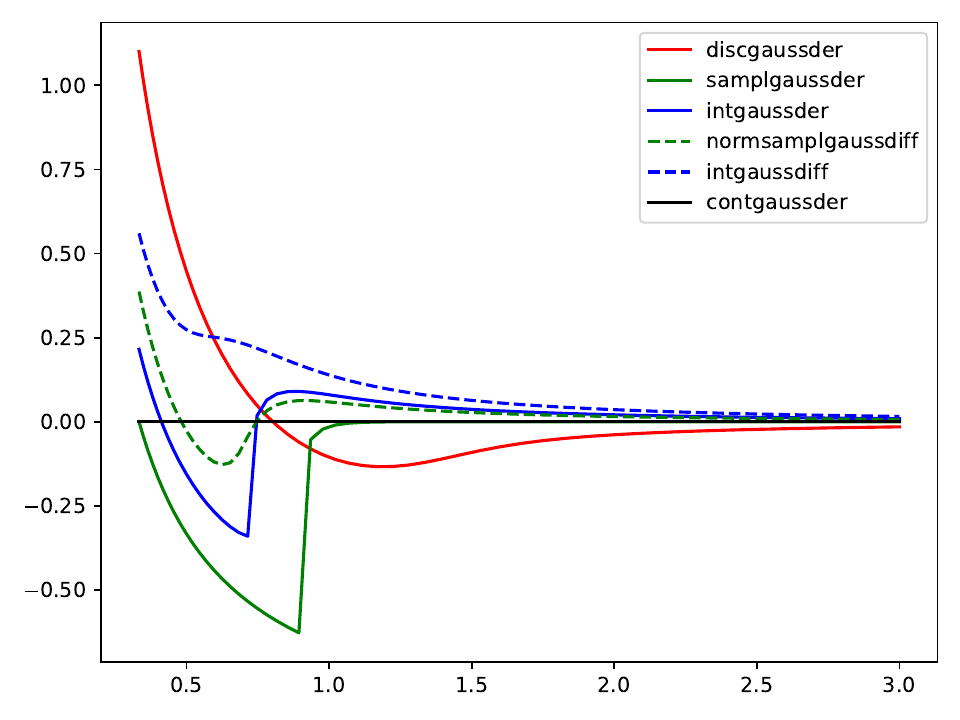}
      \end{tabular}
      \caption{Relative scale estimation error \, $E_{\scaleestrel}(\sigma)$.}
     %  \caption{Relative error
    %$E_{\scaleestrel}(\sigma)$ from $\hat{\sigma} = \argmax {\cal H}_{\norm} L$.}
       \label{fig-rel-sc-err-detHess-scsel-gauss-blob}
    \end{subfigure}%
    ~\quad\quad~%
    \begin{subfigure}[t]{0.45\textwidth}
       \begin{tabular}{c}
          {\em Relative scale error for principal curvature scale selection\/} \\
           \includegraphics[width=0.97\textwidth]{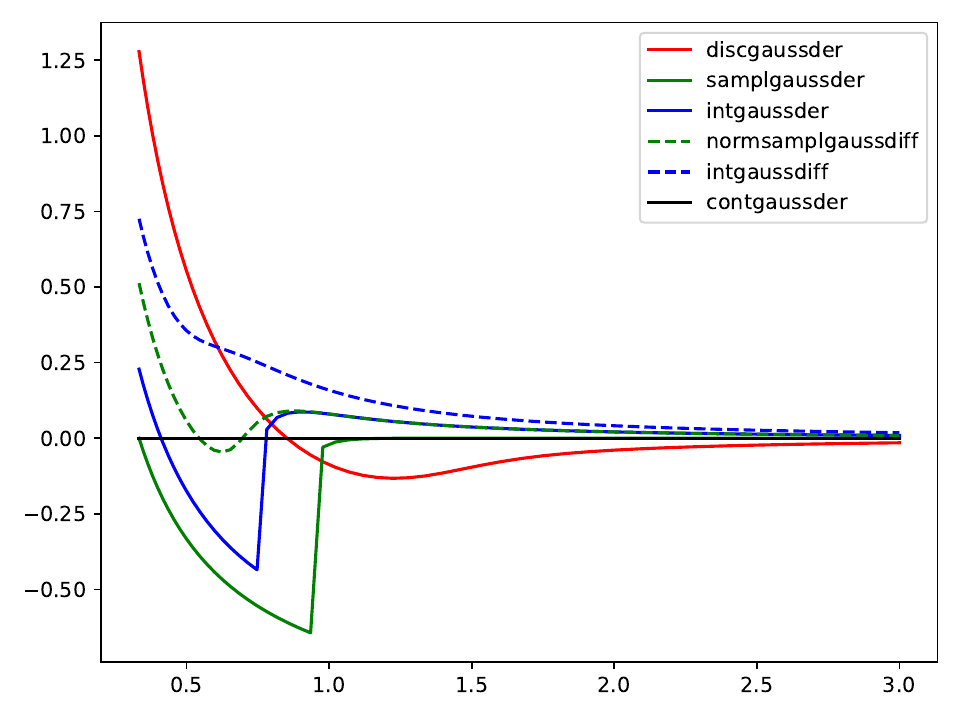}
    \end{tabular}
      \caption{Relative scale estimation error \, $E_{\scaleestrel}(\sigma)$.}
     % \caption{Relative error
     %  $E_{\scaleestrel}(\sigma)$ from $\hat{\sigma} = \argmin L_{pp,\norm}$.}
    \label{fig-rel-sc-err-gradmagn-scsel-diff-edge}
     \end{subfigure}
  \end{center}
  
  \caption{Graphs of the {\em selected scales\/} $\hat{\sigma} =
    \sqrt{\hat{s}}$ as well as the {\em relative scale estimation error\/}
    $E_{\scaleestrel}(\sigma)$, according to (34) %(\ref{eq-def-sc-sel-rel-sc-err}),
    when {\bf (left column)} applying scale selection from local extrema over scale of the
   {\em scale-normalised determinant of Hessian\/} response according
   to (22) %(\ref{eq-sc-norm-det-hess})
    to a set of Gaussian blobs of different size $\sigma_{\scaleref} =  \sigma_0$,
    for  different discrete approximations of the Gaussian
    derivative kernels or {\bf (right column)}  when applying scale selection from local extrema over scale of the
    {\em scale-normalised principal curvature\/} response according to
    (30) %(\ref{eq-sc-norm-ridge-meas})
    to a set of Gaussian ridges of different width $\sigma_{\scaleref} =  \sigma_0$,
    for either {\bf (i)} discrete analogues of Gaussian
    derivative kernels $T_{\disc,x^{\alpha}}(n;\; s)$
    according to (10), % (\ref{eq-disc-der-gauss}),
    {\bf (ii)} sampled Gaussian
    derivative kernels $T_{\sampl,x^{\alpha}}(n;\, s)$ according to
    (8), %(\ref{eq-sampl-gauss-der}),
    {\bf (iii)} integrated Gaussian derivative
    kernels $T_{\intdisc,x^{\alpha}}(n;\, s)$ according to
    (9), %(\ref{eq-def-int-gauss-der}),
    {\bf (iv)} the hybrid discretisation method
    corresponding the equivalent convolution kernels
    $T_{\hybrnormsampl,x^{\alpha}}(n;\; s)$ according to
    (16) %(\ref{eq-hybr-normsampl-disc-der})
    or
    {\bf (v)} the hybrid discretisation method
    corresponding the equivalent convolution kernels
    $T_{\hybrint,x^{\alpha}}(n;\; s)$ according to
    (17). %(\ref{eq-hybr-int-disc-der}).
    ({\bf Horizontal axes:} Reference scale
    $\sigma_{\scaleref} = \sigma_0 \in [1/3, 3]$.)}
  \label{fig-rel-sc-err-dethess-princcurv-scsel-gauss-blob}
\end{figure*}

Figure~\ref{fig-rel-sc-err-Lapl-edge-scsel-gauss-blob}(a) shows a graph of the
scale estimates obtained for the second-order
Laplacian interest point detector in this way,
with the corresponding relative scale errors in
Figure~\ref{fig-rel-sc-err-Lapl-edge-scsel-gauss-blob}(c).
Figures~\ref{fig-rel-sc-err-dethess-princcurv-scsel-gauss-blob}(a) and \ref{fig-rel-sc-err-dethess-princcurv-scsel-gauss-blob}(c) 
show corresponding results for the non-linear determinant
of the Hessian interest point detector.

As can be seen from these graphs, the consistency errors in
the scale estimates obtained for the hybrid discretisation methods,
based on either the normalised Gaussian kernel or the
integrated Gaussian kernel with central differences,
are for larger values of the scale
parameter higher than the corresponding consistency
errors in the regular discretisation methods based on either
sampled Gaussian derivatives or integrated Gaussian derivatives.
For smaller values of the scale parameter, there is, however,
a range of scale values, where the consistency errors are lower
for the hybrid discretisation methods than for underlying
corresponding regular discretisation methods.

Notably, the consistency errors for the discretisation methods
involving central differences are also
generally lower for the genuinely discrete method, based on convolution
with the discrete analogue of the Gaussian kernel followed by central
differences, than for the hybrid methods.
%In this regard, it should, however, be remarked that the scale normalisation
%operation for the discrete analogues of Gaussian derivatives has here been
%performed based on a fully continuous model, while one could more generally
%consider deriving genuinely discrete scale normalisation factors for the
%discretisation approach based on discrete analogues of Gaussian derivatives.

Figure~\ref{fig-rel-sc-err-Lapl-edge-scsel-gauss-blob}(b) shows the selected
scale levels for the first-order gradient-magnitude-based edge detection operation,
with the corresponding relative error measures shown in
Figure~\ref{fig-rel-sc-err-Lapl-edge-scsel-gauss-blob}(d).
As can be seen from these graphs, the consistency errors are notably
higher for the hybrid discretisation approaches, compared to
their underlying regular methods.
In these experiments, the consistency errors are also higher for the
hybrid discretisation methods than for the genuinely discrete approach,
based on discrete analogues of Gaussian derivatives.

Finally, Figures~\ref{fig-rel-sc-err-dethess-princcurv-scsel-gauss-blob}(b) and~\ref{fig-rel-sc-err-dethess-princcurv-scsel-gauss-blob}(d) 
show corresponding results
for the second-order principal-curvature-based ridge detector,
which are structurally similar to the previous results for
the second-order Laplacian and determinant of the Hessian interest point detectors.

\section{Summary and discussion}

In this paper, we have extended the in-depth treatment of different discretisations
of Gaussian derivative operators in terms of explicit convolution operations
in (Lindeberg \citeyear{Lin24-JMIV}) to two more
discretisation methods, based on hybrid combinations of either convolutions
with normalised sampled Gaussian kernels or convolutions with integrated
Gaussian kernels with central difference operators.

The results from the treatment show that it is possible to characterise
general properties of these hybrid discretisation methods in terms of
the effective amount of spatial smoothing that they imply.
Specifically, for very small values of the scale parameter, the
results obtained after the spatial discretisation may differ significantly
from the results obtained from the fully continuous scale-space theory,
as well as between the different types of discretisation methods.

The results from this treatment are intended to be generically applicable
in situations, when scale-space operations are to be applied at scale
levels below the otherwise rule of thumb in classical computer vision,
of not going below a certain
minimum scale level, corresponding to a standard deviation of the
Gaussian kernel of the order of $1/\sqrt{2}$ or 1.
We argue that the presented quantitative performance characterisations
should have a predictive ability, for how the different
types of discrete derivative approximation methods could have
comparative advantages in other multi-scale settings.

One direct application domain for these results is when
implementing deep networks in terms of Gaussian derivatives,
where empirical evidence indicates that deep networks often
tend to benefit from using finer scale levels than as indicated
by the previous rule of thumb in classical computer vision,
and which we will address in future work.

\section*{Acknowledgments}

Python code, that implements the discretisation methods for
Gaussian smoothing and Gaussian derivatives in this paper, is available
in the pyscsp package, available at GitHub:
\begin{quote}
  \url{https://github.com/tonylindeberg/pyscsp}
\end{quote}
as well as through PyPi:
\begin{quote}
  \tt pip install pyscsp
\end{quote}
\appendix

\section*{Appendix}

\section{Theoretical background of the notion of scale-space
  representation}
\label{app-theor-scsp}

Concerning the notion of multi-scale signal and image representations,
the concept of scale-space representation stands out as
a canonical type of multi-scale representation, in the respect that its
basic image operators in terms of Gaussian smoothing and Gaussian
derivatives can be {\em uniquely\/} determined from symmetry properties
that reflect structural properties of the environment in combination
with requirements to guarantee consistent treatment of image structures
between different spatial scales
(Iijima \citeyear{Iij62}; Koenderink \citeyear{Koe84};
Koenderink and van Doorn \citeyear{KoeDoo92-PAMI};
Lindeberg \citeyear{Lin96-ScSp}, \citeyear{Lin10-JMIV}, \citeyear{Lin21-Heliyon};
Weickert {\em et al.\/}\ \citeyear{WeiIshImi99-JMIV}).

This appendix section develops the basic theoretical background
concerning scale-space representation, regarding:
\begin{itemize}
\item
  the definitions of a scale-space representation
  with its Gaussian derivative responses,
\item
  the fundamental scale
  covariance property of the scale-normalised Gaussian derivative
  responses and homogeneous polynomial differential invariants
  constructed from scale-normalised Gaussian derivative responses,
\end{itemize}

which are then used in Sections~\ref{sec-meth-consist-scales} and~\ref{sec-results-scsel} for:
\begin{itemize}
\item
  quantifying the
  performance of different discrete approximations of Gaussian in terms
  of scale selection properties for different types of differential
  feature detectors.
\end{itemize}
Compared to previous presentations of these topics, mainly
devoted to 2-D images, and which this
presentation builds upon, the treatment in this section is performed
for general $D$-dimensional images.
Note specifically that in order to handle the general $D$-dimensional
case in this supplementary appendix,
we throughout refer denote the $D$-dimensional
image coordinates as $x = (x_1, x_2, \dots, x_D)^T$
as opposed to as the use of $(x, y)$ for the 2-D image coordinates in Section~\ref{sec-meth-consist-scales}.

\subsection{Scale-space representation and Gaussian derivative responses}

Given a $D$-dimensional signal $f \colon \bbbr \rightarrow \bbbr$,
a scale-space representation $L \colon \bbbr \times \bbbr_+ \rightarrow \bbbr$
is the one-parameter family of signals obtained by convolving the
input signal
\begin{equation}
  \label{eq-def-raw-scsp-D-dims}
  L(x;\; s)
  = (g(\cdot;\; s) * f(\cdot))(x;\; s)
  = \int_{u \in \bbbr^D} f(x - u) \, g(u;\; s) \, du
\end{equation}
with $D$-dimensional Gaussian kernels at different scales $s \in \bbbr_+$
\begin{equation}
  \label{eq-gauss-kern-D-dim}
  g(x;\; s) = \frac{1}{\sqrt{2\pi s}^D} \, e^{-\frac{x^T x}{2s}},
\end{equation}
where $x = (x_1, x_2, \dots, x_D)^T$ % in this supplementary material
denotes $D$-dimensional image coordinates.

From this representation, Gaussian derivative responses of order
$\alpha$, where $\alpha = ({\alpha_1, \alpha_2, \dots, \alpha_D})$
is a multi-index notation, are, in turn, defined as
\begin{equation}
  L_{x^{\alpha}}(x;\; s)
  = \partial_{x^{\alpha}} L(x;\; s)
  =  (g_{x^{\alpha}}(\cdot;\; s) * f(\cdot))(x;\; s),
\end{equation}
where
\begin{equation}
  g_{x^{\alpha}}(x;\; s)
  = \partial_{x^{\alpha}} (g(x;\; s))
  = \partial_{x_1^{\alpha_1} x_2^{\alpha_2} \dots x_D^{\alpha_D}} (g(x;\; s))
\end{equation}
denotes the $D$-dimensional Gaussian derivative kernel of order
$\alpha$.
In the area of scale-space theory for image processing and
computer vision, it has been shown that such
Gaussian derivative responses can be used as a powerful basis for
expressing a large number of visual operations, including feature detection,
image matching, classification and object recognition,
see Lindeberg (\citeyear{Lin94-SI}, \citeyear{Lin99-CVHB},
\citeyear{Lin08-EncCompSci}) for overviews.

\subsection{Scale covariance property of Gaussian derivative
  responses}

With regard to handling image structures at different spatial scales,
a very useful property of this notion of scale-space representation is
that it is closed under spatial scaling transformations of the form
\begin{equation}
  f'(x') = f(x) \quad\quad\mbox{for}\quad\quad x' = S \, x,
\end{equation}
where $S > 0$ is a spatial scaling factor. Specifically, if we define
scale-normalised Gaussian derivative responses according to
(Lindeberg \citeyear{Lin97-IJCV} Equation~(18))
\begin{equation}
   L_{\xi^{\alpha}}(x;\; s) = s^{\gamma \, |\alpha|/2} \, L_{x^{\alpha}}(x;\; s),
\end{equation}
where $\gamma > 0$ is a scale normalisation power and
$|\alpha| = \alpha_1 + \alpha_2 + \dots + \alpha_D$
denotes the total order of differentiation, then it follows that if we
additionally define the scale-normalised Gaussian derivative response
of the rescaled signal $f'(x')$ according to
\begin{equation}
  L'_{x^{\alpha}}(x';\; s')
%  = \partial_{x^{\alpha}} L'(x;\; s)
  =  ({s'}^{\gamma \, |\alpha|/2} g_{x^{\alpha}}(\cdot;\; s') * f'(\cdot))(x';\; s'),
\end{equation}
then the corresponding scale-normalised responses between
matching image positions $x$ and $x'$ in the two domains
will be equal up to a uniform scaling factor 
(see Lindeberg (\citeyear{Lin97-IJCV}) Equation~(20), although with different
notation for the scale parameters and the scaling factor used there) 
\begin{equation}
  L_{x^{\alpha}}(x;\; s) = S^{|\alpha| \, (1 - \gamma)} L'_{x^{\alpha}}(x';\; s'),
\end{equation}
provided that the values of the scale parameters $s$ and $s'$ in the
two domains are matched according to
\begin{equation}
  \label{eq-sc-cov-prop-s-transf}
  s' = S^2 \, s.
\end{equation}
This scale covariance property of Gaussian derivative responses has
been of fundamental importance to classical computer vision, for
formulating scale-invariant algorithms that can handle {\em a  priori\/}
unknown scaling variabilities in image data in a fully automatic
manner
(Lindeberg \citeyear{Lin97-IJCV}, \citeyear{Lin98-IJCV},
\citeyear{Lin12-JMIV}, \citeyear{Lin15-JMIV}, \citeyear{Lin21-EncCompVis};
Bretzner and Lindeberg \citeyear{BL97-CVIU};
Chomat {\em et al.\/} \citeyear{ChoVerHalCro00-ECCV};
Lowe \citeyear{Low04-IJCV};
Bay {\em et al.\/} \citeyear{BayEssTuyGoo08-CVIU}).
When using Gaussian derivative responses as mathematical primitives for
formulating the layers in deep networks, this scale covariance property has also been
of fundamental importance for formulating provably scale-covariant
and scale-invariant deep networks
(Lindeberg \citeyear{Lin20-JMIV}, \citeyear{Lin22-JMIV};
Sangalli {\em et al.\/} \citeyear{SanBluVelAng22-BMVC};
Yang {\em et al.\/} \citeyear{YanDasMah23-arXiv}).

Specifically, with regard to feature detection methods formulated in
terms of homogeneous polynomial combinations of such scale-normalised
Gaussian derivative responses (Lindeberg \citeyear{Lin97-IJCV} Equation~(21))
\begin{equation}
  \label{eq-def-hom-pol-diff-inv}
  ({\mathcal D}_{\norm} L)(x;\; s)
  = \sum_{i = 1}^I c_i \prod_{j = 1}^J L_{\xi^{\alpha_{ij}}}(x;\; s)
\end{equation}
where the sum of the orders of differentiation in a certain term
\begin{equation}
  \sum_{j = 1}^J |\alpha_{ij}| = M
\end{equation}
is the same for all the terms in (\ref{eq-def-hom-pol-diff-inv}) with
indices $i$, then the corresponding transformed feature
\begin{equation}
  ({\mathcal D}_{\norm} L')(x';\; s')
  = \sum_{i = 1}^I c_i \prod_{j = 1}^J L'_{\xi^{\alpha_{ij}}}(x';\; s')
\end{equation}
will for matching spatial positions $x$ and $x'$ be equal
up to a uniform scaling factor 
(see Lindeberg (\citeyear{Lin97-IJCV}) Equation~(25), although again
with different notation of the spatial scale parameters and the
spatial scaling factor used there)
\begin{equation}
  ({\mathcal D}_{\norm} L)(x;\; s) = S^{M (1 - \gamma)} ({\mathcal D}_{\norm} L')(x';\; s'),
\end{equation}
provided that the scale parameters $s$ and $s'$ are matched according
to (\ref{eq-sc-cov-prop-s-transf}).

Examples of such scale-normalised differential feature detectors for
the case of 2-D spatial images are
given in Equations~(\ref{eq-sc-norm-lapl}) and~(\ref{eq-sc-norm-det-hess}) regarding interest
point detection, in Equation~(\ref{eq-sc-norm-grad-magn}) regarding edge
detection and in Equation~(\ref{eq-sc-norm-ridge-meas}) regarding ridge detection.

\subsection{Evaluating discretisation methods based on consistency
  requirement relative to scale selection properties}

In the experiments reported in Section~\ref{sec-results}, we use
consistency measurements relative
to this transformation property of Gaussian-derivative-based feature
detectors under spatial scaling transformations, to
quantitatively evaluate to quality of different discrete approximations of
Gaussian derivative kernels, that are used for approximating Gaussian
derivative responses in discrete implementations for 2-D image data.

Specifically, we will evaluate feature detectors in terms of interest
point detection, edge detection and ridge detection, formulated
with the systematic scale selection methodology proposed in
(Lindeberg \citeyear{Lin97-IJCV}, \citeyear{Lin98-IJCV}),
based on detecting local extrema over scales of scale-normalised
differential entities of the form (\ref{eq-def-hom-pol-diff-inv})
according to
\begin{equation}
  \{ \hat{s}(x) \}
  = \operatorname{argmaxminlocal}_s   ({\mathcal D}_{\norm} L)(x;\; s),
\end{equation}
where $\operatorname{argmaxminlocal}_s ({\mathcal D}_{\norm} L)(x;\; s)$
denotes the set of all local extrema over scale of the differential
entity ${\mathcal D}_{\norm} L$ at the spatial point $x$.
In the ideal continuous case, it holds that these scale
levels $\hat{s}$, where local extrema
of differential invariants of the form (\ref{eq-def-hom-pol-diff-inv})
assume local maxima over scales,
will transform according to (\ref{eq-sc-cov-prop-s-transf})
\begin{equation}
  \label{eq-sc-cov-prop-s-transf-scsel}
  \hat{s'} = S^2 \, \hat{s},
\end{equation}
and are in this respect provably scale covariant.

Due to deviations between the continuous theory and the
discrete implementation, based on spatially
discretised input functions combined with different spatial discretisations of
the Gaussian derivative kernels, we cannot, however, expect
corresponding properties to hold exactly in a numerical
implementation. In this respect, the deviations between the results
obtained from a discrete implementation of these computational steps
in relation to the corresponding fully continuous results can be used
for evaluating the approximation properties of different approaches
for discretising Gaussian derivative operators.

With such a paradigm for quantitative evaluation, we can, in
particular, quantify the accuracy of different spatial
discretisations relative to a primary desirable property of a
multi-scale representation.

In methods for feature detection with automatic scale selection,
it is specifically assumed that the scale levels $\hat{s}$, where
appropriately designed differential feature detectors assume their
maxima over scales, may correspond to interesting structures in the
original image data. By experimental investigations for different types
of feature detectors, for tasks such as edge detection,
blob detection, corner detection, interest point detection and
dense local scale estimation, it has been demonstrated that this approach
allows for robust feature detectors that very well handle substantial
variabilities in scale for natural image data
(Lindeberg \citeyear{Lin97-IJCV}, \citeyear{Lin98-IJCV}, \citeyear{Lin99-CVHB},
\citeyear{Lin12-JMIV}, \citeyear{Lin18-SIIMS}, \citeyear{Lin21-EncCompVis};
Lowe \citeyear{Low04-IJCV}).

\subsubsection{The notion of characteristic scale}
\label{app-char-scale}

In the area of computer vision, such
locally dominant scale values in the image data, where
scale-normalised differential entities assume local extrema over
scale, are often referred to as locally
``characteristic scales'', based on the empirically verified
hypothesis that such scale levels, when measured in terms of the scale
parameter expressed in units dimension of $[\mbox{length}]$ as
$\sigma = \sqrt{s}$, will for both idealised models of image structures
and natural image data with similar appearance reflect a
characteristic length of such image structures.
See the introduction of Section~4 in (Lindeberg \citeyear{Lin97-IJCV})
for a definition of the notion of
characteristic length of features extracted image data and more generally
Sections~5.1 and~6.2 in (Lindeberg \citeyear{Lin97-IJCV}),
Sections~4.5 and~5.6.1 in (Lindeberg \citeyear{Lin98-IJCV}),
Sections~3.1.1 and~3.2--3.3 in (Lindeberg \citeyear{Lin12-JMIV})
and Sections~2.2--2.3 in (Lindeberg \citeyear{Lin18-SIIMS})
for examples of mathematical analysis of scale selection properties for different
types of idealised model signals.

\section{Theoretical background for
  discretisations of the Gaussian smoothing operation and
  the Gaussian derivative operators}
\label{app-theor-disc}

This section gives mathematical definitions of the different
discretisation methods, that we will consider for 
discrete approximations of Gaussian derivative operators and then
quantitatively evaluate in Section~\ref{sec-results}. For simplicity, we
throughout restrict ourselves to convolution operations that are separable over
the $D$-dimensional image data, thereby making it possible to restrict
the following theoretical treatment to studying 1-D discrete approximations of Gaussian
derivatives, based on the separability property of the Gaussian
derivative operators in the continuous case described below.

The theoretical material in this appendix is based on an in-depth treatment of the
topic of approximating the Gaussian convolution and the Gaussian
derivative operators for discrete data in (Lindeberg
\citeyear{Lin24-JMIV}). For related treatments of the topic of
approximating scale-space operators on discrete data, see also
(Lindeberg \citeyear{Lin90-PAMI}, \citeyear{Lin93-JMIV},
\citeyear{Lin93-Dis},
{\AA}str{\"o}m and Heyden \citeyear{AstHey97-SCSPTH},
Wang \citeyear{Wan00-TIP},
Lim and Stiehl \citeyear{LimSti03-ScSp},
Tschirsich and Kuijper \citeyear{TscKui15-JMIV},
Slav{\'\i}k and Stehl{\'\i}k \citeyear{SlaSte15-JMathAnalAppl},
Rey-Otero and Delbracio \citeyear{OteDel16-IPOL}).

\subsection{Separable Gaussian smoothing and Gaussian derivative computations}
\label{app-sep-gauss-smooth-gauss-ders}

Due to the separability of the Gaussian kernel in Equation~(\ref{eq-gauss-kern-D-dim})
\begin{equation}
  g(x;\; s) = \prod_{i = 1}^D g_{1D}(x_i;\; s),
\end{equation}
where $g_{1D}(x_i;\; s)$ denotes a 1-D Gaussian kernel over the $i$:th
dimension with $x = (x_1, %x_2,
\dots, x_D$)
\begin{equation}
  \label{eq-def-1D-gauss-kern-app}
  g_{1D}(x_i;\; s) = \frac{1}{\sqrt{2 \pi s}} \, e^{-\frac{x_i^2}{2s}},
\end{equation}
it follows that the $D$-dimensional convolution integral for computing
the $D$-dimensional scale-space representation according to
(\ref{eq-def-raw-scsp-D-dims}) can be computed as a cascade of
separable Gaussian convolution steps
\begin{multline}
  L(x;\; s) %= \\
  = \int_{u_1\in \bbbr} g_{1D}(u_1;\, s) \, \dots \\ \left( \int_{u_D
       \in \bbbr} g_{1D}(u_D;\; s) \, f(x_1 - u_1, \dots, x_D - u_D) \,
     du_D \right) \dots \, du_1.
\end{multline}
In a corresponding way, from the separability of the corresponding Gaussian
convolution kernels with $\alpha = (\alpha_1, %\alpha_2,
\dots, \alpha_D)$
\begin{equation}
  g_{x^{\alpha}}(x;\; s) = \prod_{i = 1}^D g_{1D,x^{\alpha_i}}(x_i;\; s),
\end{equation}
where each 1-D Gaussian derivative kernel is given by
\begin{equation}
  g_{1D,x^{\alpha_i}}(x_i;\; s) = \partial_{x_i^{\alpha_i}}(g_{1D}(x_i;\; s)),
\end{equation}
it follows that the $D$-dimensional Gaussian derivative kernel of order
$\alpha$ can be computed with a separable cascade with a set of 1-D
Gaussian derivative kernels of order $\alpha_i$ along each of the
respective dimensions
according to
\begin{align}
   \begin{split}
      & L_{x^{\alpha}}(x;\; s) %= \\
   \end{split}\nonumber\\
   \begin{split}
     &
      = \int_{u_1\in \bbbr} g_{1D,x^{\alpha_1}}(u_1;\, s) %\times
      \dots
   \end{split}\nonumber\\
   \begin{split}
       & \quad\quad\quad\quad
      \left( \int_{u_D \in \bbbr}
        g_{1D,x^{\alpha_D}}(u_D;\; s) \,
      \right.
  \end{split}\nonumber\\
   \begin{split}
     & \quad\quad\quad\quad\quad\quad\quad\quad
     \left.\vphantom{\int_{u_D \in \bbbr}}
          f(x_1 - u_1, \dots, x_D - u_D) \, du_D
             \right) \dots \, du_1.
    \end{split}
\end{align}

\subsection{Definition of the 1-D modelling problems studied regarding
  discrete approximations of Gaussian smoothing and Gaussian
  derivatives}
\label{app-def-1D-model-probl}

In the following, we will henceforth%
\footnote{In fact, for the principled axiomatic theory of discrete
  scale-space representations in (Lindeberg \citeyear{Lin90-PAMI}, \citeyear{Lin93-Dis}),
  it has been shown that more accurate approximations to rotational
  symmetry can be obtained by instead using discrete kernels that are not
  formally restricted to being computed with a cascade of only using
  one convolution kernel along each dimension
  (Lindeberg \citeyear{Lin93-Dis} Proposition~4.16).
  Since those more isotropic discretisations, however, require either the combination
  of two separable convolutions along the Cartesian coordinate
  directions with two separable diagonal convolutions in the 2-D case, or a
  Fourier-based implementation,
  and we are for our intended applications interested in
  achieving high computational efficiency with publicly available
  software libraries for discrete implementations based on purely discrete
  convolution operations, we will, however, not
  consider those discretisations in this work.}
consider spatial discretisations of 1-D Gaussian convolution integrals
of the form
\begin{equation}
  \label{eq-1D-gauss-conv-int}
  L(x;\; s) 
  = \int_{u \in \bbbr} g_{1D}(u;\, s) \, f(x - u) \, du
\end{equation}
with discrete approximations of Gaussian kernels
of the form
\begin{equation}
  \label{eq-1D-gauss-conv-disc}  
  L(x;\; s) 
  = \sum_{n \in \bbbz} T(n;\, s) \, f(x - n),
\end{equation}
as well as spatial discretisations of 1-D Gaussian derivative
convolution integrals of the form
\begin{equation}
  \label{eq-1D-gaussder-conv-int}  
  L_{x^{\alpha}}(x;\; s) 
  = \int_{u \in \bbbr} g_{1D,{x^{\alpha}}}(u;\, s) \, f(x - u) \, du
\end{equation}
with discrete
approximations of Gaussian derivative kernels of
the form
\begin{equation}
  \label{eq-1D-gaussder-conv-disc}    
  L_{x^{\alpha}}(x;\; s) 
  = \sum_{n \in \bbbr} T_{x^{\alpha}}(n;\, s) \, f(x - n),
\end{equation}
where we will then specifically consider different options for
choosing the discrete approximations of Gaussian kernels $T(n;\; s)$
and the discrete approximations of Gaussian derivative kernels
$T_{x^{\alpha}}(n;\; s)$, or mathematical equivalent formulations
thereof, based on different types of criteria.

\subsection{The sampled Gaussian kernel}
\label{app-sampl-gauss-kern}

The presumably most straightforward approach to approximate the
continuous Gaussian convolution integral (\ref{eq-1D-gauss-conv-int}) by a discrete
convolution operation of the form (\ref{eq-1D-gauss-conv-disc}) is by sampling the
continuous integral (\ref{eq-1D-gauss-conv-int}) with equidistant
samples, which then corresponds
to discrete convolution with the sampled Gaussian kernel
\begin{equation}
  \label{eq-sampl-gauss-app}
  T_{\sampl}(n;\; s) = g_{1D}(n;\; s)
\end{equation}
with the sampled Gaussian kernel $g_{1D}(n;\; s)$ defined from the
corresponding continuous kernel
$g_{1D}(x;\; s)$ according to (\ref{eq-def-1D-gauss-kern-app}).

While this operation can be expected to work very well for
sufficiently large values of the scale parameter $s$, this operation
may, however, imply potential problems for very small values of the scale
parameter. As developed further in
(Lindeberg \citeyear{Lin24-JMIV} Section~2.3),
for very small values of the scale parameter $s$, the sum of the discrete
filter coefficients may significantly exceed 1, which implies that the
corresponding discrete convolution operation then cannot be regarded as a
true spatial smoothing process, since it will not leave a constant
signal essentially unchanged. For this reason, it is therefore natural
to consider using other types of discretisation method for
discretising the Gaussian convolution integral at very fine scales.

\subsection{The normalised sampled Gaussian kernel}
\label{app-norm-sampl-gauss-kern}

A straightforward, but {\em ad hoc\/}, way of avoiding that the sum of
the filter coefficients may exceed 1 is by instead normalizing the
filter coefficient in the sampled Gaussian kernel to unit $l_1$-norm
\begin{equation}
  \label{eq-def-norm-sampl-gauss-app}
  T_{\normsampl}(n;\; s) = \frac{g(n;\; s)}{\sum_{m \in \bbbz} g(m;\; s)}.
\end{equation}
Then, we are sure that a constant input signal is left
unchanged. Still, however, a potential problem that arises when to
use the resulting normalised sampled Gaussian kernel for discrete
implementation is that the discrete variance of the resulting
discrete sampled Gaussian kernel may for small values of the scale parameter be
significantly lower than the continuous variance of the continuous
Gaussian kernel (see Lindeberg (\citeyear{Lin24-JMIV}) Section~2.4, specifically
Figures~3--4 in that paper for additional details). Since the
variance of a kernel is not affected by a
renormalisation of the filter coefficients, the same  structural problem
applies also to the regular sampled Gaussian kernel $T_{\sampl}(n;\; s)$.

\subsection{The integrated Gaussian kernel}
\label{app-int-gauss-kern}

A possibly less {\em ad hoc\/} way, of avoiding that the sum of the
discrete filter coefficients could deviate from their aimed ideal
value 1, is
by instead forming a discrete approximation of the Gaussian kernel by
integrating the values of the continuous Gaussian kernel over each
pixel support region
\begin{equation}
  \label{eq-def-int-gauss-kern-app}
   T_{\intdisc}(n;\; s) = \int_{x = n - 1/2}^{n + 1/2} g(x;\; s) \, dx.
 \end{equation}
Then, it clearly follows that the sum of the filter coefficient in the
resulting discrete kernel will be guaranteed to be equal to 1. As further described in
(Lindeberg \citeyear{Lin24-JMIV} Section~2.5 and Appendix~A.2),
the effect of this approach is
equivalent to first defining a continuous input function $f_c(x)$, that over each
pixel support region is equal to the integral of the original function
$f(x)$ over that pixel support region, and then subjecting that
continuous signal $f_c(x)$ to continuous Gaussian convolution.
As further described in
(Lindeberg \citeyear{Lin24-JMIV} Section~2.5), the spatial
integration operation in the definition of the derived continuous
signal $f_c(x)$ from the original input signal $f(x)$ does, however, introduce
a scale bias, in the sense that the discrete variance of the integrated
Gaussian kernel may be significantly larger than the continuous
variance of the continuous Gaussian kernel.

\subsection{The discrete analogue of the Gaussian kernel}
\label{app-disc-gauss}

A more principled approach to defining a discrete approximation of the
Gaussian kernel is by instead formulating similar assumptions
(referred to as scale-space axioms in the area of scale-space theory)
over a discrete domain as uniquely single out the choice of the
continuous Gaussian kernel over a continuous domain.

In (Lindeberg \citeyear{Lin90-PAMI,Lin93-Dis}) it is shown that, if we
assume that (i)~the discrete kernels should be guaranteed to not never
introduce new local extrema or equivalently new zero-crossings from
any finer to any coarser level of scale, then if combined with (ii)~a semi-group
property
\begin{equation}
  T(\cdot;\; s_1) * T(\cdot;\; s_2) = T(\cdot;\; s_1 + s_2)
\end{equation}
(iii)~spatial symmetry through the origin $T(-n;\; s) = T(-n;\; s)$,
and (iv)~normalisation to unit $l_1$ norm $\sum_{n} T(n;\; s) = 1$,
then the special family of kernels
\begin{equation}
  \label{eq-disc-gauss-app}
   T_{\disc}(n;\; s) = e^{-s} I_n(s),
\end{equation}
referred to as the discrete analogue of the Gaussian kernel is singled
out {\em uniquely\/}, where $I_n(s)$ denotes the modified Bessel functions of
integer order. This discrete kernel $T_{\disc}(n;\; s)$ is specifically the natural heat kernel over a 1-D
discrete domain, in the respect that it satisfies the semi-discrete
diffusion equation
\begin{equation}
  \partial_s L = \frac{1}{2} \, \delta_{xx} L 
 \end{equation}
with initial condition $L(n;\; 0) = f(n)$, where $\delta_{xx}$
denotes the second-order discrete difference operator
$\delta_{xx} = (+1, -2, +1)$.

With regard to handling 
multiple spatial scales in signals, an additional attractive property
of this discrete kernel is that the discrete variance of the discrete
analogue of the Gaussian kernel is exactly equal to the value of the
scale parameter
\begin{equation}
   V(T(\cdot;\; s)) = \sum_{n \in \bbbz} n^2 \, T(n;\; s) = s.
\end{equation}
In this respect, this discrete analogue of the Gaussian kernel
conceptually avoids the structural
problems with the sampled Gaussian kernel and the normalised sampled
Gaussian kernel at very fine scales, as described in
Sections~\ref{app-sampl-gauss-kern}
and~\ref{app-norm-sampl-gauss-kern} above.
Additionally, the discrete analogue of the Gaussian kernel also avoids
the scale offset of the integrated Gaussian kernel, as described in
Section~\ref{app-int-gauss-kern} above.

\subsection{The sampled Gaussian derivative kernel}
\label{app-sampl-gaussder-kern}

In analogy with the definition of the sampled Gaussian kernel
in Section~\ref{app-sampl-gauss-kern} above,
the presumably most straightforward approach of approximating the
Gaussian derivative convolution integral
(\ref{eq-1D-gaussder-conv-disc}) is by sampling that convolution
integral with equidistant samples, which then directly leads to a
convolution with the sampled Gaussian derivative kernel
\begin{equation}
  \label{eq-sampl-gauss-der-app}
  T_{\sampl,x^{\alpha}}(n;\, s) = g_{x^{\alpha}}(n;\, s).
\end{equation}
As for the previously mentioned problems regarding the sampled
Gaussian kernel, as described in Section~\ref{app-sampl-gauss-kern} above,
a conceptual problem with the sampled Gaussian derivative kernel,
however, is that that sampling
process may lead to substantial discretisation artefacts at very fine
scales, which then, as previously described in
(Lindeberg \citeyear{Lin24-JMIV} Section~3.3, specifically
Figures~10--11 in that paper),
implies that the discrete variance of
(the absolute value of) the sampled Gaussian derivative kernel may
vary significantly from the
continuous variance of (the absolute value of) the continuous Gaussian
kernel. Due to this property, the value of the scale parameter $s$ may
not fully reflect the amount of discrete blur in the corresponding
discrete derivative approximation kernel at very fine scales.

Another practical computational problem, that arises with the sampled
Gaussian derivative approach, is furthermore
that in use cases, where multiple derivatives of different orders
$\alpha$ are to be computed from the same input signal $f$,
then the computation of each discrete approximation of the Gaussian
derivative response for each order requires a new discrete convolution
operation with a full size discrete derivative approximation kernel.
Thereby, the computational work
increases linearly, with the same factor, for increasing numbers
of orders of differentiation, which can be problem in situations where
very large numbers of Gaussian derivative responses are to be
computed, such as in the application of Gaussian derivative operators
as computational primitives in deep networks.

For these reasons, it is therefore warranted to also consider other
types of discretisation methods for computing discrete approximations
of Gaussian derivative responses.

\subsection{The integrated Gaussian derivative kernel}
\label{app-int-gaussder-kern}

In analogy with the definition of the integrated Gaussian kernel in
Section~\ref{app-int-gauss-kern} above,
we can also with the potential aim of trying to reduce
the severe discretisation artefacts of the sampled Gaussian kernel at
very fine scales, define an integrated Gaussian kernel by integrating
the values of the continuous Gaussian derivative kernel over each
pixel support region
\begin{equation}
  \label{eq-def-int-gauss-der-app}
  T_{\intdisc,x^{\alpha}}(n;\, s) = \int_{x = n - 1/2}^{n + 1/2} g_{x^{\alpha}}(x;\; s) \, dx.
\end{equation}
As described in (Lindeberg \citeyear{Lin24-JMIV} Section~3.4 and
Appendix~A.2), the effect of this operation
is equivalent to first defining a continuous input signal $f_c(x)$,
by integrating the original input signal $f(x)$ over each pixel support
region and then letting the new continuous input signal $f_c(x)$
be equal to that average value over each pixel support region.
As further described in
(Lindeberg \citeyear{Lin24-JMIV} Section~3.4, specifically
Figures~10--11 in that paper),
a problem with this
approach is, however, is that the spatial integration operation, used for
defining the new continuous input signal $f_c(x)$ from the original
input signal $f(x)$, introduces a scale bias.

When to use the integrated Gaussian derivative kernel for computing
discrete approximations of multiple orders, this approach also suffers
from the same problem as the sampled Gaussian derivative kernel, as
described in Section~\ref{app-sampl-gaussder-kern} above,
in that a new large support convolution operation needs to
be invoked for each order of differentiation, thus implying that the
computational work will increase with the same linear factor for
increasing numbers of orders of differentiation, which may be a
substantial problem in situations when a very large number of Gaussian
derivative responses are to be computed, such as in deep networks
based on using combinations of Gaussian derivative operators as the
computational primitives.

\subsection{Discrete derivative approximations constructed from the
  combination of the discrete analogue of the Gaussian kernel with central
  differences}
\label{app-disc-gauss-ders}

According to the theory for discrete derivative approximations with
scale-space properties in (Lindeberg \citeyear{Lin93-JMIV,Lin93-Dis})
(see also Lindeberg (\citeyear{Lin24-JMIV}) Section~3.5), discrete
approximations of Gaussian derivative responses can be computed by applying
central difference operators $\delta_{x^{\alpha}}$ to the discrete
scale-space representation obtained by convolution with the discrete
analogue of the Gaussian kernel (\ref{eq-disc-gauss-app}).
This operation then leads to
equivalent discrete derivative approximation kernels of the form
\begin{equation}
  \label{eq-disc-der-gauss-app}
  T_{\disc,x^{\alpha}}(n;\; s) = (\delta_{x^{\alpha}} T_{\disc})(n;\; s),
\end{equation}
although with the important distinction that such equivalent discrete
approximation kernels should never be used for practical numerical
computations, but instead
just be considered as a conceptual modelling tool, for representing the
mathematically equivalent effect of the discrete derivative
approximation method.

The reason for this
is that in situations, when multiple Gaussian derivative responses are
to be computed from the same input signal $f$ for different orders of
differentiation $\alpha$, it is computationally much more efficient to instead
just perform the spatial smoothing operation with the discrete
analogue of the Gaussian kernel once and for all,
and then apply the small-support
central difference operators $\delta_{x^{\alpha}}$ for the different
orders of differentiation. In this way, a substantial amount of
computation can be saved in situations when several Gaussian
derivative responses are to be computed for different orders of
differentiation. This implementation approach thereby has substantial
computational advantages, when using discrete approximations of
Gaussian derivative responses as the basic computational primitives
for defining the layers in deep networks.

Other conceptual advantages of this discretisation method for Gaussian
derivative computation are also that the resulting discrete derivative
approximations obey a cascade smoothing property over scales
\begin{equation}
  L_{\disc}(\cdot;\; s_2) = T_{\disc}(\cdot;\; s_2 - s_1) * L_{\disc}(\cdot;\; s_1)
\end{equation}
between any pair of scale levels $s_2 > s_1$,
and which then also implies that the scale-space properties of the
zero-order discrete analogue of the Gaussian kernel, in terms of
{\em e.g.\/}, non-creation of new local extrema or new zero-crossings
from finer to coarser levels of scale, are also carried over to the
resulting discrete approximations of Gaussian derivative responses.
In these respects, this discretisation approach, of combining the
discrete analogue of the Gaussian kernel with central differences,
has specific both theoretical and computational advantages at very
fine scale levels.

\subsection{Hybrid discretisation based on the normalised sampled
  Gaussian kernel in combination with central differences}
\label{app-hybr-disc-normsamplgaussdiff}

Inspired by the computational efficiency of the discrete derivative
approximations described in the previous section, and obtained
by convolution with the discrete analogue of the Gaussian kernel
according to (\ref{eq-disc-gauss-app}) followed by central
differences $\delta_{x^{\alpha}}$, we will in this paper analyse the
effects of computing corresponding discrete derivative approximations,
by applying central difference operators $\delta_{x^{\alpha}}$ to the
result of smoothing a discrete signal with the normalised sampled
Gaussian kernel according to (\ref{eq-def-norm-sampl-gauss-app}).

The motivation for this
operation is again to obtain substantially higher computational
efficiency compared to explicit convolutions with multiple sampled
Gaussian derivative kernels according to
(\ref{eq-sampl-gauss-der-app}), when to compute
Gaussian derivative responses for multiple orders of differentiation
$\alpha$.

The equivalent convolution kernels for this discretisation method will
then be of the form
\begin{equation}
  \label{eq-hybr-normsampl-disc-der-app}
  T_{\hybrnormsampl,x^{\alpha}}(n;\; s)
  = (\delta_{x^{\alpha}} T_{\normsampl})(n;\; s).
\end{equation}
As for the above discretisation method with convolution with the
discrete analogue of the Gaussian kernel followed by central
differences, the intention with this definition of a corresponding
discrete equivalent convolution kernel is, however, again only to be
used a theoretical modelling tool, and not for any actual
implementations.
The reason for this is that, in situations when multiple Gaussian
derivative responses are to be computed for the same input signal $f$,
it is computationally much more efficient to first perform the smoothing with
the normalised sampled Gaussian kernel once and for all, and then apply each central
difference operator $\delta_{x^{\alpha}}$.

This hybrid discretisation method was proposed among suggestions to
future work in Section~6 in (Lindeberg \citeyear{Lin24-JMIV}) and with
theoretical properties of that discretisation method analysed in
Footnote~13 in that paper. The quantitive performance of this
hybrid discretisation method was, however, not evaluated experimentally,
which is a main subject of this paper.

\subsection{Hybrid discretisation based on the integrated Gaussian kernel
  in combination with central differences}
\label{app-hybr-disc-intgaussdiff}

In analogy with the treatment above concerning the discretisation
method, based on discrete convolution with the normalised
sampled Gaussian kernel (\ref{eq-def-norm-sampl-gauss-app})
followed by central differences $\delta_{x^{\alpha}}$,
we can also in a corresponding manner
instead use the integrated Gaussian kernel
(\ref{eq-def-int-gauss-kern-app})
as an initial smoothing
step and then apply central differences $\delta_{x^{\alpha}}$ to the
smoothed data, corresponding to the following equivalent discrete
derivative approximation kernel
\begin{equation}
  \label{eq-hybr-int-disc-der-app}
  T_{\hybrint,x^{\alpha}}(n;\; s)
  = (\delta_{x^{\alpha}} T_{\intdisc})(n;\; s).
\end{equation}
Again, this equivalent convolution kernel is, however,
not intended to be used for practical
implementations, since in situations when Gaussian derivative
responses are to be computed for different orders of differentiation
$\alpha$ for the same input signal, it is computationally much more
efficient to just perform the spatial smoothing operation once and for
all and then apply the different central difference operators
$\delta_{x^{\alpha}}$ to the smoothed signal. In this way, when
Gaussian derivative responses are to be computed from the same signal
$f$ for different orders of differentiation $\alpha$, it will be
computationally much more efficient to perform the computations in
this way, compared to {\em e.g.\/}\ performing explicit convolutions
with a corresponding set of integrated Gaussian derivative kernels
according to (\ref{eq-def-int-gauss-der}).

This hybrid discretisation method was also proposed among suggestions to
future work in Section~6 in (Lindeberg \citeyear{Lin24-JMIV}) and with
theoretical properties of that discretisation method analysed in
Footnote~13 in that paper. The quantitive performance of this
hybrid discretisation method was, however, not evaluated experimentally,
which is a main subject of this paper.

\section{Definition of discrete model signals for quantitatively
  evaluating the performance of the different discrete
  approximations of Gaussian derivative operator}
\label{sec-def-model-signals}

This appendix section describes how the discrete input model signals
(Gaussian blobs for interest point detection, diffuse edges for edge detection,
or Gaussian ridges for ridge detection) are generated for 
evaluating the scale selection properties for the different types
of feature detectors (Laplacian interest point detection, determinant of the
Hessian interest point detection, gradient magnitude edge detection or principal
curvature ridge detection) experimentally evaluated
in Section~\ref{sec-results-scsel}. Specifically,
this appendix section describes how the methods for discrete
approximations of the continuous input signals are determined for the
different methods for discrete approximations of the Gaussian
derivative operators.

\subsection{Methodology for defining the discrete input data for the scale
selection experiments}

When generating input data for different values of the reference scale $\sigma_0$,
we 
\begin{itemize}
\item
  for the purpose of Laplacian or determinant of the Hessian interest point detection,
  for generating a discrete blob, to be used as input for computing the
  scale-normalised differential entities $(\nabla_{\norm}^2 L)(x, y;\; s)$
  and $(\det {\mathcal H}_{\norm} L)(x, y;\; s)$
  for interest point detection at multiple scales, 
  we use a discrete approximation of Gaussian blob
  (Equation~(\ref{eq-def-2D-gauss-blob-s0})),
  %(\ref{eq-def-2D-gauss-blob-s0})
  obtained by convolving a 2-D discrete delta function
  \begin{equation}
    \delta(x, y)
    =
    \left\{
      \begin{array}{ll}
        1 & \mbox{if $(x, y) = (0, 0)$}, \\
        0 & \mbox{otherwise},
      \end{array}
    \right.
  \end{equation}
  with
  a discrete approximation of a Gaussian kernel along each dimension,
\item
  for the purpose of edge detection, for generating a discrete diffuse
  edge, to the used as input for computing the scale-normalised
  gradient magnitude $L_{v,\norm}(0, 0;\; s)$ at multiple scales,
  we use a discrete approximation of a diffuse edge
  % (\ref{eq-def-ideal-blur-edge-s0}),
  (Equation~(\ref{eq-def-ideal-blur-edge-s0})),
  obtained by convolving an ideal step edge function
    \begin{equation}
    H(x)
    =
    \left\{
      \begin{array}{ll}
        + 1/2 & \mbox{if $x > 0$}, \\
        0 & \mbox{if $x = 0$}, \\
        - 1/2 & \mbox{if $x < 0$}, \\
      \end{array}
    \right.
  \end{equation}
  with a discrete approximation of a Gaussian kernel along the $x$-direction,
\item
  for the purpose of ridge detection, for generating a discrete
  Gaussian ridge, to be used as input for computing the
  scale-normalised principal curvature measure
  $L_{pp,\norm}(0, 0;\; s)$ at multiple scales,
  we use a discrete approximation of a Gaussian ridge
  (\ref{eq-def-ideal-ridge-s0}),
  obtained by convolving the 2-D extension of
   a 1-D discrete delta function
  \begin{equation}
    \delta(x)
    =
    \left\{
      \begin{array}{ll}
        1 & \mbox{if $x = 0$}, \\
        0 & \mbox{otherwise},
      \end{array}
    \right.
  \end{equation}
  along the $y$-direction with
  a discrete approximation of a Gaussian kernel along the $x$-direction.
\end{itemize}
Since a main purpose of the experiments defined in Section~\ref{sec-meth-consist-scales} and reported in Section~\ref{sec-results-scsel}  is to measure
the {\em consistency\/}
between the characteristic scales in the input with the characteristic
scale values obtained by multi-scale processing of the input data,
we choose the discretisation method for generating the input data for
the scale estimation algorithm in
relation to the discretisation method used for defining the discrete
derivative approximations in the following way:
\begin{itemize}
\item
  when evaluating the discretisation based on the discrete analogue of
  Gaussian derivatives, we
  use the discrete analogue of the Gaussian kernel as the discrete
  convolution kernel when generating the input data,
\item
  when evaluating the discretisation method based on sampled Gaussian derivatives,
  we use the sampled Gaussian kernel as the discrete
  convolution kernel when generating the input data,
\item
  when evaluating the discretisation method based on integrated Gaussian derivatives,
  we use the integrated Gaussian kernel as the discrete
  convolution kernel when generating the input data,
\item
  when evaluating the hybrid method, based on convolution with the normalised
  sampled Gaussian kernel followed by central differences, we use the normalised
  sampled Gaussian kernel as the discrete convolution kernel when generating
  the input data, and
\item
  when evaluating the hybrid method, based on convolution with the 
  integrated Gaussian kernel followed by central differences, we use the
  integrated Gaussian kernel as the discrete convolution kernel when generating
  the input data.
\end{itemize}
Thus, since the purpose of the experiments in Section~3.2 in the main
text is to measure the mutual
consistency between the discrete approximations of Gaussian derivative
responses between different scales, the intention is to use an as conceptually
similar discretisation method for generating
the input data, as will be used when analysing the same data by the
feature detection algorithms.
Since the generation of the input data will, however, hence differ between the
different discretisation methods for Gaussian derivatives, some care therefore
needs to be taken when interpreting the experimental results.

\subsection{Motivation for the experimental paradigm in terms of
  internal consistency of scale estimates}

The motivation for evaluating the performance of the different use
cases for feature detection with automatic scale
selection with regard to the internal consistency between the discrete
approximations of the input and the discrete approximations of the
Gaussian derivative kernels, as opposed to comparing all the different
discretisation methods relative to the same ground truth, is that the
definition of a common ground truth itself becomes very problematic at very fine
scale levels, in turn, because of the quantitative measurement of the
amount of effective discrete blur for a given value of the scale parameter $s$ may differ
significantly between the different discrete approximation methods at very fine
scales, although the value of the underlying scale parameter is the
same (see Sections~2.8.2--2.8.4 and~3.8.2 in
Lindeberg (\citeyear{Lin24-JMIV}) as well as
Figures~\ref{fig-spat-spread-meas-with-hybr-discr}
and~\ref{fig-spat-spread-offset-meas-with-hybr-discr}
in this paper
for numerical quantifications of this property).

\begin{figure*}[hbtp]
  \begin{center}
    \begin{tabular}{ccc}
      {\footnotesize\em Discrete Gaussian blob for $\sigma = 1/2$}
      & {\footnotesize\em Discrete Gaussian blob for $\sigma = 1$}
      & {\footnotesize\em Discrete Gaussian blob for $\sigma = 2$} \\
      \includegraphics[width=0.30\textwidth]{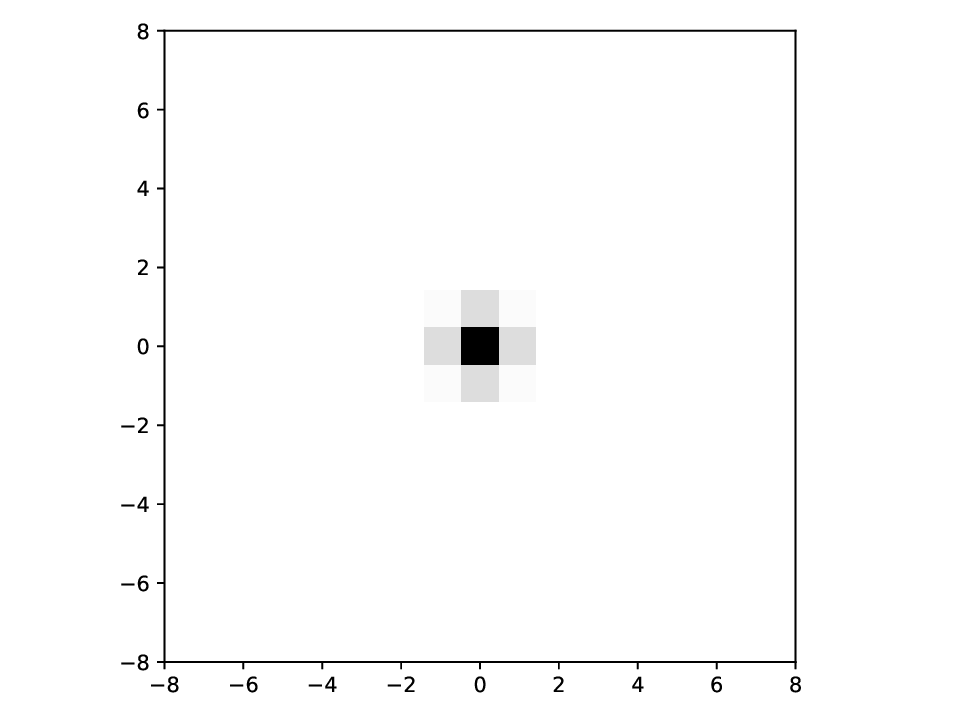}
      & \includegraphics[width=0.30\textwidth]{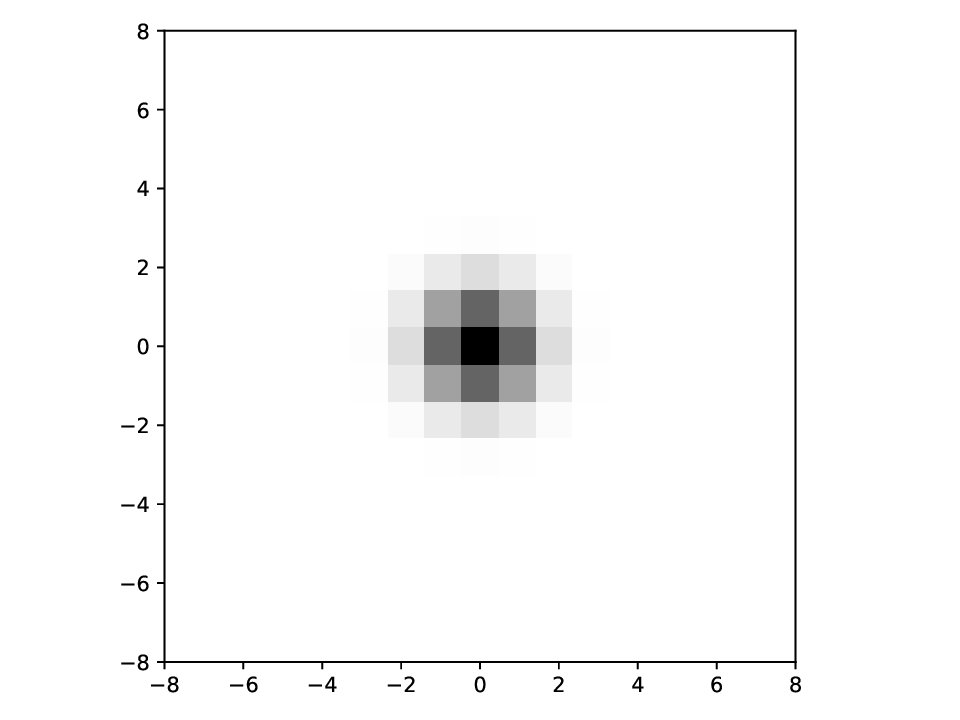}
      & \includegraphics[width=0.30\textwidth]{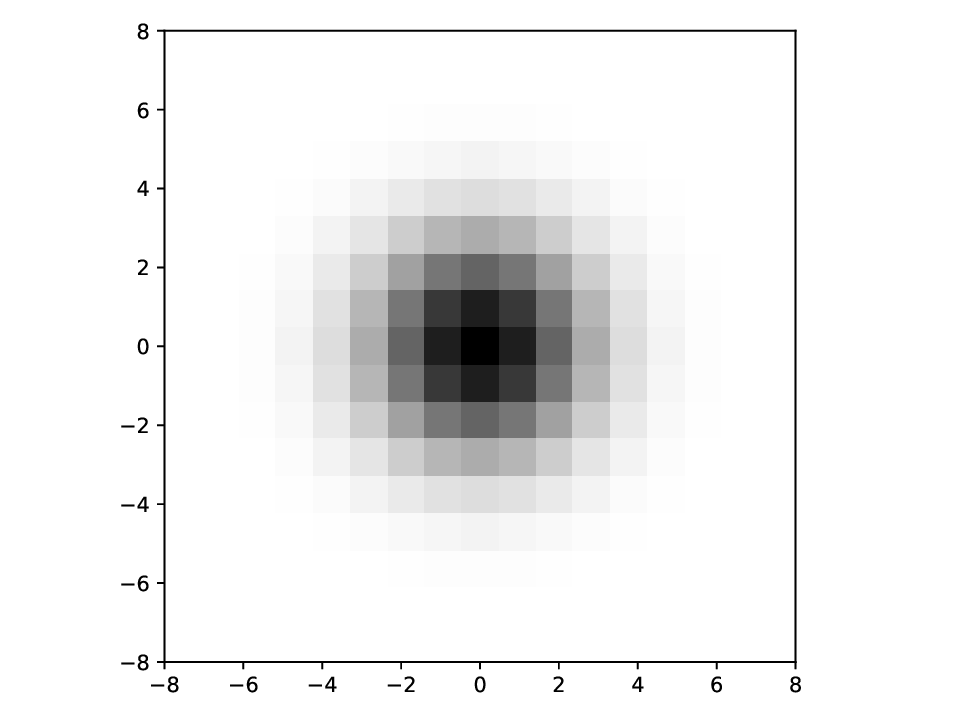}\\
      $\,$ \\
      {\footnotesize\em Laplacian response $\nabla_{\norm}^2 L$ at $\sigma = 1/2$}
      & {\footnotesize\em Laplacian response $\nabla_{\norm}^2 L$ at $\sigma = 1$}
      & {\footnotesize\em Laplacian response $\nabla_{\norm}^2 L$ at $\sigma = 2$} \\
      \includegraphics[width=0.30\textwidth]{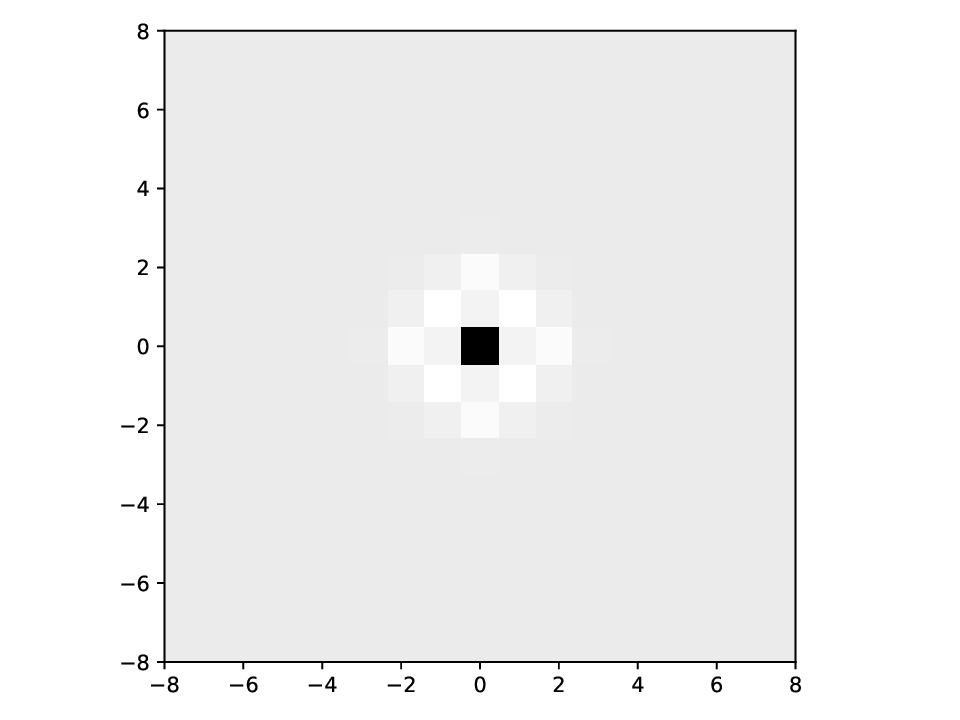}
      & \includegraphics[width=0.30\textwidth]{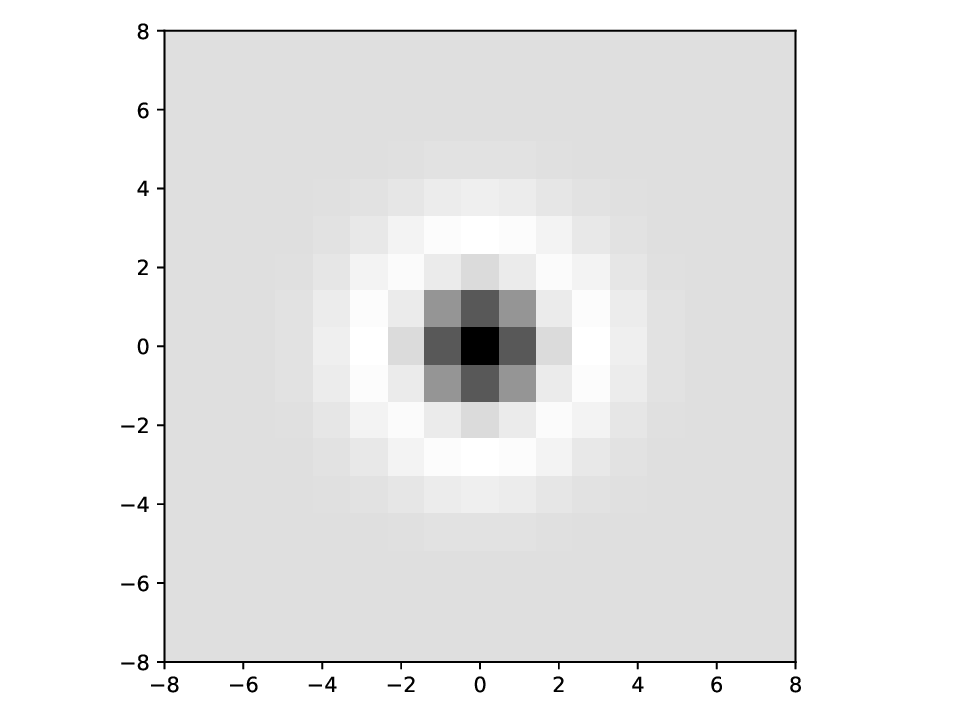}
      & \includegraphics[width=0.30\textwidth]{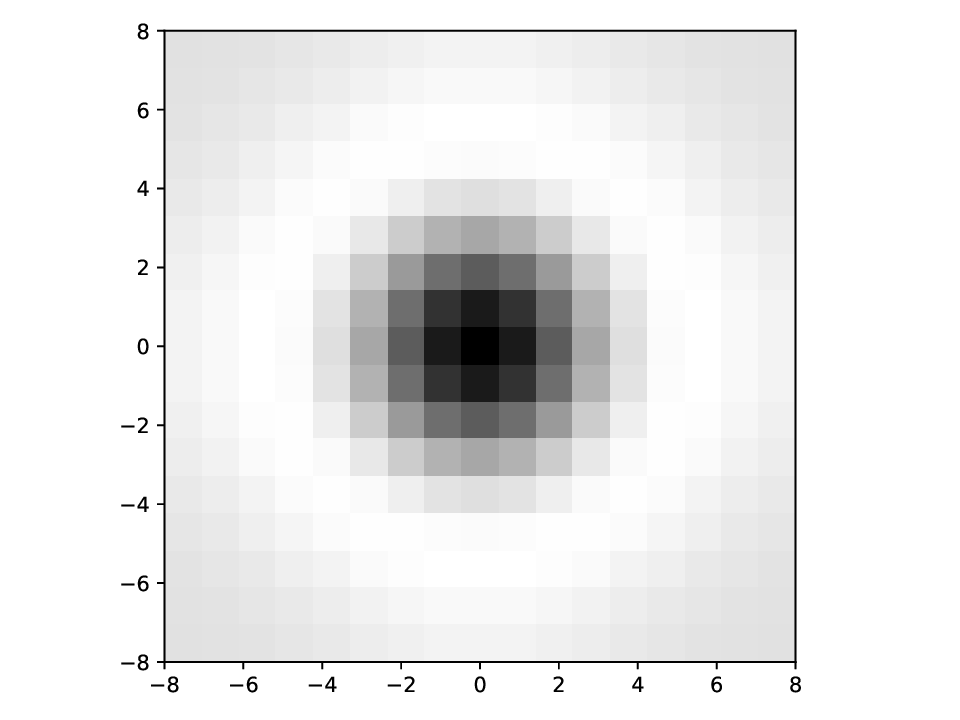}\\
      $\,$ \\
      {\footnotesize\em Variability over scale for $\nabla_{\norm}^2 L$}
      & {\footnotesize\em Variability over scale for $\nabla_{\norm}^2 L$}
      & {\footnotesize\em Variability over scale for $\nabla_{\norm}^2 L$} \\
      \includegraphics[width=0.30\textwidth]{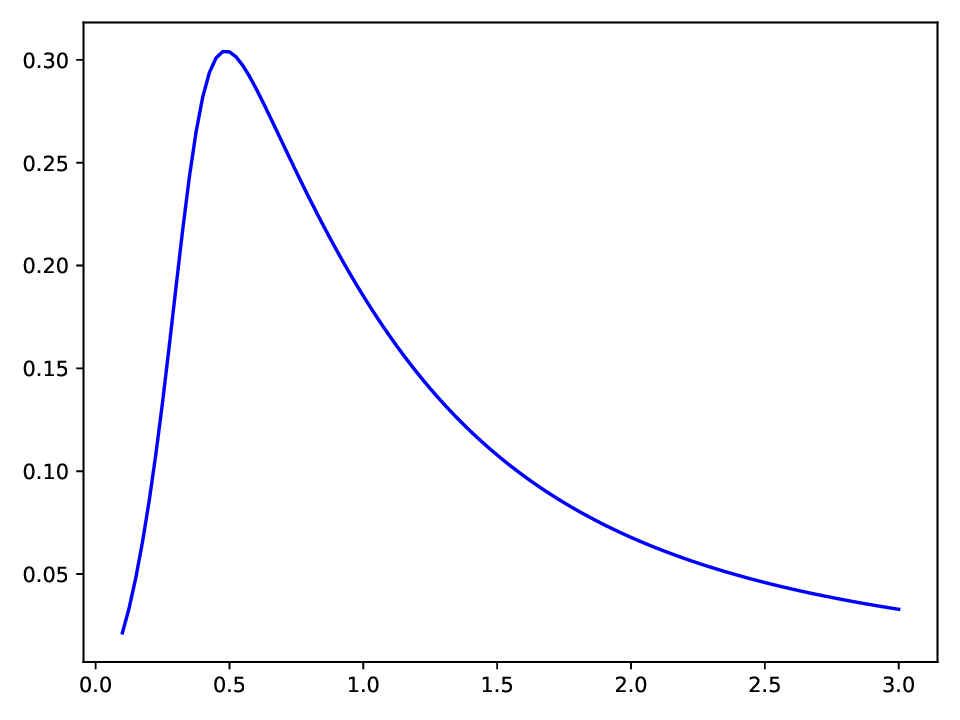}
      & \includegraphics[width=0.30\textwidth]{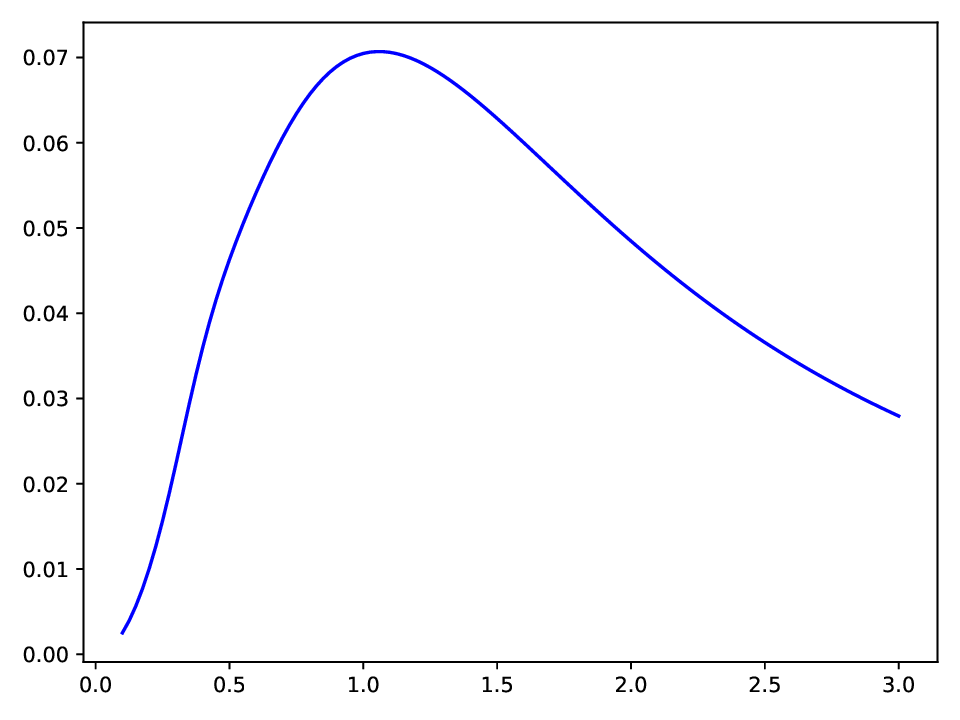}
      & \includegraphics[width=0.30\textwidth]{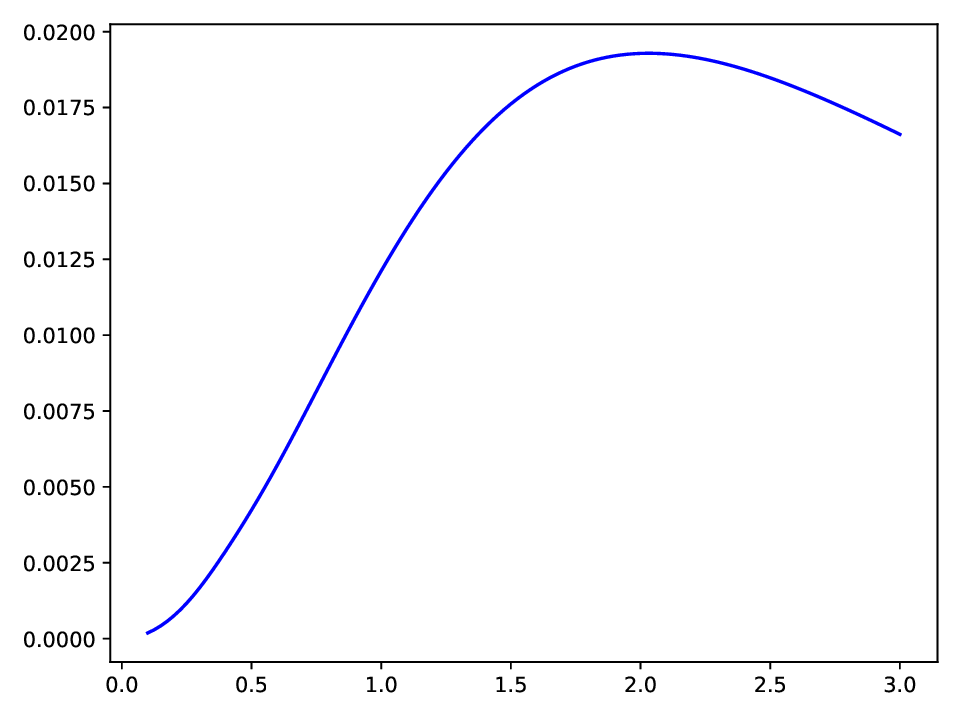}
    \end{tabular}
  \end{center}
  \caption{Visualisation of the conceptual steps involved when
    computing scale estimates using {\em Laplacian interest point detection
    with automatic scale selection\/}, here using the hybrid
    discretisation method based on convolution with normalised sampled
    Gaussian kernels followed by central differences.
    {\bf (top row)} Input images for three different sizes $\sigma_0 = 1/2$,
    $\sigma_0 = 1$ and $\sigma_0 = 2$ of the discrete Gaussian blobs
    (with the contrast reversed in the visualisation).
    {\bf (middle row)} The result of computing a discrete
    approximation of the Laplacian response
    $\nabla L = L_{xx} + L_{yy}$ for the different 
    input images in the top row, at the scale
    levels $\sigma = 1/2$, $\sigma = 1$ and $\sigma = 2$,
    respectively.
    {\bf (bottom row)} Graphs of the variability over scale $\sigma$ for (the
    negative value of) the scale-normalised Laplacian response
    $\nabla_{\norm}^2 L = \sigma^2 \, (L_{xx} + L_{yy})$ for each one of
    the images in the top row.
    Note that the selected scales $\hat{\sigma}$, where the maxima in these
    graphs are assumed, do rather well correspond to the scales
    $\sigma_0$ in the input data.}
 %    (The reason, why the peak values over scale are different for the
%     differently sized input, is because the input images have been
%     obtained by convolution with discrete normalised sampled
%     Gaussian kernels, with their normalisation to unit norm. Thereby, the
%   peak values of the discrete Gaussian blobs will be different for the
%   different values of $\sigma_0$.)}
  \label{fig-lapl-intpt-det}
\end{figure*}

\begin{figure*}[hbtp]
  \begin{center}
    \begin{tabular}{ccc}
      {\footnotesize\em Discrete Gaussian blob for $\sigma = 1/2$}
      & {\footnotesize\em Discrete Gaussian blob for $\sigma = 1$}
      & {\footnotesize\em Discrete Gaussian blob for $\sigma = 2$} \\
      \includegraphics[width=0.30\textwidth]{blob0p5.eps}
      & \includegraphics[width=0.30\textwidth]{blob1.eps}
      & \includegraphics[width=0.30\textwidth]{blob2.eps}\\
      $\,$ \\
      {\footnotesize\em detHessian $\det {\mathcal H}_{\norm} L$ at $\sigma = 1/2$}
      & {\footnotesize\em detHessian $\det {\mathcal H}_{\norm} L$ at $\sigma = 1$}
      & {\footnotesize\em detHessian $\det {\mathcal H}_{\norm}  L$ at $\sigma = 2$} \\
      \includegraphics[width=0.30\textwidth]{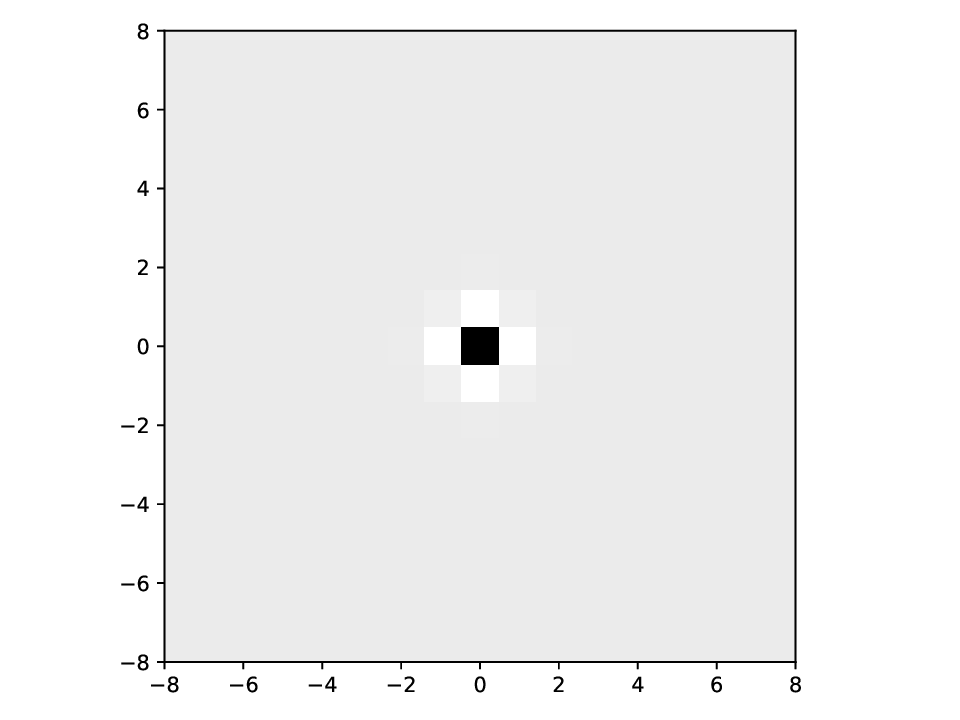}
      & \includegraphics[width=0.30\textwidth]{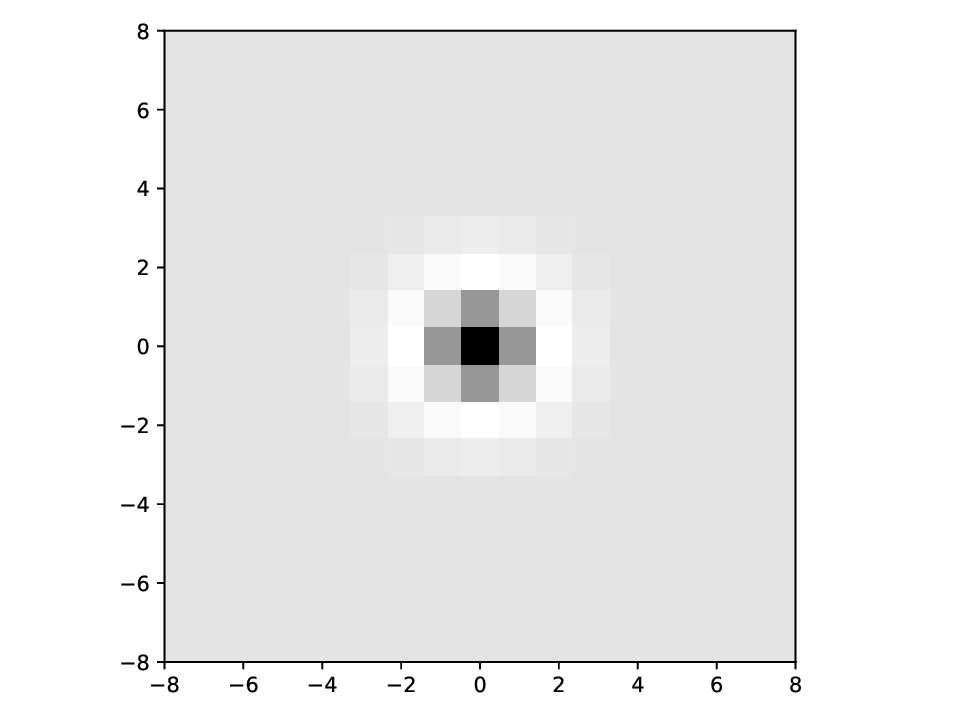}
      & \includegraphics[width=0.30\textwidth]{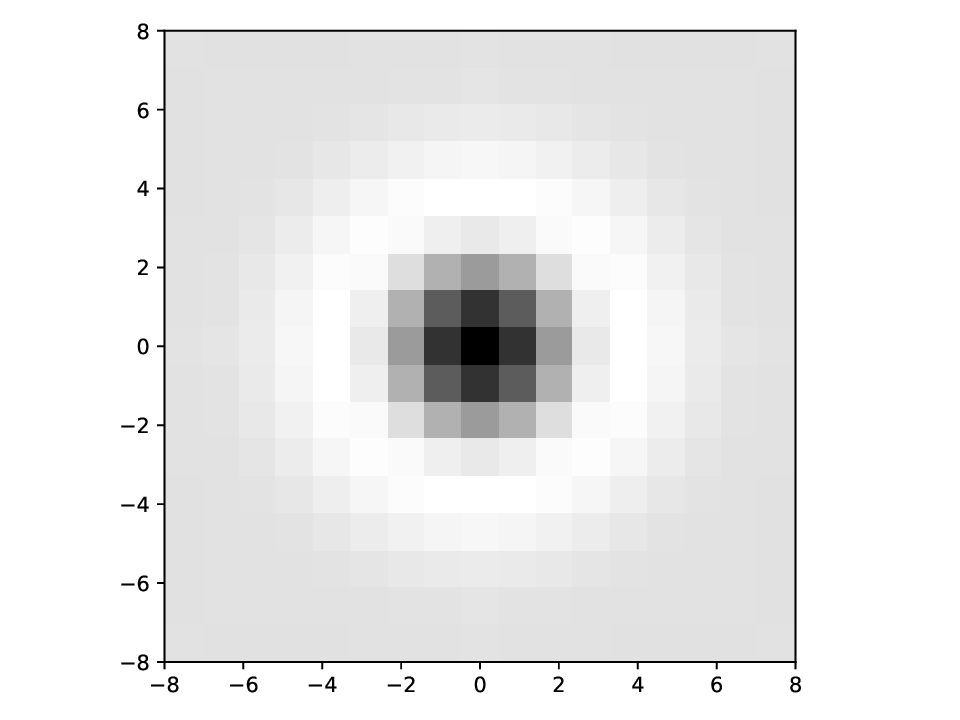}\\
      $\,$ \\
      {\footnotesize\em Variability over scale $\det {\mathcal H}_{\norm} L$}
      & {\footnotesize\em Variability over scale $\det {\mathcal H}_{\norm} L$}
      & {\footnotesize\em Variability over scale $\det {\mathcal H}_{\norm} L$} \\
      \includegraphics[width=0.30\textwidth]{blob0p5-lapl-scspsign.eps}
      & \includegraphics[width=0.30\textwidth]{blob1-lapl-scspsign.eps}
      & \includegraphics[width=0.30\textwidth]{blob2-lapl-scspsign.eps}
    \end{tabular}
  \end{center}
  \caption{Visualisation of the conceptual steps involved when
    computing scale estimates using {\em determinant of the Hessian interest point detection
    with automatic scale selection\/}, here using the hybrid
    discretisation method based on convolution with normalised sampled
    Gaussian kernels followed by central differences.
    {\bf (top row)} Input images for three different sizes $\sigma_0 = 1/2$,
    $\sigma_0 = 1$ and $\sigma_0 = 2$ of the discrete Gaussian blobs
    (with the contrast reversed in the visualisation).
    {\bf (middle row)} The result of computing a discrete
    approximation of the determinant of the Hessian response
    $\det {\mathcal H} L = L_{xx} \, L_{yy} - L_{xy}^2$ for the different 
    input images in the top row, at the scale
    levels $\sigma = 1/2$, $\sigma = 1$ and $\sigma = 2$,
    respectively (again with the contrast reversed in the visualisation).
    {\bf (bottom row)} Graphs of the variability over scale $\sigma$ for (the
    negative value of) the scale-normalised determinant of the Hessian response
    $\det {\mathcal H}_{\norm} L = \sigma^4 \, (L_{xx} \, L_{yy} - L_{xy}^2)$
    for each one of the images in the top row.
    Note that the selected scales $\hat{\sigma}$, where the maxima in these
    graphs are assumed, do rather well correspond to the scales
    $\sigma_0$ in the input data.}
 %    (The reason, why the peak values over scale are different for the
%     differently sized input, is because the input images have been
%     obtained by convolution with discrete normalised sampled
%     Gaussian kernels, with their normalisation to unit norm. Thereby, the
%   peak values of the discrete Gaussian blobs will be different for the
%   different values of $\sigma_0$.)}
  \label{fig-dethess-intpt-det}
\end{figure*}

\begin{figure*}[hbtp]
  \begin{center}
    \begin{tabular}{ccc}
      {\footnotesize\em Discrete diffuse edge for $\sigma = 1/2$}
      & {\footnotesize\em Discrete diffuse edge for $\sigma = 1$}
      & {\footnotesize\em Discrete diffuse edge for $\sigma = 2$} \\
      \includegraphics[width=0.30\textwidth]{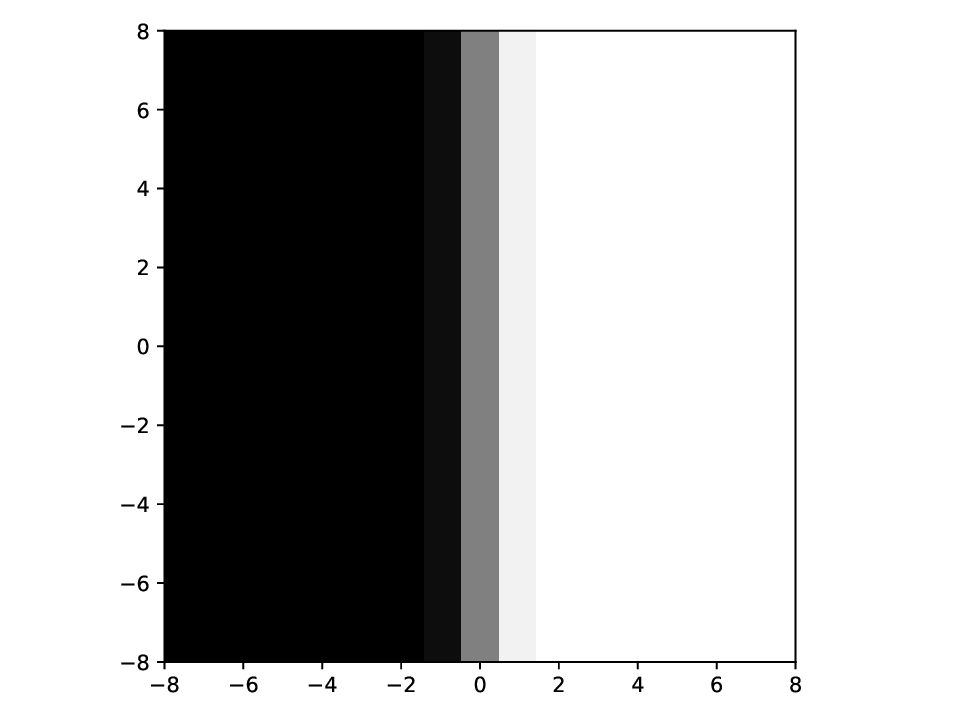}
      & \includegraphics[width=0.30\textwidth]{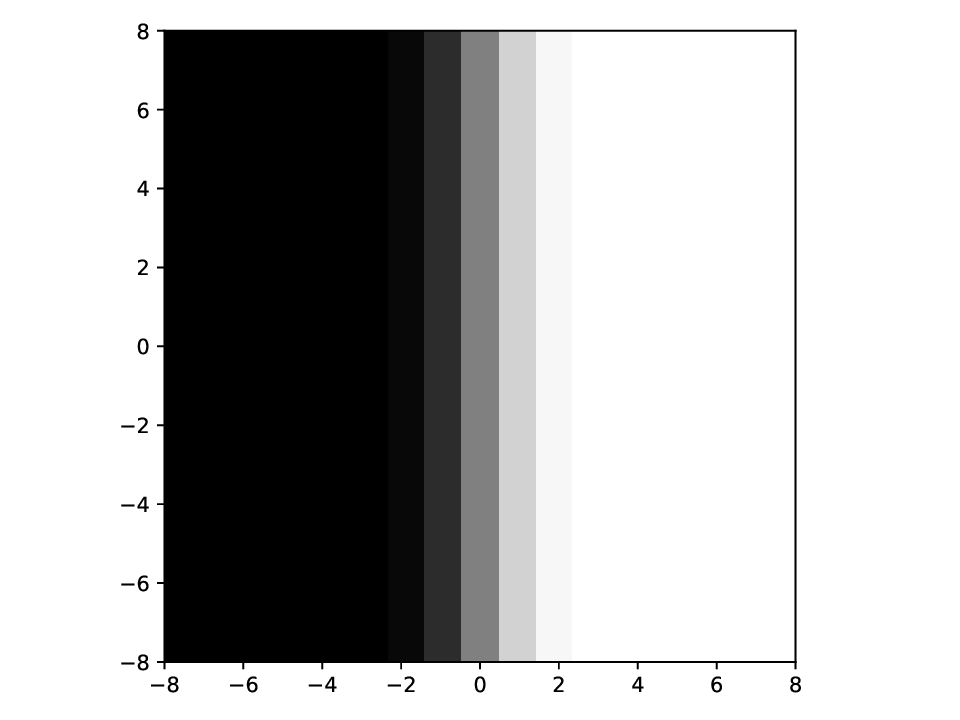}
      & \includegraphics[width=0.30\textwidth]{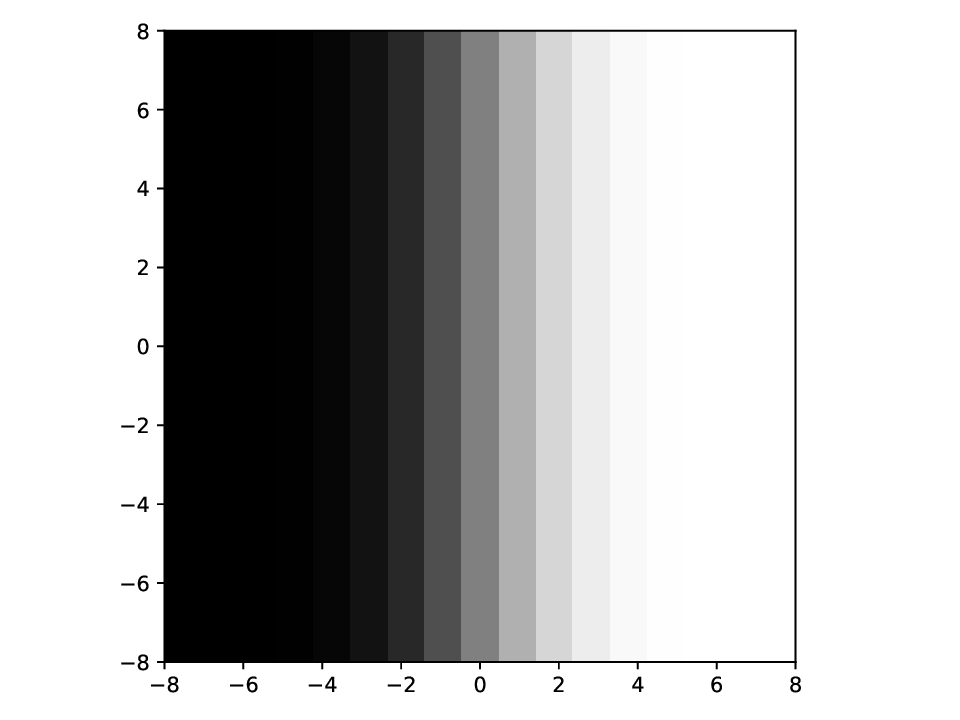}\\
      $\,$ \\
      {\footnotesize\em gradient magnitude $L_{v,\norm}$ at $\sigma = 1/2$}
      & {\footnotesize\em gradient magnitude $L_{v,\norm}$ at $\sigma = 1$}
      & {\footnotesize\em gradient magnitude $L_{v,\norm}$ at $\sigma = 2$} \\
      \includegraphics[width=0.30\textwidth]{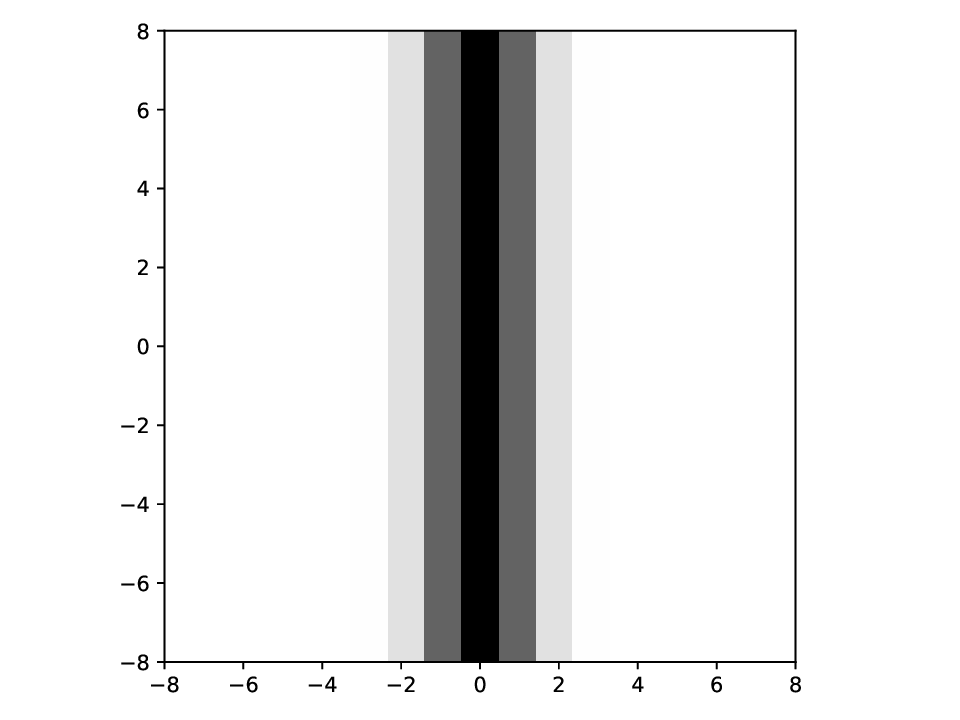}
      & \includegraphics[width=0.30\textwidth]{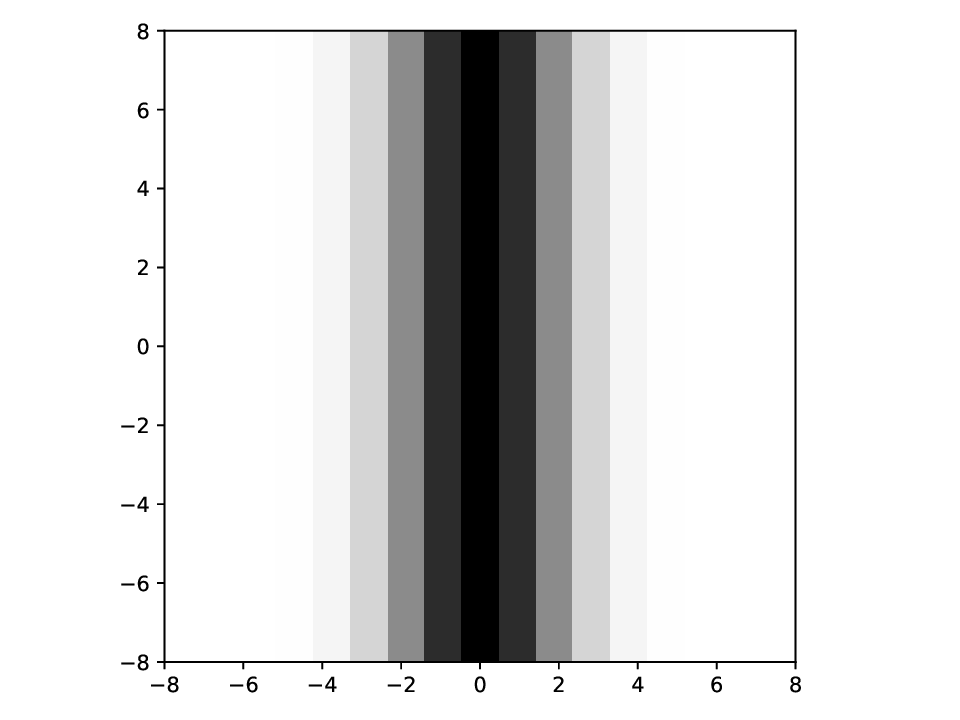}
      & \includegraphics[width=0.30\textwidth]{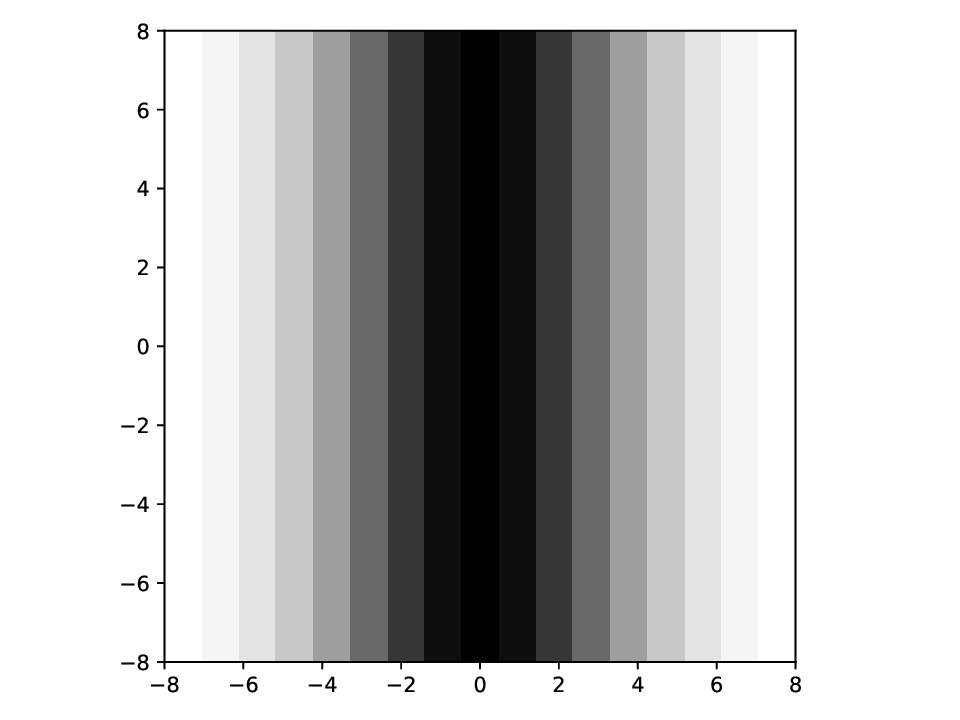}\\
      $\,$ \\
      {\footnotesize\em Variability over scale of $L_{v,\norm}$}
      & {\footnotesize\em Variability over scale of $L_{v,\norm}$}
      & {\footnotesize\em Variability over scale of $L_{v,\norm}$} \\
      \includegraphics[width=0.30\textwidth]{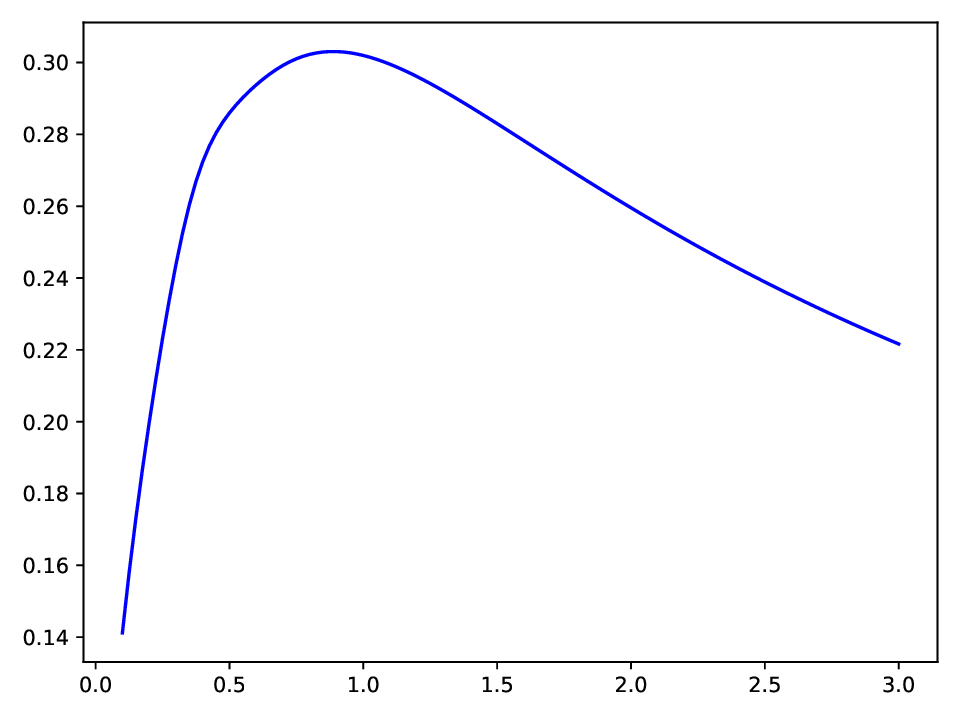}
      & \includegraphics[width=0.30\textwidth]{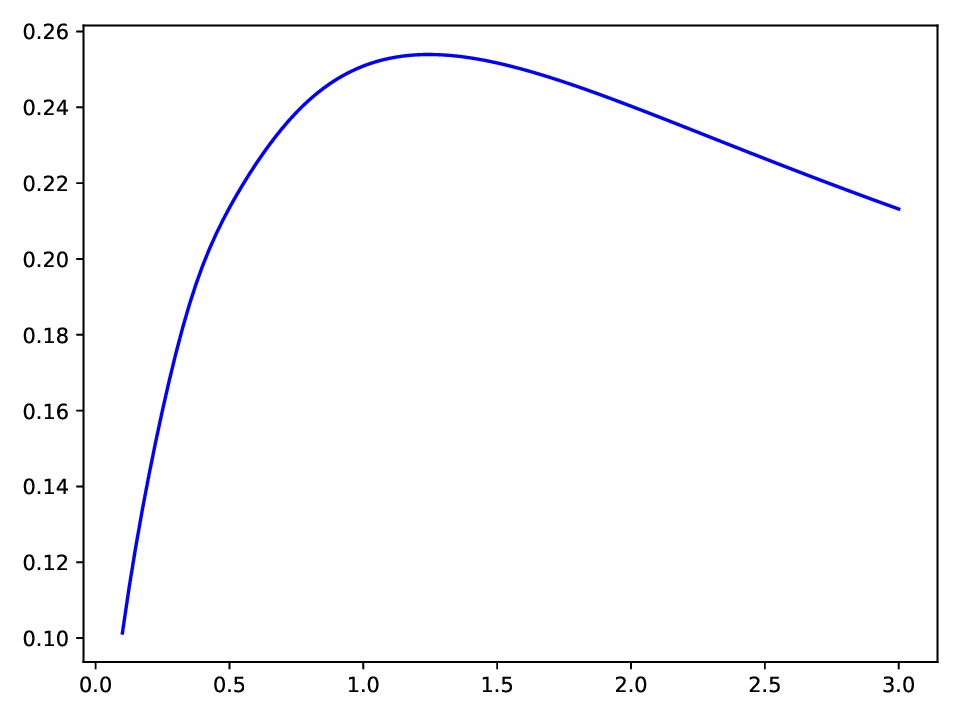}
      & \includegraphics[width=0.30\textwidth]{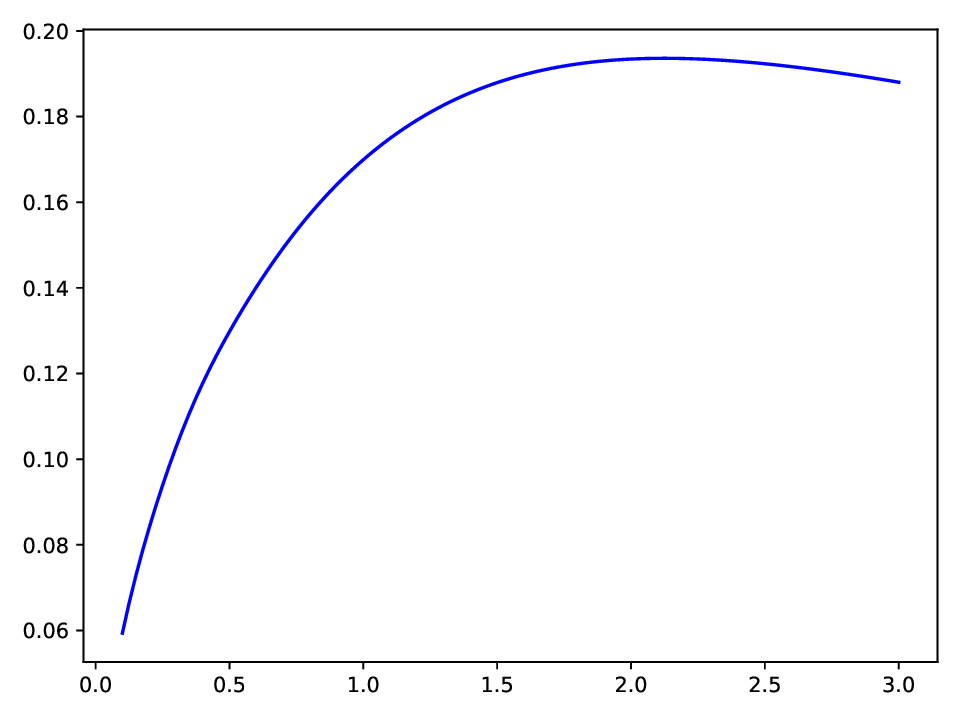}
    \end{tabular}
  \end{center}
  \caption{Visualisation of the conceptual steps involved when
    computing scale estimates using {\em gradient magnitude based edge detection
    with automatic scale selection\/}, here using the hybrid
    discretisation method based on convolution with normalised sampled
    Gaussian kernels followed by central differences.
    {\bf (top row)} Input images for three degrees of diffuseness $\sigma_0 = 1/2$,
    $\sigma_0 = 1$ and $\sigma_0 = 2$ for the discrete diffuse edges.
    {\bf (middle row)} The result of computing a discrete
    approximation of the gradient magnitude 
    $L_v = \sqrt{L_x^2 + L_y^2}$ for the different 
    input images in the top row, at the scale
    levels $\sigma = 1/2$, $\sigma = 1$ and $\sigma = 2$,
    respectively (with the contrast reversed in the visualisation).
    {\bf (bottom row)} Graphs of the variability over scale $\sigma$ for
    the scale-normalised gradient magnitude
    $L_{v,\norm} = \sigma^{\gamma} \sqrt{L_x^2 + L_y^2}$
    for each one of the images in the top row, using the scale
    normalisation power $\gamma = 1/2$.
    Note that the selected scales $\hat{\sigma}$, where the maxima in these
    graphs are assumed, do rather well correspond to the scales
    $\sigma_0$ in the input data.}
 \label{fig-gradmagn-edge-det}
\end{figure*}

\begin{figure*}[hbtp]
  \begin{center}
    \begin{tabular}{ccc}
      {\footnotesize\em Discrete Gaussian ridge for $\sigma = 1/2$}
      & {\footnotesize\em Discrete Gaussian ridge for $\sigma = 1$}
      & {\footnotesize\em Discrete Gaussian ridge for $\sigma = 2$} \\
      \includegraphics[width=0.30\textwidth]{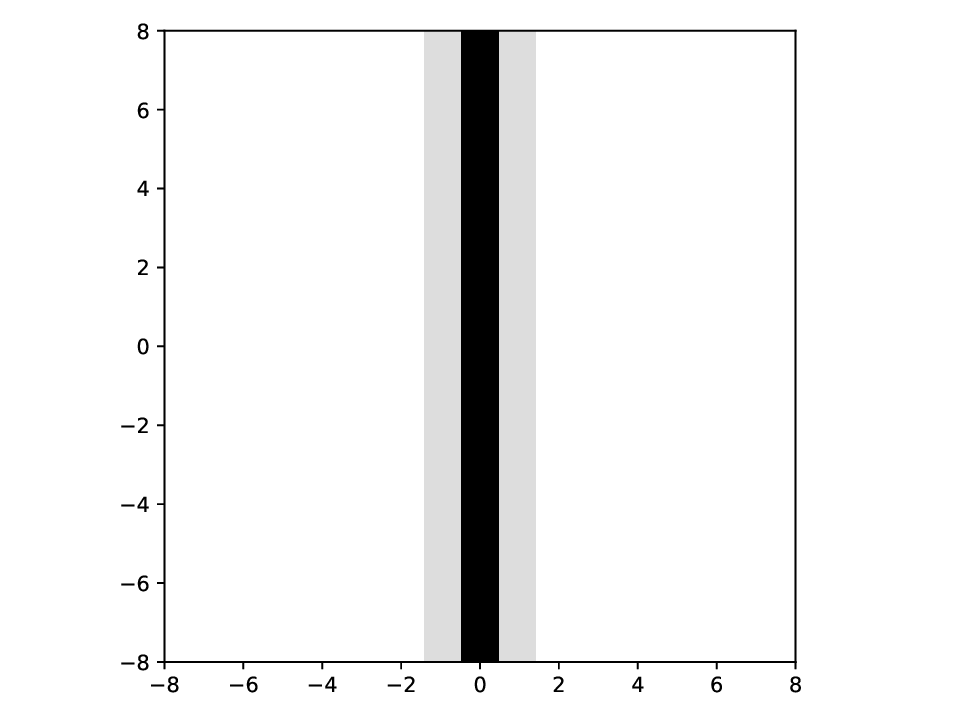}
      & \includegraphics[width=0.30\textwidth]{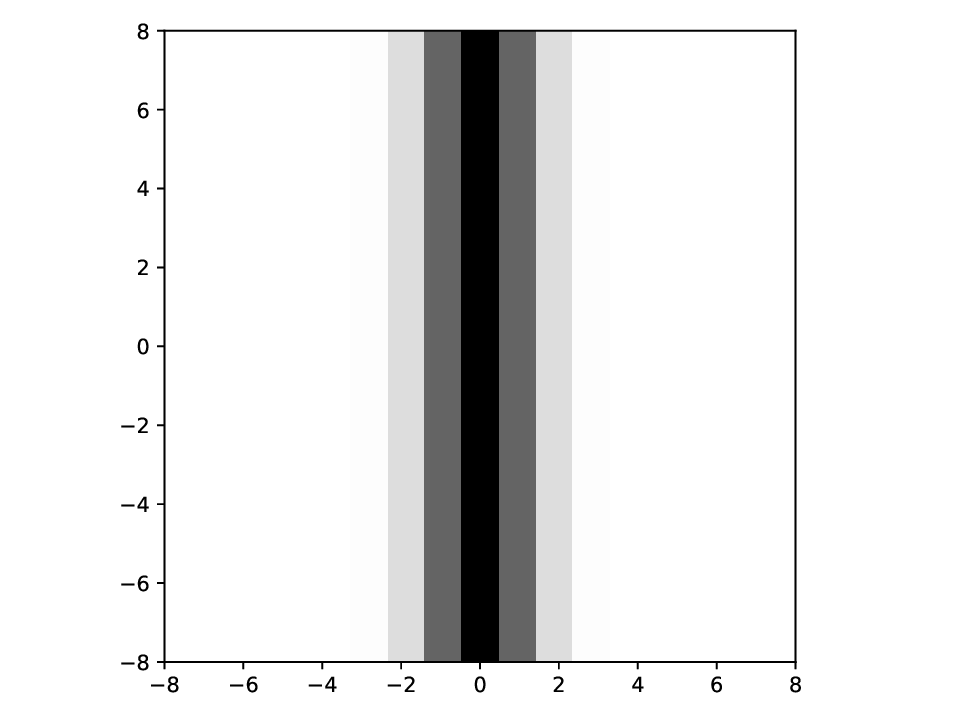}
      & \includegraphics[width=0.30\textwidth]{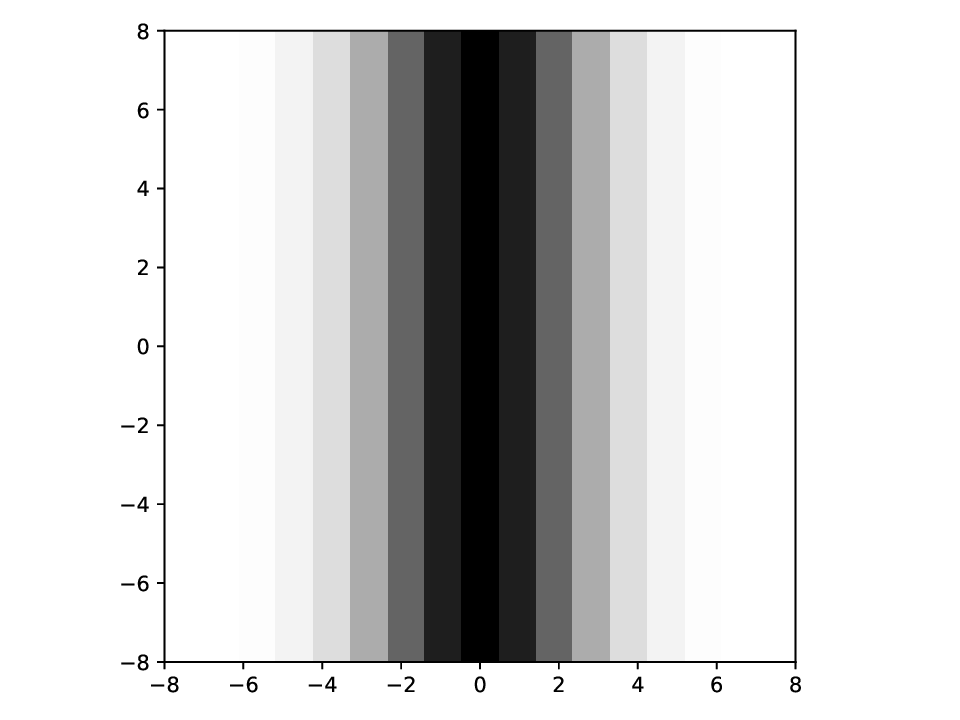}\\
      $\,$ \\
      {\footnotesize\em principal curvature $L_{pp,\norm}$ at $\sigma = 1/2$}
      & {\footnotesize\em principal curvature $L_{pp,\norm}$ at $\sigma = 1$}
      & {\footnotesize\em principal curvature $L_{pp,\norm}$ at $\sigma = 2$} \\
      \includegraphics[width=0.30\textwidth]{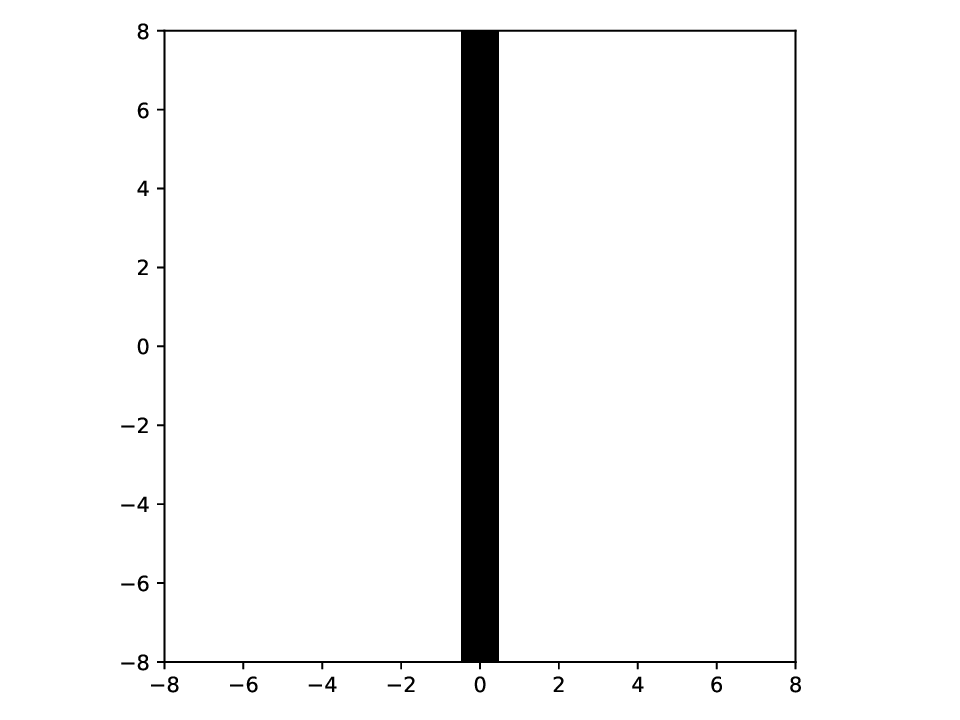}
      & \includegraphics[width=0.30\textwidth]{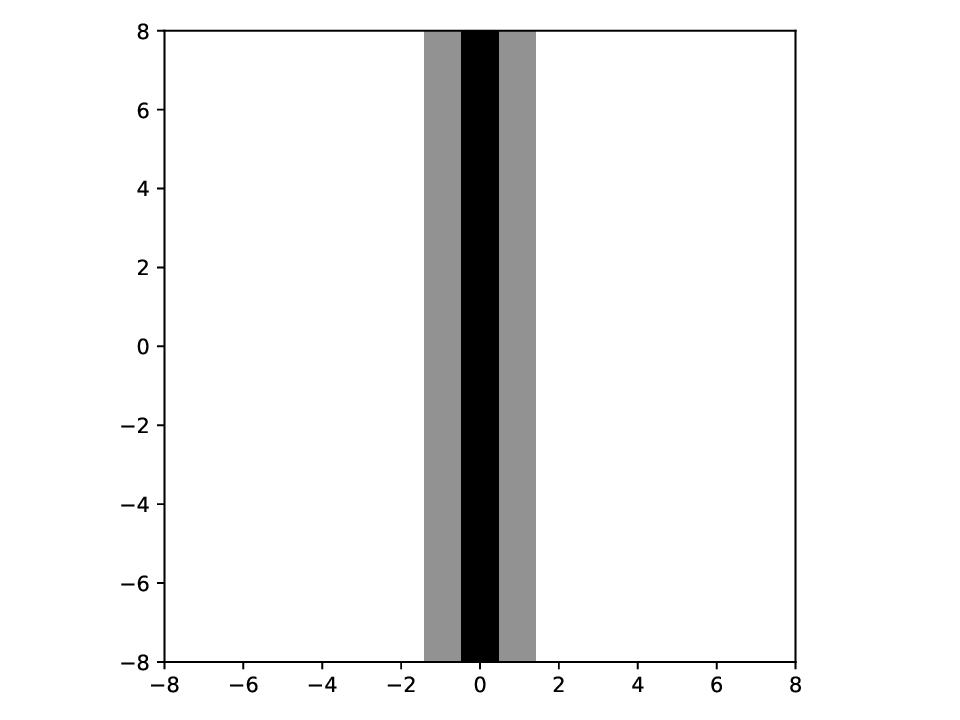}
      & \includegraphics[width=0.30\textwidth]{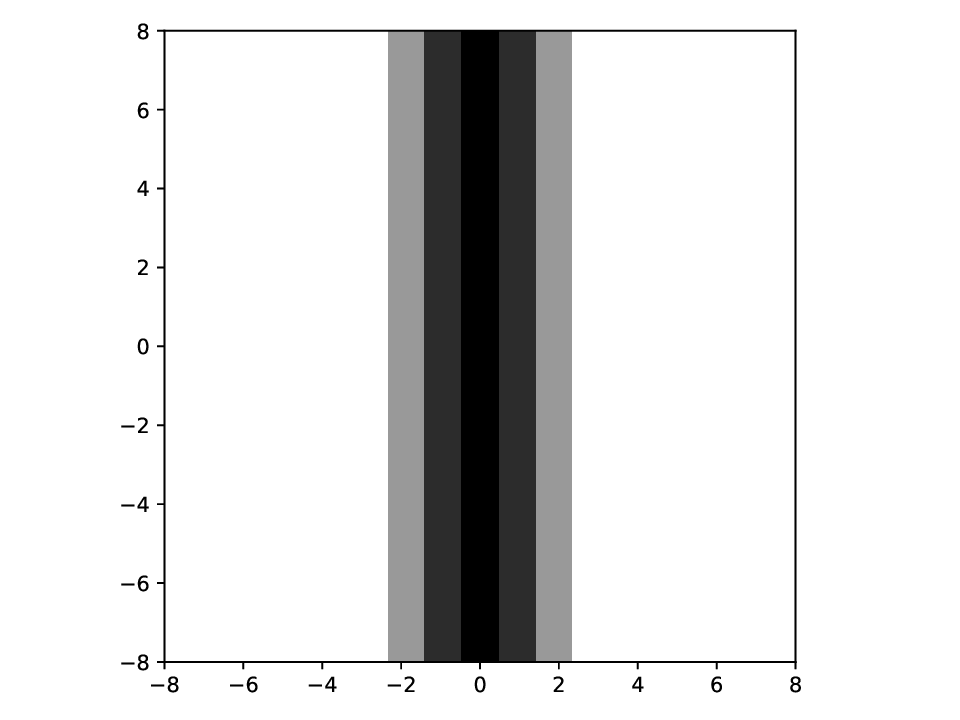}\\
      $\,$ \\
      {\footnotesize\em Variability over scale of $L_{pp,\norm}$}
      & {\footnotesize\em Variability over scale of $L_{pp,\norm}$}
      & {\footnotesize\em Variability over scale of $L_{pp,\norm}$} \\
      \includegraphics[width=0.30\textwidth]{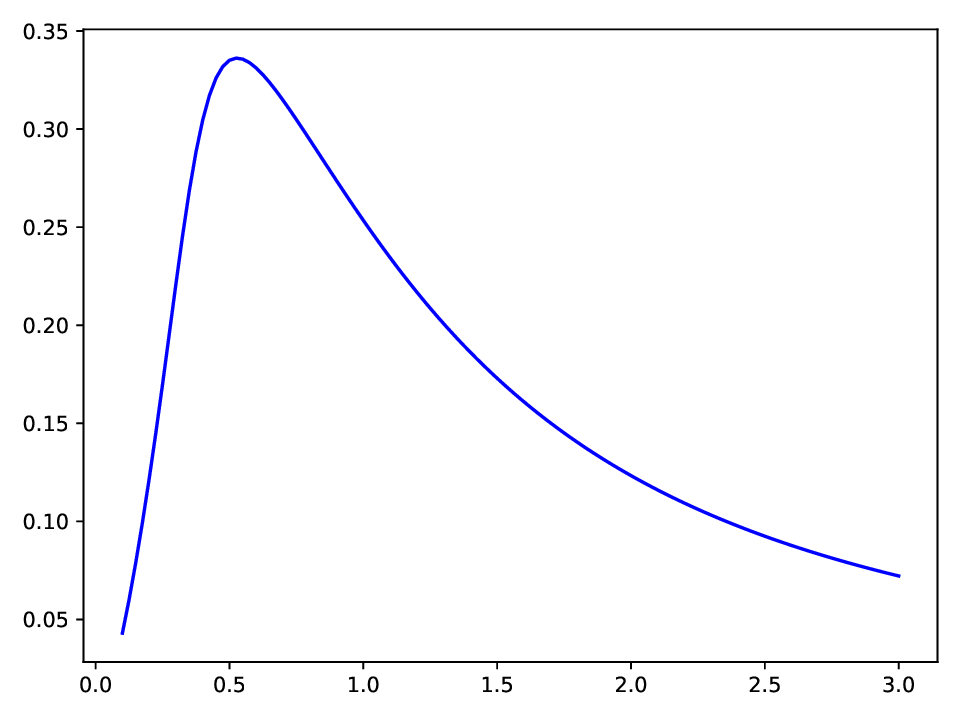}
      & \includegraphics[width=0.30\textwidth]{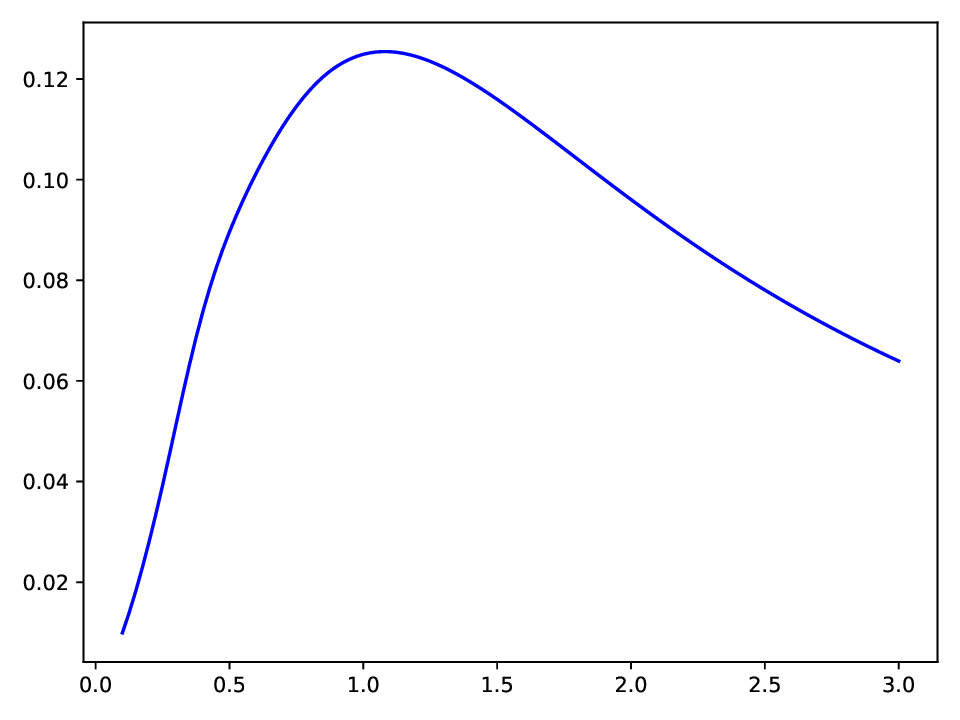}
      & \includegraphics[width=0.30\textwidth]{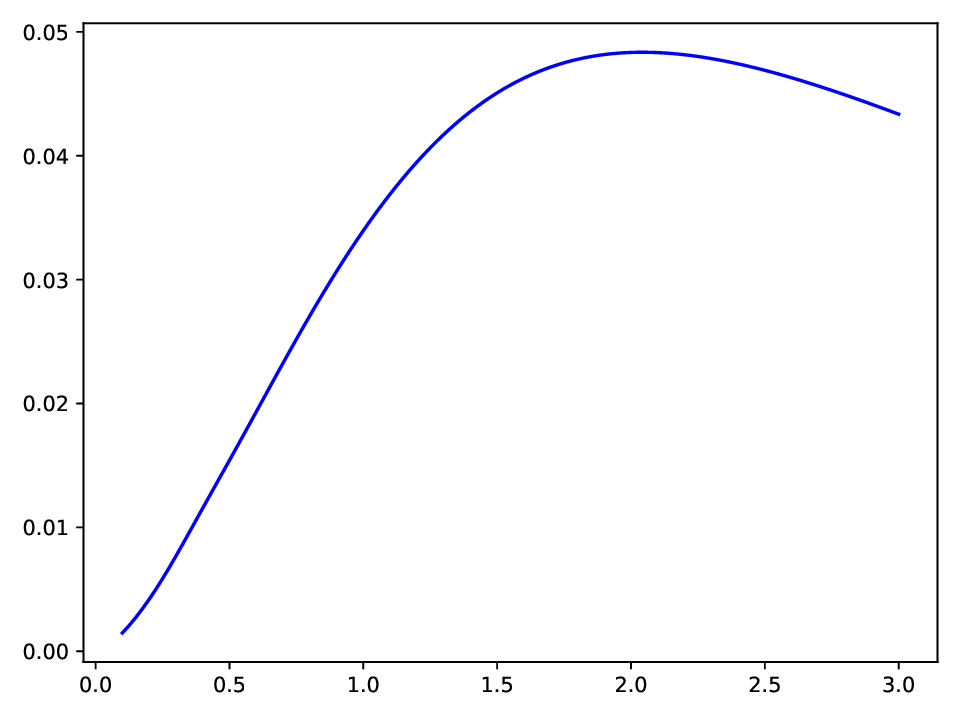}
    \end{tabular}
  \end{center}
  \caption{Visualisation of the conceptual steps involved when
    computing scale estimates using {\em principal curvature based ridge detection
    with automatic scale selection\/}, here using the hybrid
    discretisation method based on convolution with normalised sampled
    Gaussian kernels followed by central differences.
    {\bf (top row)} Input images for three different widths $\sigma_0 = 1/2$,
    $\sigma_0 = 1$ and $\sigma_0 = 2$ for the discrete Gaussian ridges.
    {\bf (middle row)} The result of computing a discrete
    approximation of the principal curvature based ridge strength measure
    $L_{pp} = L_{xx} + L_{yy} - \sqrt{(L_{xx} - L_{yy})^2 + 4 L_{xy}^2}$
    for the different  input images in the top row, at the scale
    levels $\sigma = 1/2$, $\sigma = 1$ and $\sigma = 2$,
    respectively (with the contrast reversed in the visualisation).
    {\bf (bottom row)} Graphs of the variability over scale $\sigma$ for
    the scale-normalised principal curvature based ridge strength measure
    $L_{pp,\norm} = \sigma^{2 \gamma} \, (L_{xx} + L_{yy} - \sqrt{(L_{xx} - L_{yy})^2 + 4 L_{xy}^2})$
    for each one of the images in the top row, using the scale
    normalisation power $\gamma = 3/4$.
    Note that the selected scales $\hat{\sigma}$, where the maxima in these
    graphs are assumed, do rather well correspond to the scales
    $\sigma_0$ in the input data.}
%     (One reason, why the peak values over scale are different for the
%     different degrees of diffuseness in the input, is because
%     the input images have been
%     obtained by convolution with discrete normalised sampled
%     Gaussian kernels, with their normalisation to unit norm. Thereby, the
%   peak values of the discrete Gaussian ridges will be different for the
%   different values of $\sigma_0$. Additionally, the magnitudes of the
%     scale-normalised differential derivative responses can neither be
%     expected to be equal, when using a value of the scale normalisation power
%     $\gamma \neq 1$.)}
  \label{fig-princcurv-ridge-det}
\end{figure*}

\section{Visualisations of image feature models, differential feature
  responses and variabilities in scale-normalised derivative responses
  over scale, which underlie the computation of scale estimates for
  different inherent scale values $s_0$ in the input data}
\label{sec-visualisation}

Figures~\ref{fig-lapl-intpt-det}--\ref{fig-princcurv-ridge-det}
provide visualisations of the conceptual steps involved, when defining
the scale estimates in Section~\ref{sec-meth-consist-scales}, that are then
quantitatively evaluated in Section~\ref{sec-results-scsel}.

In these figures, we
specifically show:
\begin{itemize}
\item
  {\em discretisations of the the idealised model signals,\/} that are used as
  input:
  \begin{itemize}
  \item
    Gaussian blobs, according to Equation~(\ref{eq-def-2D-gauss-blob-s0}),
    for Laplacian and determinant of the Hessian
    interest point detection,
  \item
    diffuse edges, according to Equation~(\ref{eq-def-ideal-blur-edge-s0}), for gradient magnitude based
    edge detection, and
  \item
    Gaussian ridges, according to Equation~(\ref{eq-def-ideal-ridge-s0}), for principal curvature based
    ridge detection,
  \end{itemize}
  for different choices of the scale values
  $\sigma_0 \in \{\tfrac{1}{2}, 1, 2 \}$,
\item
  {\em discrete approximations of the differential feature strength
  measures\/} at the scales $\hat{\sigma} = \sigma_0$, where the
  corresponding fully continuous differential feature strength
  measures ${\mathcal D}_{\norm} L$
  would assume their local extrema over scale $\sigma$ at the image
  center for each type of feature detector:
  \begin{itemize}
  \item
    the scale-normalised Laplacian operator $\nabla_{\norm}^2 L$,
    according to Equation~(\ref{eq-sc-norm-lapl}),
  \item
    the scale-normalised determinant of the Hessian $\det {\mathcal H}_{\norm} L$,
    according to Equation~(\ref{eq-sc-norm-det-hess}),
  \item
    the scale-normalised gradient magnitude $L_{v,\norm}$,
    according to Equation~(\ref{eq-sc-norm-grad-magn}), and
  \item
    the scale-normalized principal curvature based ridge strength
    measure $L_{pp,\norm}$, according to Equation~(\ref{eq-sc-norm-ridge-meas}),
  \end{itemize}
  for the different choices of the scale values
  $\sigma_0 \in \{\tfrac{1}{2}, 1, 2 \}$ for the input data, as well as
\item
  {\em graphs of the variability over scale $\sigma$ for discrete approximations of
  these scale-normalised measures of feature strength\/}, demonstrating
  that the selected scale values, obtained from the local extrema over
  scales, do reasonably well reflect the inherent characteristic scales
  $\sigma_0$ in the input data for the resulting scale selection methods.
\end{itemize}
For the visualisations in this appendix, we have in all these figures used: 
\begin{itemize}
\item
  the hybrid discretisation method, based on smoothing with the normalised
  sampled Gaussian kernel followed by central differences,
  according to Appendix~\ref{app-hybr-disc-normsamplgaussdiff}.
\end{itemize}
In the actual quantitative performance evaluation experiments
reported in Section~\ref{sec-meth-consist-scales},
corresponding computations are additionally performed for the four other
discrete derivative approximation methods based on:
\begin{itemize}
\item
  the sampled Gaussian derivative kernels,
  according to Appendix~\ref{app-sampl-gaussder-kern},
\item
  the integrated Gaussian derivative kernels,
  according to Appendix~\ref{app-int-gaussder-kern},
\item
  the discrete analogue of the Gaussian kernel combined with central
  difference operators, according to
  Appendix~\ref{app-disc-gauss-ders}, and
\item
   smoothing with the
   integrated Gaussian kernel followed by central difference
  operators, according to Appendix~\ref{app-hybr-disc-intgaussdiff}.
\end{itemize}
The graphs in Figures~\ref{fig-rel-sc-err-Lapl-edge-scsel-gauss-blob} and~\ref{fig-rel-sc-err-dethess-princcurv-scsel-gauss-blob} do additionally show the results of numerical scale
estimates for a much denser set of scale values
$\sigma_0 \in \{ \tfrac{1}{3}, 2 \}$.

\footnotesize
\bibliographystyle{abbrvnat}
\bibliography{bib/defs,bib/tlmac}

\end{document}